\documentclass{article}

\parskip \medskipamount
\usepackage[T1]{fontenc}
\usepackage{geometry}
\usepackage{titling}
\usepackage{setspace}
\usepackage{xcolor}
\usepackage{authblk}
\usepackage{amsmath,amssymb,amsthm,amsfonts,bm,bbm}
\DeclareMathOperator*{\esssup}{ess\,sup}
\DeclareMathOperator*{\essinf}{ess\,inf}

\DeclareMathOperator*{\Var}{Var}

\usepackage{mathrsfs,stmaryrd}
\usepackage{upgreek}
\usepackage{stmaryrd}
\SetSymbolFont{stmry}{bold}{U}{stmry}{m}{n}
\usepackage{indentfirst}
\usepackage{subdepth}
\usepackage{graphicx}
\usepackage{diagbox}
\usepackage{pdfpages}
\usepackage[latin1]{inputenc}
\usepackage{url}
\usepackage{float}
\usepackage{appendix}
\usepackage{etoolbox}
\usepackage{microtype}
\usepackage[centerlast]{caption}
\usepackage[plainpages=false,hyperfootnotes=false]{hyperref}
\hypersetup{colorlinks=true,urlcolor=blue,linkcolor=red,citecolor=blue}
\usepackage{cleveref}
\crefname{equation}{\hspace{-0.4em}}{\hspace{-0.4em}}
\crefrangeformat{equation}{(#3#1#4)--(#5#2#6)}
\usepackage{bookmark}
\usepackage{framed}
\usepackage{enumitem}
\usepackage{apacite}
\usepackage{booktabs}

\newtheorem{theorem}{Theorem}[section]
\newtheorem{lemma}[theorem]{Lemma}
\newtheorem{proposition}[theorem]{Proposition}
\newtheorem{remark}[theorem]{Remark}
\newtheorem{definition}[theorem]{Definition}

\renewenvironment{proof}{\noindent {\bf Proof.}}{\hfill $\Box$}
\allowdisplaybreaks[4]

\hypersetup{
  pdftitle={Time-consistent portfolio selection with monotone mean-variance preferences},
  pdfauthor={Y. Wang}
}

\title{Time-consistent portfolio selection \\ with monotone mean-variance preferences
\thanks{This work is supported by the National Natural Science Foundation of China [Grants 12401611, 12501663 and 12571520]; 
                                  Fundamental Research Funds for the Central Universities [Grant No. JBK202511002];
                                  Major Program of the Key Research Institute on Humanities and Social Science of China Ministry of Education [Grant 22JJD790091];
                                  the 111 Project [Grant B17050];
                              and CTBU [Grant 2355010].}}
\date{\vspace{-6ex}}
\author{
    Yike Wang\thanks{School of Finance, Chongqing Technology and Business University, Chongqing, China},
    Yusha Chen\thanks{School of Finance, Southwestern University of Finance and Economics, Chengdu, China.},
    Jingzhen Liu\thanks{China Institute of Actuarial Science, Central University of Finance and Economics, Beijing, China.}
}

\begin{document}
\maketitle

\begin{abstract}
We investigate time-inconsistent portfolio problems under a broader class of monotone mean-variance (MMV) preferences. 
Since the optimal strategies for MMV and mean-variance (MV) preferences coincide, the MMV optimal strategies at different initial times are necessarily time-inconsistent. 
To address this time-inconsistency, we consider Nash equilibrium controls of both open-loop and closed-loop types, and characterize them within a random parameter setting. 
The two control problems reduce to solving a flow of forward-backward stochastic differential equations and a system of extended Hamilton-Jacobi-Bellman equations, respectively. 
In particular, we derive semi-closed-form solutions for both types of equilibria under a deterministic parameter setting, and both solutions share the same representation, which is independent of the wealth state and the random path. 
We show that the investment amount under the MMV equilibrium exceeds that under the MV equilibrium, and the gap narrows over time. 
Furthermore, under a constant parameter setting, we find that the derived closed-loop Nash equilibrium control is a strong equilibrium strategy only when the interest rate is sufficiently large, 
whereas the derived open-loop Nash equilibrium control is necessarily a strong equilibrium strategy.
\end{abstract}

\vspace{3mm}

\vspace{3mm}

\section{Introduction}

The study of the mean-variance (MV) problem has been going on for decades (e.g., \cite{Markowitz-1952,Li-Ng-2000,Zhou-Li-2000}), becoming a cornerstone of stochastic control and mathematical finance. 
However, the conventional MV preference presents a lack of monotonicity, which implies that it may fail to align with state-wise dominance and, 
consequently, may not always yield rational choices, as pointed out by \cite{Maccheroni-Marinacci-Rustichini-Taboga-2009}. 
Intuitively, in many situations, the MV preference suggests choosing a strategy that is obviously worse, because the variance causes the penalty for positive deviations of a random return. 
To address this limitation, \cite{Maccheroni-Marinacci-Rustichini-Taboga-2009} built on their earlier work \cite{Maccheroni-Marinacci-Rustichini-2006} to propose a minor yet meaningful modification to the variational representation of the MV preference, 
giving rise to the monotone mean-variance (MMV) preference. 
Also, \cite{Cerny-2020} formulate the MMV preference by considering the monotone hull of the MV preference, and obtained the same result as the aforementioned variational modification method.
This refinement has since inspired a growing body of research on portfolio selection under MMV preferences, highlighting the distinctions and advantages over the classical MV formulation 
(e.g., \cite{Trybula-Zawisza-2019,Strub-Li-2020,Li-Guo-2021,Li-Liang-Pang-2022,Shen-Zou-2022,Hu-Shi-Xu-2023,Li-Guo-Tian-2023,Li-Liang-Pang-2023,Shi-Xu-2024,Du-Strub-2024,Li-Liang-Pang-2025,Cerny-Ruf-Schweizer-2025}).

Although the MMV preferences have been widely explored, \cite{Wang-Chen-Liu-Cui-2025} pointed out that the MMV preferences do not consistently withstand repeated testing.
More concretely, for the example given in \cite[Section 1]{Maccheroni-Marinacci-Rustichini-Taboga-2009}, the MMV preference suggests that the two choices are indistinguishable,
but every rational agent would make the same choice due to the presence of strict state-wise dominance. 
In other words, if numerous rational agents have encountered this selection problem or if a single rational agent is expected to address it repeatedly, the resulting choice will deviate from that in \cite{Maccheroni-Marinacci-Rustichini-Taboga-2009}. 
See details in \cite[Section 1]{Wang-Chen-Liu-Cui-2025}.
To handle this drawback, \cite{Wang-Chen-Liu-Cui-2025} introduce a class of strictly monotone mean-variance (SMMV) preferences based on a further minor modification on the MMV preference. 
More specifically, these preferences assign the following performance function to an uncertain prospect $X$:
\begin{equation*}
{V}_{ \theta, \zeta } ( X ) =
\inf_{ Y \in \mathbb{L}^{2}_{ \mathcal{F}_{T} } ( \Omega; \mathbb{R} ), Y \ge \zeta } 
  \bigg\{ \mathbb{E} [ X Y ] 
        - \inf_{ f \in \mathbb{L}^{2}_{ \mathcal{F}_{T} } ( \Omega; \mathbb{R} ) } 
          \Big\{ \mathbb{E} [ f Y ] - \Big( \mathbb{E} [ f ] - \frac{\theta }{2} \Var [f] \Big) \Big\} \bigg\}
\end{equation*}
where $\zeta $ is an actively chosen ``gradient'' to capture the strict monotonicity (see also \cite[Theorem 2.4]{Wang-Chen-Liu-Cui-2025}).
Conversely, the MMV preference is recovered as a special case by by setting $\zeta \equiv 0$.
In other words, the (time-consistent) MMV portfolio strategy is included in the present paper.

\begin{figure}[H]
  \centering
  \includegraphics[width=8cm]{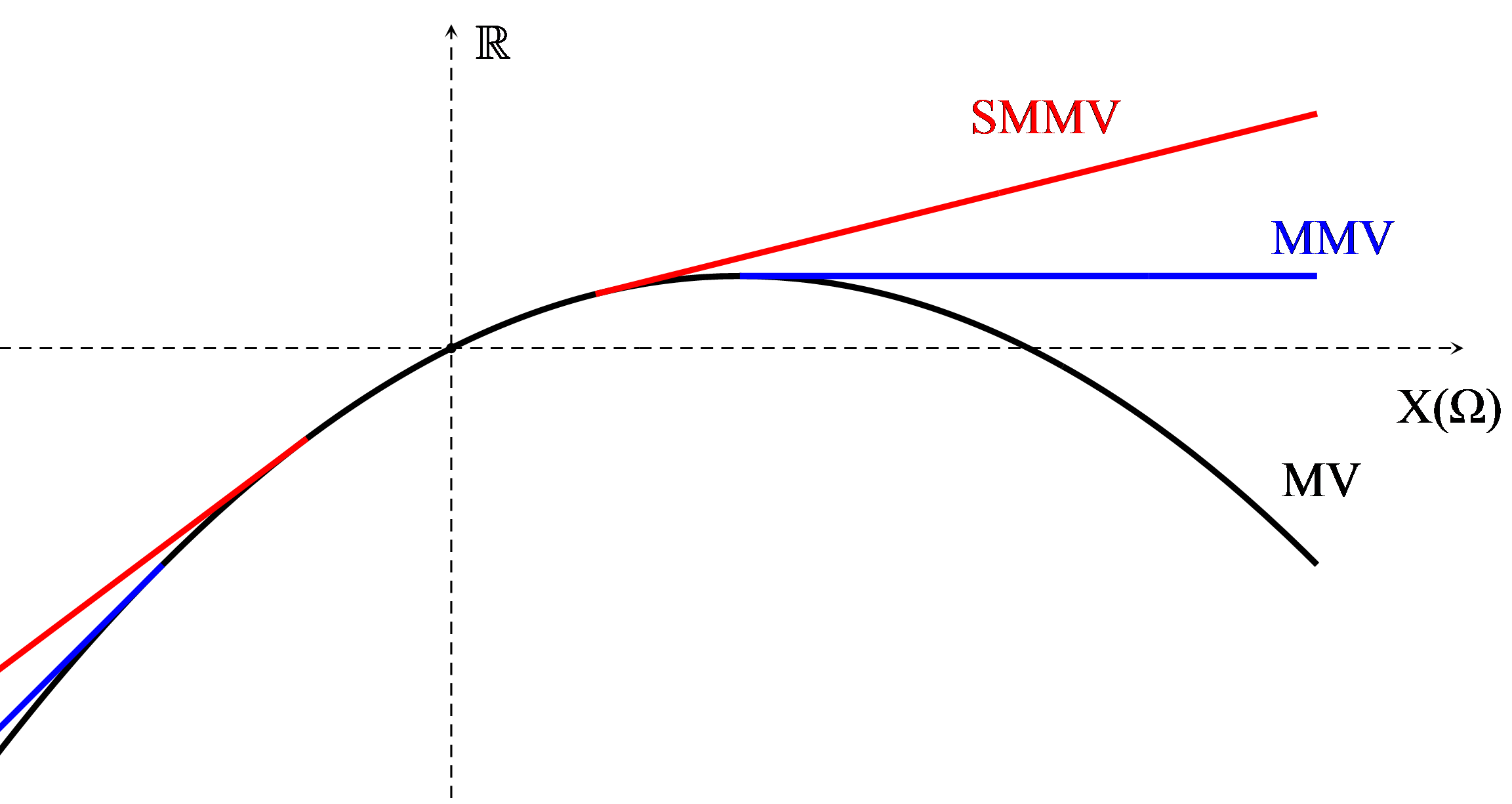}
  \caption{Example 2.8 given by \protect\cite{Wang-Chen-Liu-Cui-2025} illustrating the differences among MV, MMV, and SMMV Preferences.}
  \label{fig: illustrate MV MMV SMMV}
\end{figure}

An intuitive illustration for the differences among the MV, MMV and SMMV preferences is provided in \cite[Example 2.8, therein]{Wang-Chen-Liu-Cui-2025}, see also Figure \ref{fig: illustrate MV MMV SMMV}. Roughly speaking, the MMV preference merely flattens the right half of the MV parabola into a horizontal line, while the SMMV preference further tilts this line obliquely upward.

On the one hand, it was found that the optimal strategies under the MV, MMV, and SMMV preferences coincide in many cases (see \cite[Theorem 4.3]{Wang-Chen-Liu-Cui-2025}).
As a consequence, the MMV and SMMV optimal strategies for different initial epochs are time-inconsistent, inheriting this property from the well-established time-inconsistency of MV optimal strategies (see, e.g., \cite{Basak-Chabakauri-2010}).
On the other hand, like MV problems, the objective function for evaluating an SMMV preference is not standard for applying dynamic programming principle, and the optimal feedback control explicitly depends on initial time-state pairs.
In conclusion, an investor with MV/MMV/SMMV preferences is motivated to realize the optimality from present to future by changing her/his previous optimal plan,
and some incarnation of her/his future selves also has motivation to deviate from the current best plan.
In this paper, we focus on the time-inconsistency in the SMMV portfolio problems, following decades of research that can be traced back to \cite{Strotz-1955,Pollak-1968,Peleg-Yaari-1973,Goldman-1980}.

To tackle the time-inconsistency in dynamic maximization problems, the game-theoretical framework and spike variation method have gained increasing popularity in recent years. These approaches explore the ``optimality'' from a different perspective---seeking a subgame perfect Nash equilibrium in the game between the agent and the incarnation of her/his future selves---rather than aiming to realize a global maximum.
Any spike deviation from this equilibrium at any time instant will be worse off, so that no one has incentive to deviate from the equilibrium strategy.
Different formulations of spike perturbation leads to different type of equilibrium strategies, including the open-loop Nash equilibrium control (ONEC, in short) and the closed-loop Nash equilibrium control (CNEC, in short).
In particular, the equilibrium strategy is an ONEC (resp. a CNEC), if the spike perturbation is applied to stochastic process (resp. feedback function) of control.
For MV portfolio problems, the Nash equilibrium controls of both the open-loop and closed-loop types have been well studied, 
see, e.g., \cite{Hu-Jin-Zhou-2012,Hu-Jin-Zhou-2017,Huang-Li-Wang-2017} for the ONEC and \cite{Bjork-Murgoci-Zhou-2014,Bjork-Khapko-Murgoci-2017} for the CNEC.
One of the main objectives of this paper is to characterize both the ONEC and CNEC for SMMV portfolio problems.
Interested readers can set $\zeta \equiv 0$ in this paper to see the result for time-consistent portfolio problems with MMV preferences, which have not been studied elsewhere yet.

However, the mathematical definition of the aforementioned Nash equilibrium controls cannot ensure that the spike deviation is strictly worse off.
In other words, the equilibrium does not necessarily correspond to a maximum as mentioned in \cite[Remark 3.5]{Bjork-Khapko-Murgoci-2017}. 
In this case, the resulting equilibrium solution (i.e. the time-consistent strategy) may be questionable, as the agent could still have an incentive to deviate from it.
Once the equilibrium is merely a saddle point, some deviations from it are profitable in the discrete-time approximation with some time lattice.
In this situation, the corresponding Nash equilibrium control is not considered as a strong equilibrium strategy. 
\cite{He-Jiang-2021} studied the closed-loop strong equilibrium strategy for several concerned problems (including mean-variance and non-exponential discounting utility maximization).
In this paper, we investigate whether the derived Nash equilibrium controls in both types for the SMMV portfolio problem are strong equilibrium strategies, when all the parameters and the spike deviation are constant.

Summing up, in this paper, we formulate the time-consistent portfolio selection problem with SMMV preferences and seek Nash equilibrium controls.
The main contributions of this paper are as follows.

\emph{Firstly}, we reduce the ONEC seeking problem with random parameters to solving a flow of forward-backward stochastic differential equations (FBSDEs) under an equilibrium condition; see \Cref{thm: verification theorem for ONEC}.
It is worth mentioning that different from the classical time-inconsistent problems,
there is an implicit nonlinear law-dependent term $\lambda ( \cdot )$ (see \cref{eq: non-linear objective functional :eq} and \cref{eq: lambda functional :eq}) causing additional difficulties in the present paper.
Moreover precisely, we refer to \cite{Liang-Xia-Yuan-2025} and re-express the objective function as a nonlinear law-dependent preference;
however, we need to validate more properties due to the presence of $\lambda ( \cdot )$; see \Cref{lem: properties for ONEC} and \Cref{thm: asymptotic estimate for ONEC}.
Notably, our equilibrium condition \cref{eq: equilibrium condition for ONEC :eq} is formulated only by the diagonal processes, and does not rely on assuming the continuity of the BSDE solution, a condition imposed in, e.g., \cite{Liang-Xia-Yuan-2025}.
We shall validate both the sufficiency and necessity of the equilibrium condition, while \cite{Liang-Xia-Yuan-2025} only considered the sufficiency.
In addition, the variable separation method developed by \cite{Hu-Jin-Zhou-2017} is not applicable to deliver to the elegant equilibrium condition \cref{eq: equilibrium condition for ONEC :eq}
from the limit-form \cref{eq: primal equilibrium condition for ONEC :eq}.
To bridge this gap, we conduct more asymptotic estimates, an approach that distinguishes our perturbation analysis from those in the existing literature, e.g., \cite{Hu-Jin-Zhou-2017,Liang-Xia-Yuan-2025}.

\emph{Secondly}, we reduce the CNEC seeking problem to solving a system of extended Hamilton-Jacobi-Bellman (HJB) equations for CNECs; see \Cref{thm: verification theorem for CNEC}. 
Due to the presence of randomness of parameters, this system contains backward stochastic partial differential equations (BSPDEs).
We also tackle the new difficulty brought by the aforementioned law-dependent term $\lambda ( \cdot )$ in the perturbation analysis, 
and extend existing results represented in, e.g., \cite{Bjork-Khapko-Murgoci-2017,He-Jiang-2021,Wang-Liu-Bensoussan-Yiu-Wei-2025}, where only classical partial differential equations (PDEs) are employed.

\emph{Thirdly}, we derive the semi-closed-form expressions of an ONEC and a CNEC in the case with deterministic parameters; see \Cref{thm: path-independent ONEC,thm: state-independent CNEC}.
It is interesting to find that both types of Nash equilibrium controls deliver identical solutions, which are independent of random path and wealth state.
Moreover, the solution can be expressed as the conventional time-consistent MV portfolio strategy multiplied by a time-varying factor ${\beta }_{t} \ge 1$, which decreases over time.
This implies that although the optimal strategies under the MV, MMV, SMMV preferences coincide, their time-consistent strategies are not the same.
We also conduct a numerical experiment to illustrate the difference among the MV, MMV, SMMV time-consistent strategies.
It visually confirms that the SMMV equilibrium entails a larger investment amount than the MV equilibrium, and that the difference between them diminishes over time.

\emph{Finally}, given that the parameters are constants, we show that the aforementioned state-independent CNEC is not a strong equilibrium strategy unless the interest rate is sufficiently large,
while the aforementioned path-independent ONEC is necessarily a strong equilibrium strategy; see \Cref{thm: equilibrium value function and strong equilibria}.
This interesting finding implies that the equilibrium strategies (albeit expressed in the same form) exhibit different properties depending on the definition of equilibrium.
Notably, exploiting the derived explicit expression of SMMV preference, the validation is straightforward and more clear than using the higher-order conditions in \cite{He-Jiang-2021}.

The rest of this paper is organized as follows.
In \Cref{sec: Model and problem formulation}, we formulate our SMMV portfolio problems for both types of Nash equilibrium controls.
In \Cref{sec: Characterization for ONEC,sec: Characterization for CNEC}, we characterize the open-loop and the closed-loop Nash equilibrium control, respectively.
In \Cref{sec: closed-form solution}, we consider deterministic parameters and derive the semi-closed-form solution of the both two types.
In \Cref{sec: Numerical example}, we provide some numerical results with constant parameters.
In \Cref{sec: Concluding remark}, we make a brief concluding remark.

\section{Model and problem formulation}
\label{sec: Model and problem formulation}

\subsection{Notation}

Let $T \in ( 0, + \infty )$ be a fixed time horizon, and $\{ {W}_{t} \}_{ t \in [0,T ] }$ be a one-dimensional standard Brownian motion on $( \Omega, \mathcal{F}, \mathbb{F}, \mathbb{P} )$ that satisfying the usual condition, 
where $\mathbb{F} := \{ \mathcal{F}_{t} \}_{ t \in [ 0,T ] }$ is generated by $W$.
Denote the expectation operator by $\mathbb{E}$, and write $\mathbb{E}_{t} [ \cdot ] := \mathbb{E} [ \cdot | \mathcal{F}_{t} ]$, which is supposed to be modified to right-continuous.
To express the objective function as a law-dependent preference, we let $\mathcal{B} ( \mathbb{R} )$ be the Borel $\sigma $-algebra generated by all open intervals in the real line $\mathbb{R}$,
$\mathcal{P}$ be the Banach space of all finite signed measures on $( \mathbb{R} \times \mathbb{R}, \mathcal{B} ( \mathbb{R} ) \times \mathcal{B} ( \mathbb{R} ) )$,
and $\mathbb{P}^{t}_{X,Z} \in \mathcal{P}$ be the regular conditional law of random variables $( X,Z )$ given $\mathcal{F}_{t}$ under $\mathbb{P}$, so that $\mathbb{E}_{t} [ f ( X,Z ) ] = \int_{\Omega } f ( x,z ) \mathbb{P}^{t}_{X,Z} ( dx, dz )$.
For ease of reference, the functional spaces employed in this paper are listed below. Here $( n,m )$ are non-negative integers, $( t,p ) \in [ 0,T ) \times [ 1, \infty )$ and $\mathcal{R} \subseteq \mathbb{R}$.
\begin{itemize}
\item $\mathbb{L}^{p}_{ \mathcal{F}_{t} } ( \Omega; \mathcal{R} )$ denotes the set of all $\mathcal{F}_{t}$-measurable random variables $f: \Omega \to \mathcal{R}$ with $\mathbb{E} [ | f |^{p} ] < \infty $.
\item $\mathbb{L}^{2}_{\mathbb{F}} ( 0,T ; \mathbb{L}^{p} ( \Omega; \mathcal{R} ) )$ denotes the set of all $\mathbb{F}$-progressively measurable processes $f: [ 0,T ] \times \Omega \to \mathcal{R}$ 
      with $\mathbb{E} [ ( \int_{0}^{T} | f ( s, \cdot ) |^{2} ds )^{ \frac{p}{2} } ] < \infty $.
\item ${C}_{\mathbb{F}} ( 0,T; \mathbb{L}^{2} ( \Omega; \mathcal{R} ) )$ denotes the set of all $\mathbb{P}$-a.s. continuous processes $f \in \mathbb{L}^{2}_{\mathbb{F}} ( 0,T ; \mathbb{L}^{2} ( \Omega; \mathcal{R} ) )$
      with $\mathbb{E} [ \sup_{ s \in [ 0,T ] } | f ( s, \cdot ) |^{2} ] < \infty $.
\item ${C}_{\mathbb{F}} ( 0,T; \mathbb{L}^{2} ( \Omega; {C}^{n,m} ( \mathbb{R} \times \mathbb{R}; \mathcal{R} ) ) )$ denotes the set of all random fields $f: [ 0,T ] \times \Omega \times \mathbb{R} \times \mathbb{R} \to \mathcal{R}$ 
      such that $f ( \cdot, \cdot, x,y ) \in {C}_{\mathbb{F}} ( 0,T; \mathbb{L}^{2} ( \Omega; \mathcal{R} ) )$ 
                and $f ( t, \omega, x,y )$ is continuously differentiable in $x$ (resp. $y$) up to the $n$-th order (resp. $m$-th order) for every $t \in [ 0,T )$ and $\mathbb{P}$-a.e. $\omega \in \Omega $.
\item $\mathbb{L}^{2}_{\mathbb{F}} ( 0,T; \mathbb{L}^{2} ( \Omega; {C}^{n,m} ( \mathbb{R} \times \mathbb{R}; \mathcal{R} ) ) )$ denotes the set of all random fields $f: [ 0,T ] \times \Omega \times \mathbb{R} \times \mathbb{R} \to \mathcal{R}$ 
      such that $f ( \cdot, \cdot, x,y ) \in \mathbb{L}^{2}_{\mathbb{F}} ( 0,T ; \mathbb{L}^{2} ( \Omega; \mathcal{R} ) )$ 
                and $f ( t, \omega, x,y )$ is continuously differentiable in $x$ (resp. $y$) up to the $n$-th order (resp. $m$-th order) for every $t \in [ 0,T )$ and $\mathbb{P}$-a.e. $\omega \in \Omega $.
\end{itemize}
Similarly, we shall define ${C}_{\mathbb{F}} ( 0,T; \mathbb{L}^{2} ( \Omega; {C}^{2} ( \mathbb{R}; \mathbb{R} ) ) )$ 
and $\mathbb{L}^{2}_{\mathbb{F}} ( 0,T; \mathbb{L}^{2} ( \Omega; {C}^{1} ( \mathbb{R}; \mathbb{R} ) ) )$ for random fields $f: [ 0,T ] \times \Omega \times \mathbb{R} \to \mathbb{R}$.
For the sake of brevity, hereafter we suppress the statement of sample path $\omega $ and write ${f}_{t} = f ( t, \cdot )$, $f ( t,x ) = f ( t, \cdot, x )$, etc., unless otherwise mentioned.
In addition, we let $\Phi : \mathbb{R} \to ( 0,1 )$ denote the cumulated distribution function of standard normal distribution, $\Psi : ( 0, + \infty ) \to \mathbb{R}$ denote the inverse function of $f(x) = \int_{ - \infty }^{x} \Phi (z) dz$,
and $\Phi \circ \Psi ( \cdot ) := \Phi ( \Psi ( \cdot ) )$.

\subsection{Wealth dynamics}

Suppose that there is a bond and a stock in the financial market, the price dynamics of which evolve according to the following stochastic differential equations (SDEs):
\begin{equation*}
\left\{ \begin{aligned}
& d {B}_{t} = {B}_{t} {r}_{t} dt, && s.t. \quad {B}_{0} = {b}_{0} > 0 ~ \text{fixed a priori}; \\
& d {S}_{t} = {S}_{t} {r}_{t} dt + {S}_{t} {\sigma }_{t} ( d {W}_{t} + {\vartheta }_{t} dt ), && s.t. \quad {S}_{0} = {s}_{0} > 0 ~ \text{fixed a priori}.
\end{aligned} \right.
\end{equation*}
Here $( r, \sigma, \vartheta )$ denote the processes of the interest rate, volatility and risk premium, 
which are assumed to be $\mathbb{F}$-predictable, essentially bounded on $[ 0,T ] \times \Omega $ with $\essinf_{ [ 0,T ] \times \Omega } ( | \sigma | \wedge | \vartheta | ) > 0$.
Notably, the parameters may rely on the bond/stock price explicitly, e.g., ${r}_{t}$ could be a simplified form of $r ( t, \omega, B ( t, \omega ), S ( t, \omega ) )$.
Corresponding to the investment amount process $\pi \in \mathbb{L}^{2}_{\mathbb{F}} ( 0,T ; \mathbb{L}^{2} ( \Omega; \mathbb{R} ) )$ on the stock,
the wealth dynamics of a ``small'' investor is supposed to evolve as
\begin{equation}
\label{eq: controlled SDE with open-loop control :eq}
d {X}^{\pi }_{t} = \frac{ {\pi }_{t} }{ {S}_{t} } d {S}_{t} + \frac{ {X}^{\pi }_{t} - {\pi }_{t} }{ {B}_{t} } d {B}_{t}
                 = {X}^{\pi }_{t} {r}_{t} dt + {\pi }_{t} {\sigma }_{t} ( d {W}_{t} + {\vartheta }_{t} dt ),
\end{equation}  
with the initial amount ${X}^{\pi }_{0} = {x}_{0}$ fixing a priori.
Additionally, in order to formulate feedback control problems, we introduce the following controlled SDEs indexed by the initial pair $( t,x ) \in [ 0,T ) \times \mathbb{R}$:
\begin{equation}
\label{eq: controlled SDE with closed-loop control :eq}
d {X}^{ t,x, \Pi }_{s} = {X}^{ t,x, \Pi }_{s} {r}_{s} ds + \Pi ( s, {X}^{ t,x, \Pi }_{s} ) {\sigma }_{s} ( d {W}_{s} + {\vartheta }_{s} ds ), \quad s.t. \quad {X}^{ t,x, \Pi }_{t} = x.
\end{equation}  
Thus, $\Pi $ represents a feedback random field such that ${\pi }_{t} = \Pi ( t, {X}^{ 0, {x}_{0}, \Pi }_{t} )$.
Denote by $\mathcal{U}$ the set of all the feedback random fields $\Pi: [ 0,T ] \times \Omega \times \mathbb{R} \to \mathbb{R}$ such that 
$\{ \Pi ( t, {X}^{ 0, {x}_{0}, \Pi }_{t} ) \}_{ t \in [ 0,T ] } \in \mathbb{L}^{2}_{\mathbb{F}} ( 0,T ; \mathbb{L}^{2} ( \Omega; \mathbb{R} ) )$. 

\subsection{SMMV objective function}

Referring to \cite{Wang-Chen-Liu-Cui-2025}, 
for investors with SMMV preferences indexed by $( \theta, \zeta ) \in ( 0, + \infty ) \times \mathbb{L}^{2}_{ \mathcal{F}_{T} } ( \Omega; [ 0, 1 ] )$ with $\mathbb{E} [ \zeta ] < 1$,
the primal object is to maximize the SMMV performance functional:
\begin{align*}
{V}_{ \theta, \zeta } ( {X}^{\pi }_{T} )
& = \inf_{ Y \in \mathbb{L}^{2}_{ \mathcal{F}_{T} } ( \Omega; \mathbb{R} ), Y \ge \zeta } 
    \bigg\{ \mathbb{E} [ {X}^{\pi }_{T} Y ] 
          - \inf_{ f \in \mathbb{L}^{2}_{ \mathcal{F}_{T} } ( \Omega; \mathbb{R} ) } 
            \Big\{ \mathbb{E} [ f Y ] - \Big( \mathbb{E} [ f ] - \frac{\theta }{2} \mathbb{E} [ ( f - \mathbb{E} [ f ] )^{2} ] \Big) \Big\} \bigg\} \\
& = \theta \int_{ - \infty }^{ {\lambda }_{ {X}^{\pi }_{T}, \theta, \zeta } } s \mathbb{P} \Big( {X}^{\pi }_{T} + \frac{\zeta }{\theta } \le s \Big) ds
  + \mathbb{E} [ {X}^{\pi }_{T} \zeta ] + \frac{1}{ 2 \theta } \mathbb{E} [ {\zeta }^{2} ] - \frac{1}{ 2 \theta } \\
& = \mathbb{E} \bigg[ \theta \int_{ - \infty }^{ {\lambda }_{ {X}^{\pi }_{T}, \theta, \zeta } } s {1}_{\{ {X}^{\pi }_{T} + \frac{\zeta }{\theta } \le s \}} ds
                    + {X}^{\pi }_{T} \zeta + \frac{1}{ 2 \theta } ( {\zeta }^{2} - 1 ) \bigg],
\end{align*}
where the deterministic scalar ${\lambda }_{ {X}^{\pi }_{T}, \theta, \zeta } \in ( \essinf \{ {X}^{\pi }_{T} + \frac{\zeta }{\theta } \}, \mathbb{E} [ {X}^{\pi }_{T} ] + \frac{1}{\theta } ]$ is the unique solution of 
\begin{equation*}
1 = \mathbb{E} [ \zeta ] + \theta \int_{ - \infty }^{ {\lambda }_{ {X}^{\pi }_{T}, \theta, \zeta } } \mathbb{P} \Big( {X}^{\pi }_{T} + \frac{\zeta }{\theta } \le s \Big) ds
  \equiv \mathbb{E} \Big[ \zeta + \theta \Big( {\lambda }_{ {X}^{\pi }_{T}, \theta, \zeta } - {X}^{\pi }_{T} - \frac{\zeta }{\theta } \Big)_{+} \Big].
\end{equation*}
Notably, the existence and uniqueness of ${\lambda }_{ {X}^{\pi }_{T}, \theta, \zeta }$ arise from the continuity and strict monotonicity of $\int_{ - \infty }^{\lambda } \mathbb{P} ( {X}^{\pi }_{T} + \frac{\zeta }{\theta } \le s ) ds$
on $\lambda \in ( \essinf \{ {X}^{\pi }_{T} + \frac{\zeta }{\theta } \}, + \infty )$.
In fact, with the ``free'' parameter $\zeta $ for capturing the monotonicity, ${V}_{ \theta, \zeta } ( {X}^{\pi }_{T} )$ is a monotone modification of the MV performance functional
\begin{align*}
& \mathbb{E} [ {X}^{\pi }_{T} ] - \frac{\theta }{2} \mathbb{E} [ ( {X}^{\pi }_{T} - \mathbb{E} [ {X}^{\pi }_{T} ] )^{2} ] \\
& = \inf_{ Y \in \mathbb{L}^{2}_{ \mathcal{F}_{T} } ( \Omega; \mathbb{R} ) } 
    \bigg\{ \mathbb{E} [ {X}^{\pi }_{T} Y ] 
          - \inf_{ f \in \mathbb{L}^{2}_{ \mathcal{F}_{T} } ( \Omega; \mathbb{R} ) } 
            \Big\{ \mathbb{E} [ f Y ] - \Big( \mathbb{E} [ f ] - \frac{\theta }{2} \mathbb{E} [ ( f - \mathbb{E} [ f ] )^{2} ] \Big) \Big\} \bigg\},
\end{align*}
where the equality is due to Fenchel-Moreau theorem.

In this paper, we consider $\zeta \in \mathbb{L}^{2}_{ \mathcal{F}_{T} } ( \Omega; [ 0, 1 - \delta ] )$ for a sufficiently small $\delta > 0$, so that $\mathbb{E}_{t} [ \zeta ] < 1$.
Then, by mirroring the derivation for \cite[Theorem 2.9]{Wang-Chen-Liu-Cui-2025}, 
we introduce the SMMV objective function conditioned on $\mathcal{F}_{t}$ with $t \in [ 0,T )$ as the following:
\begin{align}
\notag
  J ( t, \pi )
& := \inf_{ Y \in \mathbb{L}^{2}_{ \mathcal{F}_{T} } ( \Omega; \mathbb{R} ), Y \ge \zeta } 
     \bigg\{ \mathbb{E}_{t} [ {X}^{\pi }_{T} Y ] 
           - \inf_{ f \in \mathbb{L}^{2}_{ \mathcal{F}_{T} } ( \Omega; \mathbb{R} ) } 
             \Big\{ \mathbb{E}_{t} [ f Y ] - \Big( \mathbb{E}_{t} [ f ] - \frac{\theta }{2} \mathbb{E}_{t} [ ( f - \mathbb{E}_{t} [ f ] )^{2} ] \Big) \Big\} \bigg\} \\
\label{eq: SMMV preference at epoch t :eq}
& = \mathbb{E}_{t} \bigg[ \theta \int_{ - \infty }^{ \lambda ( \mathbb{P}^{t}_{ {X}^{\pi }_{T}, \zeta } ) } s {1}_{\{ {X}^{\pi }_{T} + \frac{\zeta }{\theta } \le s \}} ds
                        + {X}^{\pi }_{T} \zeta + \frac{1}{ 2 \theta } ( {\zeta }^{2} - 1 ) \bigg],
\end{align}
where $\lambda ( \mathbb{P}^{t}_{ {X}^{\pi }_{T}, \zeta } ) \in ( \essinf \{ {X}^{\pi }_{T} + \frac{\zeta }{\theta } \}, \mathbb{E}_{t} [ {X}^{\pi }_{T} ] + \frac{1}{\theta } ]$ is $\mathcal{F}_{t}$-measurable and uniquely solves
\begin{equation*}
1 = \mathbb{E}_{t} \Big[ \zeta + \theta \Big( \lambda ( \mathbb{P}^{t}_{ {X}^{\pi }_{T}, \zeta } ) - {X}^{\pi }_{T} - \frac{\zeta }{\theta } \Big)_{+} \Big]
  \equiv \mathbb{E}_{t} [ \zeta ] 
       + \theta \int_{ - \infty }^{ \lambda ( \mathbb{P}^{t}_{ {X}^{\pi }_{T}, \zeta } ) } \mathbb{E}_{t} [ {1}_{\{ s \ge {X}^{\pi }_{T} + \frac{\zeta }{\theta } \}} ] ds.
\end{equation*}
Obviously, $\lambda ( \mathbb{P}^{T}_{ {X}^{\pi }_{T}, \zeta } ) = {X}^{\pi }_{T} + \frac{1}{\theta } \ge {X}^{\pi }_{T} + \frac{\zeta }{\theta }$, and hence
\begin{equation*}
J ( T, \pi ) := \theta \int_{ {X}^{\pi }_{T} + \frac{\zeta }{\theta } }^{ \lambda ( \mathbb{P}^{T}_{ {X}^{\pi }_{T}, \zeta } ) } s ds
             + {X}^{\pi }_{T} \zeta + \frac{1}{ 2 \theta } ( {\zeta }^{2} - 1 )
             = {X}^{\pi }_{T}.
\end{equation*}
Moreover, one can obtain the uniform integrability of $\{ \lambda ( \mathbb{P}^{t}_{ {X}^{\pi }_{T}, \zeta } ) \}_{ t \in [ 0,T ] }$, as the following.

\begin{lemma}\label{lem: moment estimate for lambda}
Let $\pi \in \mathbb{L}^{2}_{\mathbb{F}} ( 0,T ; \mathbb{L}^{2} ( \Omega; \mathbb{R} ) )$.
Then, $\sup_{ t \in [ 0,T ] } | \lambda ( \mathbb{P}^{t}_{ {X}^{\pi }_{T}, \zeta } ) | \in \mathbb{L}^{1}_{ \mathcal{F}_{T} } ( \Omega; \mathbb{R} )$.
\end{lemma}

\begin{proof}
See \Cref{pf-lem: moment estimate for lambda}.
\end{proof}

\subsection{Time-inconsistent control problem}

By mirroring the proof of \cite[Theorem 4.3]{Wang-Chen-Liu-Cui-2025}, one can conclude that
the unique maximizer for the SMMV objective function \cref{eq: SMMV preference at epoch t :eq} is 
the MV optimal portfolio strategy for $\mathbb{E}_{t} [ {X}^{\pi }_{T} ] - \frac{\theta }{2} \mathbb{E}_{t} [ ( {X}^{\pi }_{T} - \mathbb{E}_{t} [ {X}^{\pi }_{T} ] )^{2} ]$,
provided that 
\begin{equation*}
\zeta \le \frac{ \exp \big( - \int_{t}^{T} {r}_{v} dv - \int_{t}^{T} {\theta }_{v} d {W}_{v} - \frac{1}{2} \int_{t}^{T} | {\theta }_{v} |^{2} dv \big) }
                { \mathbb{E}_{t} \big[ \exp \big( - \int_{t}^{T} {r}_{v} dv - \int_{t}^{T} {\theta }_{v} d {W}_{v} - \frac{1}{2} \int_{t}^{T} | {\theta }_{v} |^{2} dv \big) \big] }.
\end{equation*}
Consequently, the optimal SMMV portfolio strategies for different initial epochs $t$ are time-inconsistent, inherit this property from the well-established time-inconsistency of MV optimal strategies (see, e.g., \cite{Basak-Chabakauri-2010}).
To tackle the time-inconsistency, we formulate the control problem as a game among the present and future incarnations of the investor.
Referring to \cite{Hu-Jin-Zhou-2012,Hu-Jin-Zhou-2017} for open-loop equilibria and \cite{Bjork-Khapko-Murgoci-2017} for closed-loop equilibria, we introduce the following definitions.

\begin{definition}\label{def: open-loop Nash equilibrium control}
$\bar{\pi } \in \mathbb{L}^{2}_{\mathbb{F}} ( 0,T ; \mathbb{L}^{2} ( \Omega; \mathbb{R} ) )$ is an open-loop Nash equilibrium control (ONEC), if
\begin{equation}
\label{eq:open-loop Nash equilibrium control :eq}
0 \le \liminf_{\varepsilon \downarrow 0}
      \frac{1}{\varepsilon } \big( J ( t, \bar{\pi } ) - J ( t, \bar{\pi }^{ t, \varepsilon, \xi } ) \big), \quad \mathbb{P}-a.s., ~ a.e. ~ t \in [ 0,T ),
\end{equation}
for any $\xi \in \mathbb{L}^{2}_{ \mathcal{F}_{t} } ( \Omega; \mathbb{R} )$,
where $\bar{\pi }^{ t, \varepsilon, \xi }$ is a spike variation of $\bar{\pi }$ given by $\bar{\pi }^{ t, \varepsilon, \xi }_{s} = \bar{\pi }_{s} + \xi {1}_{\{ s \in [ t, t + \varepsilon ) \}}$.
\end{definition}

\begin{definition}\label{def: closed-loop Nash equilibrium control}
$\bar{\Pi } \in \mathcal{U}$ is a closed-loop Nash equilibrium control (CNEC), if
\begin{equation}
\label{eq:closed-loop Nash equilibrium control :eq}
0 \le \liminf_{\varepsilon \downarrow 0}
      \frac{1}{\varepsilon } \big( \bar{J} ( t, {X}^{ 0, {x}_{0}, \bar{\Pi } }_{t}, \bar{\Pi } ) - \bar{J} ( t, {X}^{ 0, {x}_{0}, \bar{\Pi } }_{t}, \bar{\Pi }^{ t, \varepsilon, \xi } ) \big), \quad \mathbb{P}-a.s., ~ a.e. ~ t \in [ 0,T ),
\end{equation}
for any $\xi \in \mathbb{L}^{2}_{ \mathcal{F}_{t} } ( \Omega; \mathbb{R} )$ satisfying $\bar{\Pi }^{ t, \varepsilon, \xi } \in \mathcal{U}$,
where 
\begin{equation}
\label{eq: SMMV preference with closed-loop control :eq}
   \bar{J} ( t, x, \Pi )
:= \mathbb{E}_{t} \bigg[ \theta \int_{ - \infty }^{ \lambda ( \mathbb{P}^{t}_{ {X}^{ t,x, \Pi }_{T}, \zeta } ) } s {1}_{\{ {X}^{ t,x, \Pi }_{T} + \frac{\zeta }{\theta } \le s \}} ds 
                       + {X}^{ t,x, \Pi }_{T} \zeta + \frac{1}{ 2 \theta } ( {\zeta }^{2} - 1 ) \bigg]
\end{equation}
with the $\mathcal{F}_{t}$-measurable random variable $\lambda ( \mathbb{P}^{t}_{ {X}^{ t,x, \Pi }_{T}, \zeta } )$ being the unique solution of 
\begin{equation*}
1 = \mathbb{E}_{t} \Big[ \zeta + \theta \Big( \lambda ( \mathbb{P}^{t}_{ {X}^{ t,x, \Pi }_{T}, \zeta } ) - {X}^{ t,x, \Pi }_{T} - \frac{\zeta }{\theta } \Big)_{+} \Big],
\end{equation*}
and $\bar{\Pi }^{ t, \varepsilon, \xi }$ is a spike variation of $\bar{\Pi }$ given by
$\bar{\Pi }^{ t, \varepsilon, \xi } ( s,y ) = \bar{\Pi } ( s,y ) {1}_{\{ s \notin [ t, t + \varepsilon ) \}} + \xi {1}_{\{ s \in [ t, t + \varepsilon ) \}}$.
\end{definition}

In short, we aim to characterizing the aforementioned ONEC and CNEC in general, and find their explicit expressions in some specific cases.

\section{Characterization for ONEC}
\label{sec: Characterization for ONEC}

In this section, inspired by \cite{Liang-Xia-Yuan-2025},
we re-express the SMMV preferences as nonlinear law-dependence preferences, and then reduce our control problem for ONEC to solving a flow of FBSDEs with a simplified equilibrium condition.
We firstly rewrite \cref{eq: SMMV preference at epoch t :eq} as 
\begin{align}
\notag
  J ( t, \pi ) 
& = \lambda ( \mathbb{P}^{t}_{ {X}^{\pi }_{T}, \zeta } ) \theta \int_{ - \infty }^{ \lambda ( \mathbb{P}^{t}_{ {X}^{\pi }_{T}, \zeta } ) } \mathbb{E}_{t} [ {1}_{\{ {X}^{\pi }_{T} + \frac{\zeta }{\theta } \le s \}} ] ds \\
\notag
& \quad - \theta \mathbb{E}_{t} \bigg[ \int_{ - \infty }^{ \lambda ( \mathbb{P}^{t}_{ {X}^{\pi }_{T}, \zeta } ) } 
                                       ( \lambda ( \mathbb{P}^{t}_{ {X}^{\pi }_{T}, \zeta } ) - s ) {1}_{\{ {X}^{\pi }_{T} + \frac{\zeta }{\theta } \le s \}} ds \bigg]
        + \mathbb{E}_{t} \Big[ {X}^{\pi }_{T} \zeta + \frac{1}{ 2 \theta } ( {\zeta }^{2} - 1 ) \Big] \\
\label{eq: re-expression of SMMV preference :eq}
& = \lambda ( \mathbb{P}^{t}_{ {X}^{\pi }_{T}, \zeta } ) ( 1 - \mathbb{E}_{t} [ \zeta ] )
  - \frac{\theta }{2} \mathbb{E}_{t} \bigg[ \Big| \Big( \lambda ( \mathbb{P}^{t}_{ {X}^{\pi }_{T}, \zeta } ) - {X}^{\pi }_{T} - \frac{\zeta }{\theta } \Big)_{+} \Big|^{2} \bigg]
  + \mathbb{E}_{t} \Big[ {X}^{\pi }_{T} \zeta + \frac{1}{ 2 \theta } ( {\zeta }^{2} - 1 ) \Big],
\end{align}
and denote by $\mathcal{P}_{0}$ the set of all regular conditional laws $p \in \mathcal{P}$ such that 
$\int_{ \mathbb{R} \times \mathbb{R} } ( {x}^{2} + {z}^{2} ) p ( dx, dz ) < \infty $ and $\int_{ \mathbb{R} \times \mathbb{R} } z p ( dx, dz ) < 1$.
Then, let us introduce the following functional defined on the convex set $\mathcal{P}_{0}$:
\begin{equation}
\label{eq: non-linear objective functional :eq}
g (p) = \int_{ \mathbb{R} \times \mathbb{R} } \bigg( \lambda (p) ( 1 - z ) - \frac{\theta }{2} \Big| \Big( \lambda (p) - x - \frac{z}{\theta } \Big)_{+} \Big|^{2} + x z + \frac{1}{ 2 \theta } ( {z}^{2} - 1 ) \bigg) p ( dx, dz ),
\end{equation}
where the deterministic scalar $\lambda (p)$ is the solution of 
\begin{equation}
\label{eq: lambda functional :eq}
1 = \int_{ \mathbb{R} \times \mathbb{R} } \bigg( z + \theta \Big( \lambda (p) - x -\frac{z}{\theta } \Big)_{+} \bigg) p ( dx, dz ),
\end{equation}
so that $J ( t, \pi ) = g ( \mathbb{P}^{t}_{ {X}^{\pi }_{T}, \zeta } )$.
Notably, for any $p \in \mathcal{P}_{0}$, the existence and uniqueness of $\lambda (p)$,
as well as $\int_{\{ x + \frac{z}{\theta } < \lambda (p) \}} p ( dx, dz ) > 0$, can be seen from the following equivalent expression of \cref{eq: lambda functional :eq}:
\begin{equation}
\label{eq: lambda functional equivalent expression :eq}
0 < 1 - \int_{ \mathbb{R} \times \mathbb{R} } z p ( dx, dz ) = \theta \int_{ - \infty }^{ \lambda (p) } ds \int_{\{ x + \frac{z}{\theta } \le s \}} p ( dx, dz ).
\end{equation}

\begin{lemma}\label{lem: properties for ONEC}
For the functions $( g, \lambda )$ defined by \cref{eq: non-linear objective functional :eq} and \cref{eq: lambda functional :eq}, the following properties holds.
\begin{itemize}
\item For any $p,q \in \mathcal{P}_{0}$,
      \begin{equation*}
        \lim_{ \varepsilon \downarrow 0 } \frac{1}{\varepsilon } \Big( g \big( p + \varepsilon (q-p) \big) - g (p) \Big) 
      = d g ( p, q-p ) := \int_{ \mathbb{R} \times \mathbb{R} } \nabla g ( p,x,z ) \big( q ( dx, dz ) - p ( dx, dz ) \big),
      \end{equation*}
      where $\nabla g : \mathcal{P}_{0} \times \mathbb{R} \times \mathbb{R} \to \mathbb{R}$ is given by
      \begin{equation*}
      \nabla g ( p,x,z ) := \lambda (p) ( 1 - z ) - \frac{\theta }{2} \Big| \Big( \lambda (p) - x - \frac{z}{\theta } \Big)_{+} \Big|^{2} + x z + \frac{1}{ 2 \theta } ( {z}^{2} - 1 ).
      \end{equation*}
      Thus, $g: \mathcal{P}_{0} \to \mathbb{R}$ is continuously G\^ateaux differentiable, 
      and $d g ( p, \cdot )$ as the G\^ateaux differential of $g$ at $p$ is a continuous linear functional determined by $\nabla g ( p, \cdot, \cdot )$,
      which satisfies 
      \begin{equation*}
      \int_{ \mathbb{R} \times \mathbb{R} } | \nabla g ( p,x,z ) | q ( dx, dz ) < \infty, \quad \forall p,q \in \mathcal{P}_{0}.
      \end{equation*}
\item For any $p,q \in \mathcal{P}_{0}$, $| g(q) - g(p) - d g ( p, q-p ) | \le \frac{1}{2} | \lambda (q) - \lambda (p) |^{2}$.
\item For any $( t, \chi ) \in [ 0,T ) \times \mathbb{L}^{2}_{ \mathcal{F}_{T} } ( \Omega; \mathbb{R} )$,
      \begin{equation*}
          \frac{ \mathbb{E}_{t} [ \chi {1}_{ {R}_{1} } ] }{ \mathbb{E}_{t} [ {1}_{ {R}_{1} } ] }
      \le \lambda ( \mathbb{P}^{t}_{ {X}^{\pi }_{T} + \chi, \zeta } ) - \lambda ( \mathbb{P}^{t}_{ {X}^{\pi }_{T}, \zeta } ) 
      \le \frac{ \mathbb{E}_{t} [ \chi {1}_{ {R}_{0} } ] }{ \mathbb{E}_{t} [ {1}_{ {R}_{0} } ] },
      \end{equation*}
      where ${R}_{0} = \{ {X}^{\pi }_{T} + \frac{\zeta }{\theta } \le \lambda ( \mathbb{P}^{t}_{ {X}^{\pi }_{T}, \zeta } ) \}$ and ${R}_{1} = \{ {X}^{\pi }_{T} + \chi + \frac{\zeta }{\theta } \le \lambda ( \mathbb{P}^{t}_{ {X}^{\pi }_{T} + \chi, \zeta } ) \}$.
\end{itemize}
\end{lemma}

\begin{proof}
See \Cref{pf-lem: properties for ONEC}.
\end{proof}

So far, we have shown that $g$ defined by \cref{eq: non-linear objective functional :eq} and \cref{eq: lambda functional :eq} satisfies \cite[Assumption 3.1]{Liang-Xia-Yuan-2025}, 
so that \cite[Theorem 3.2]{Liang-Xia-Yuan-2025} can be applied for deriving a sufficient condition for ONEC.
More than that, we shall provide the necessity of the derived sufficient condition, serving for the uniqueness of equilibrium. 
In addition, we shall remove the strong assumption imposed in \cite[Theorem 3.2]{Liang-Xia-Yuan-2025}, which assumes that $\mathcal{Y}^{t}_{s}$ as a part of the BSDE solution in the present paper is right-continuous in $s$.
To achieve these goal, we establish the following \Cref{thm: asymptotic estimate for ONEC} for expanding $J ( t, \bar{\pi } ) - J ( t, \bar{\pi }^{ t, \varepsilon, \xi } )$
and isolate the main result of this section into the upcoming \Cref{thm: verification theorem for ONEC}.

\begin{theorem}\label{thm: asymptotic estimate for ONEC}
For any $\bar{\pi } \in \mathbb{L}^{2}_{\mathbb{F}} ( 0,T ; \mathbb{L}^{2} ( \Omega; \mathbb{R} ) )$, $t \in [ 0,T )$, $\varepsilon \in ( 0, T-t ]$ and $\xi \in \mathbb{L}^{2}_{ \mathcal{F}_{t} } ( \Omega; \mathbb{R} )$, 
the following asymptotic estimate holds:
\begin{align}
\notag
  J ( t, \bar{\pi }^{ t, \varepsilon, \xi } ) - J ( t, \bar{\pi } )
& = \xi \int_{t}^{ t + \varepsilon } \mathbb{E}_{t} [ ( {Y}^{t}_{s} {\vartheta }_{s} + \mathcal{Y}^{t}_{s} ) {\sigma }_{s} ] ds \\
\label{eq: asymptotic estimate for ONEC :eq}
& \quad 
  - \frac{\theta }{2} {\xi }^{2} \mathbb{E}_{t} \bigg[ \bigg( \int_{t}^{ t + \varepsilon } {\sigma }_{s} d {W}_{s} \bigg)^{2} 
                                                       {e}^{ 2 \int_{t}^{T} {r}_{v} dv } {1}_{\{ {X}^{ \bar{\pi } }_{T} + \frac{\zeta }{\theta } \le \lambda ( \mathbb{P}^{t}_{ {X}^{ \bar{\pi } }_{T}, \zeta } ) \}} \bigg] 
  + o ( \varepsilon ),
\end{align}
where $( {Y}^{t}, \mathcal{Y}^{t} ) = ( Y ( \lambda ( \mathbb{P}^{t}_{ {X}^{ \bar{\pi } }_{T}, \zeta } ) ), \mathcal{Y} ( \lambda ( \mathbb{P}^{t}_{ {X}^{ \bar{\pi } }_{T}, \zeta } ) ) )$ 
with $( Y ( \lambda ), \mathcal{Y} ( \lambda ) ) \in {C}_{\mathbb{F}} ( 0,T; \mathbb{L}^{2} ( \Omega; \mathbb{R} ) ) \times \mathbb{L}^{2}_{\mathbb{F}} ( 0,T; \mathbb{L}^{2} ( \Omega; \mathbb{R} ) )$ given by the unique solution of the linear BSDE:
\begin{equation}
\label{eq: BSDE :eq}
{Y}_{s} ( \lambda ) = \zeta + \theta \Big( \lambda - {X}^{ \bar{\pi } }_{T} - \frac{\zeta }{\theta } \Big)_{+} 
            + \int_{s}^{T} {Y}_{v} ( \lambda ) {r}_{v} dv - \int_{s}^{T} \mathcal{Y}_{v} ( \lambda ) d {W}_{v}.
\end{equation}
\end{theorem}

\begin{proof}
See \Cref{pf-thm: asymptotic estimate for ONEC}.
\end{proof}

In view of the arbitrariness of $\xi $, one can conclude that $\limsup_{ \varepsilon \downarrow 0 } \frac{1}{\varepsilon } ( J ( t, \bar{\pi }^{ t, \varepsilon, \xi } ) - J ( t, \bar{\pi } ) ) \le 0$, i.e. $\bar{\pi }$ is an ONEC,
if and only if
\begin{equation}
\label{eq: primal equilibrium condition for ONEC :eq}
\lim_{ \varepsilon \downarrow 0 } \frac{1}{\varepsilon } \int_{t}^{ t + \varepsilon } \mathbb{E}_{t} [ ( {Y}^{t}_{s} {\vartheta }_{s} + \mathcal{Y}^{t}_{s} ) {\sigma }_{s} ] ds = 0, \quad \mathbb{P}-a.s., ~ a.e. ~ t \in [ 0,T ].
\end{equation}
However, the limit-form equilibrium condition \cref{eq: primal equilibrium condition for ONEC :eq} renders it of little practical value,
especially when we have not verified the right-continuity of $( \mathcal{Y}^{t}_{s}, {\vartheta }_{s}, {\sigma }_{s} )$ in $s$.
To address this concern, we reduce \cref{eq: primal equilibrium condition for ONEC :eq} to a brief form that only includes the diagonal processes $\{ ( {Y}^{s}_{s}, \mathcal{Y}^{s}_{s} ) \}_{ s \in [ 0,T ] }$ and the parameter ${\vartheta }_{t}$.
Notably, when defining the diagonal processes based on $( {Y}_{s} ( \lambda ), \mathcal{Y}_{s} ( \lambda ) )$, 
we first fix $s \in [ 0,T ]$ and then plug $\lambda = \lambda ( \mathbb{P}^{s}_{ {X}^{ \bar{\pi } }_{T}, \zeta } )$ into the result.

\begin{theorem}\label{thm: verification theorem for ONEC}
Suppose that $\bar{\pi } \in \mathbb{L}^{2}_{\mathbb{F}} ( 0,T ; \mathbb{L}^{2} ( \Omega; \mathbb{R} ) )$.
Then, $\bar{\pi }$ is an ONEC, if and only if 
\begin{equation}
\label{eq: equilibrium condition for ONEC :eq}
{Y}^{t}_{t} {\vartheta }_{t} + \mathcal{Y}^{t}_{t} = 0, \quad \mathbb{P}-a.s., ~ a.e. ~ t \in [ 0,T ].
\end{equation}
\end{theorem}

\begin{proof}
See \Cref{pf-thm: verification theorem for ONEC}.
\end{proof}

\begin{remark}
Since the terminal value ${Y}^{t}_{T} = \zeta + \theta ( \lambda ( \mathbb{P}^{t}_{ {X}^{ \bar{\pi } }_{T}, \zeta } ) - {X}^{ \bar{\pi } }_{T} - \frac{\zeta }{\theta } )_{+}$ 
cannot be fully separated as in \cite{Hu-Jin-Zhou-2012,Hu-Jin-Zhou-2017},
the constructive method developed therein cannot be applied to derive the structure of ${Y}^{t}_{s} {\vartheta }_{s} + \mathcal{Y}^{t}_{s}$.
In this sense, it is innovative in the proof of \Cref{thm: verification theorem for ONEC} (see \Cref{pf-thm: verification theorem for ONEC}) to separate the variables $( \lambda, s )$ based on the martingale representations, e.g.,
\begin{equation}
\label{eq: martingale representation for indicator :eq}
P ( t,y ) = P ( 0,y ) + \int_{0}^{t} \eta ( s,y ) d {W}_{s}, \quad where \quad P ( t,y ) = \mathbb{E}_{t} [ {1}_{\{ y \ge {X}^{ \bar{\pi } }_{T} + \frac{\zeta }{\theta } \}} ].
\end{equation}
In addition, limited by the lack of adequate integrability of the diagonal processes $\{ ( {Y}^{t}_{t}, \mathcal{Y}^{t}_{t} ) \}_{ t \in [ 0,T ] }$, which is caused by the presence of $\lambda ( \mathbb{P}^{t}_{ {X}^{ \bar{\pi } }_{T}, \zeta } )$,
the stochastic Lebesgue differentiation theorem formulated by \cite[Lemma 3.4]{Hu-Jin-Zhou-2017} cannot be used to derive the equivalence between \cref{eq: primal equilibrium condition for ONEC :eq} and \cref{eq: equilibrium condition for ONEC :eq}.
\end{remark}

Although we have reduced our ONEC seeking problem to solving a flow of FBSDEs given by 
\cref{eq: controlled SDE with open-loop control :eq} (with $\pi = \bar{\pi }$ therein), \cref{eq: BSDE :eq} and the so-called equilibrium condition \cref{eq: equilibrium condition for ONEC :eq}, 
it is still difficult to arrive at a clear expression of some ONEC in general.
One could further reduce our problem to solving a BSPDE, referring to the following \Cref{lem: property of ONEC}.

\begin{theorem}\label{lem: property of ONEC}
Assume that $r$ is deterministic, and $\bar{\pi }$ is an ONEC such that $( P, \eta )$ given by \cref{eq: martingale representation for indicator :eq} is of
${C}_{\mathbb{F}} ( 0,T; \mathbb{L}^{2} ( \Omega; {C}^{2} ( \mathbb{R}; \mathbb{R} ) ) ) \times \mathbb{L}^{2}_{\mathbb{F}} ( 0,T; \mathbb{L}^{2} ( \Omega; {C}^{1} ( \mathbb{R}; \mathbb{R} ) ) )$.
Then, the random fields 
\begin{equation}
\label{eq: rho-varrho :eq}
\rho ( t,y ) := P \big( t, y + \lambda ( \mathbb{P}^{t}_{ {X}^{ \bar{\pi } }_{T}, \zeta } ) \big) \quad and \quad
\varrho ( t,y ) := \eta \big( t, y + \lambda ( \mathbb{P}^{t}_{ {X}^{ \bar{\pi } }_{T}, \zeta } ) \big) + {\rho }_{y} ( t,y ) \frac{ {\vartheta }_{t} }{ \theta \rho ( t,0 ) }
\end{equation}
satisfy the following BSPDE on $[ 0,T ] \times \Omega \times \mathbb{R}$:
\begin{equation}
\label{eq: BSPDE of rho :eq}
\left\{ \begin{aligned}               
- d \rho ( t,y ) & = \frac{1}{2} {\rho }_{yy} ( t,y ) \bigg| \frac{ {\vartheta }_{t} }{ \theta \rho ( t,0 ) } \bigg|^{2} dt
                   - {\rho }_{y} ( t,y ) \bigg( \frac{ {\rho }_{y} ( t,0 ) }{ 2 \rho ( t,0 ) } \bigg| \frac{ {\vartheta }_{t} }{ \theta \rho ( t,0 ) } \bigg|^{2}
                                              - \frac{ \varrho ( t,0 ) }{ \rho ( t,0 ) } \frac{ {\vartheta }_{t} }{ \theta \rho ( t,0 ) } \bigg) dt \\
                 & \quad - {\varrho }_{y} ( t,y ) \frac{ {\vartheta }_{t} }{ \theta \rho ( t,0 ) } dt
                         - \varrho ( t,y ) d {W}_{t}, \\
    \rho ( T,y ) & = {1}_{\{ \zeta \le \theta y + 1 \}}.
\end{aligned} \right.
\end{equation}
Conversely, if \cref{eq: BSPDE of rho :eq} admits a solution
$( \rho, \varrho ) \in {C}_{\mathbb{F}} ( 0,T; \mathbb{L}^{2} ( \Omega; {C}^{2} ( \mathbb{R}; \mathbb{R} ) ) ) \times \mathbb{L}^{2}_{\mathbb{F}} ( 0,T; \mathbb{L}^{2} ( \Omega; {C}^{1} ( \mathbb{R}; \mathbb{R} ) ) )$ with 
\begin{align*}
& \lim_{ y \downarrow - \infty } \esssup_{ [ 0,T ] \times \Omega } \big( | \rho ( \cdot, y ) | + | {\rho }_{y} ( \cdot, y ) | + | \varrho ( \cdot, y ) | \big) = 0, \\
& \int_{ - \infty }^{\lambda } \Big( \int_{0}^{T} | \varrho ( t,y ) |^{2} dt \Big)^{ \frac{1}{2} } dy < \infty, \quad \mathbb{P}-a.s., \quad \forall \lambda \in \mathbb{R}, \\
& \mathbb{E} \bigg[ \int_{0}^{T} \bigg| \frac{ {\vartheta }_{t} }{ \theta \rho ( t,0 ) } \bigg|^{2} dt 
                  + \bigg( \int_{0}^{T} \bigg| \frac{ {\rho }_{y} ( t,0 ) }{ 2 \rho ( t,0 ) } \Big| \frac{ {\vartheta }_{t} }{ \theta \rho ( t,0 ) } \Big|^{2}
                                             - \frac{ \varrho ( t,0 ) }{ \rho ( t,0 ) } \frac{ {\vartheta }_{t} }{ \theta \rho ( t,0 ) } \bigg| dt \bigg)^{2} \bigg] < \infty.
\end{align*}
then the unique hedging portfolio $\bar{\pi } \in \mathbb{L}^{2}_{\mathbb{F}} ( 0,T ; \mathbb{L}^{2} ( \Omega; \mathbb{R} ) )$ for 
\begin{equation*}
{X}^{ \bar{\pi } }_{T} = c + \int_{0}^{T} \frac{ {\vartheta }_{t} }{ \theta \rho ( t,0 ) } d {W}_{t} 
                           + \int_{0}^{T} \bigg( \frac{ {\rho }_{y} ( t,0 ) }{ 2 \rho ( t,0 ) } \bigg| \frac{ {\vartheta }_{t} }{ \theta \rho ( t,0 ) } \bigg|^{2}
                                               - \frac{ \varrho ( t,0 ) }{ \rho ( t,0 ) } \frac{ {\vartheta }_{t} }{ \theta \rho ( t,0 ) } \bigg) dt, 
\end{equation*}
with some ${x}_{0}$-dependent constant $c \in \mathbb{R}$, is an ONEC.
\end{theorem}

\begin{proof}
See \Cref{pf-lem: property of ONEC}.
Notably, the constant $c$ is uniquely determined by ${x}_{0} = \mathbb{E} [ \frac{ d \mathbb{Q} }{ d \mathbb{P} } {X}^{ \bar{\pi } }_{T} {e}^{ - \int_{0}^{T} {r}_{v} dv } ]$, 
where $\mathbb{Q}$ is the risk-neutral measure given by $\frac{ d \mathbb{Q} }{ d \mathbb{P} } |_{ \mathcal{F}_{T} } = {e}^{ - \int_{0}^{T} {\vartheta }_{v} d {W}_{v} - \frac{1}{2} \int_{0}^{T} | {\vartheta }_{v} |^{2} dv }$.
\end{proof}

In fact, \Cref{lem: property of ONEC} also serves for the uniqueness of open-loop equilibrium. 
More precisely, if the BSPDE \cref{eq: BSPDE of rho :eq} admits a unique solution, then there exists at most one ONEC.
However, due to its terminal condition $\rho ( T,y ) = {1}_{\{ \zeta \le \theta y + 1 \}}$, the desired uniqueness remains unclear, to the best of our knowledge.
Nevertheless, we shall demonstrate the application of \Cref{lem: property of ONEC} for uniqueness (albeit limited) of the derived ONEC in the upcoming \Cref{thm: path-independent ONEC}.

\begin{remark}
In comparison, if the terminal condition was replaced by $\rho ( T,y ) = 1$, one would obtain that $( \rho, \varrho ) \equiv ( 1,0 )$ is the unique solution of \cref{eq: BSPDE of rho :eq} by applying the stochastic Feynman-Kac formula.
In this case, we indeed face with the BSDE
${Y}^{t}_{s} = \theta ( {\lambda }_{t} - {X}^{ \bar{\pi } }_{T} ) + \int_{s}^{T} {Y}^{t}_{v} {r}_{v} dv - \int_{s}^{T} \mathcal{Y}^{t}_{v} d {W}_{v}$,
where ${\lambda }_{t} = \mathbb{E}_{t} [ {X}^{ \bar{\pi } }_{T} ] + \frac{1}{\theta }$ satisfies $1 = \mathbb{E}_{t} [ \theta ( \lambda - {X}^{ \bar{\pi } }_{T} ) ]$.
This elegant case is included in \cite{Hu-Jin-Zhou-2017}, with the uniqueness of open-loop equilibrium validated by investigating a non-classical mean-field BSDE.
\end{remark}

\section{Characterization for CNEC}
\label{sec: Characterization for CNEC}

Inspired by \cite{Yong-2012,Wang-Liu-Bensoussan-Yiu-Wei-2025} that characterize CNECs by constructing martingales and taking Feynman-Kac formula, 
we introduce the following extended objective function indexed by $\lambda \in \mathbb{R}$:
\begin{equation}
\label{eq: extended preference with closed-loop control :eq}
   \tilde{J} ( t,x, \Pi, \lambda )
:= \lambda ( 1 - \mathbb{E}_{t} [ \zeta ] )
 - \frac{\theta }{2} \mathbb{E}_{t} \bigg[ \Big| \Big( \lambda - {X}^{ t,x, \Pi }_{T} - \frac{\zeta }{\theta } \Big)_{+} \Big|^{2} \bigg]
 + \mathbb{E}_{t} \Big[ {X}^{ t,x, \Pi }_{T} \zeta + \frac{1}{ 2 \theta } ( {\zeta }^{2} - 1 ) \Big],
\end{equation}
so that $\bar{J} ( t,x, \Pi ) = \tilde{J} ( t,x, \Pi, \lambda ( \mathbb{P}^{t}_{ {X}^{ t,x, \Pi }_{T}, \zeta } ) )$ follows from \cref{eq: SMMV preference with closed-loop control :eq} and \cref{eq: re-expression of SMMV preference :eq}.
Extending existing results using classical PDEs in, e.g., \cite{Bjork-Khapko-Murgoci-2017,Wang-Liu-Bensoussan-Yiu-Wei-2025}, we have the following verification theorem to characterize CNECs, 
including the so-called ``extended HJB equations'' or ``extended HJB system'' formulated by the BSPDEs \cref{eq: BSPDE :eq} and \cref{eq: BSPDE of calM :eq} with the equilibrium condition \cref{eq: equilibrium condition for CNEC :eq}.

\begin{theorem}\label{thm: verification theorem for CNEC}
Suppose that $\bar{\Pi } \in \mathcal{U}$.
Then, $\bar{\Pi }$ is a CNEC, if the following conditions hold:
\begin{itemize}
\item There is a random field pair 
      \begin{equation*}
      ( \mathcal{V}, \psi ) \in {C}_{\mathbb{F}} ( 0,T; \mathbb{L}^{2} ( \Omega; {C}^{2,1} ( \mathbb{R} \times \mathbb{R}; \mathbb{R} ) ) ) 
                         \times \mathbb{L}^{2}_{\mathbb{F}} ( 0,T; \mathbb{L}^{2} ( \Omega; {C}^{2,1} ( \mathbb{R} \times \mathbb{R}; \mathbb{R} ) ) )
      \end{equation*}
      fulfilling the following BSPDE on $[ 0,T ) \times \mathbb{R} \times \mathbb{R}$ with terminal condition on $\mathbb{R} \times \mathbb{R}$:
      \begin{equation}
      \label{eq: BSPDE :eq}
      \left\{ \begin{aligned}
      - d \mathcal{V} ( t,x, \lambda )
      & = \bigg( \frac{1}{2} \mathcal{V}_{xx} ( t,x, \lambda ) | \bar{\Pi } ( t,x ) {\sigma }_{t} |^{2} 
               + {\psi }_{x} ( t,x, \lambda ) \bar{\Pi } ( t,x ) {\sigma }_{t} \\
      & \qquad + \mathcal{V}_{x} ( t,x, \lambda ) \big( {r}_{t} x + \bar{\Pi } ( t,x ) {\sigma }_{t} {\vartheta }_{t} \big) \bigg) dt
        - \psi ( t,x, \lambda ) d {W}_{t}, \\
          \mathcal{V} ( T,x, \lambda ) 
      & = \lambda ( 1 - \zeta ) - \frac{\theta }{2} \Big| \Big( \lambda - x - \frac{\zeta }{\theta } \Big)_{+} \Big|^{2} + \zeta x + \frac{1}{ 2 \theta } {\zeta }^{2} - \frac{1}{ 2 \theta }.
      \end{aligned} \right.
      \end{equation}
\item There is a random field pair 
      \begin{equation*}
      ( \mathcal{M}, \phi ) \in {C}_{\mathbb{F}} ( 0,T; \mathbb{L}^{2} ( \Omega; {C}^{2,0} ( \mathbb{R} \times \mathbb{R}; \mathbb{R} ) ) ) 
                         \times \mathbb{L}^{2}_{\mathbb{F}} ( 0,T; \mathbb{L}^{2} ( \Omega; {C}^{2,0} ( \mathbb{R} \times \mathbb{R}; \mathbb{R} ) ) )
      \end{equation*}
      fulfilling the following BSPDE on $[ 0,T ) \times \mathbb{R} \times \mathbb{R}$ with the terminal condition on $\mathbb{R} \times \mathbb{R}$:
      \begin{equation}
      \label{eq: BSPDE of calM :eq}
      \left\{ \begin{aligned}
      - d \mathcal{M} ( t,x, \lambda )
      & = \bigg( \frac{1}{2} \mathcal{M}_{xx} ( t,x, \lambda ) | \bar{\Pi } ( t,x ) {\sigma }_{t} |^{2} 
               + {\phi }_{x} ( t,x, \lambda ) \bar{\Pi } ( t,x ) {\sigma }_{t} \\
      & \qquad + \mathcal{M}_{x} ( t,x, \lambda ) \big( {r}_{t} x + \bar{\Pi } ( t,x ) {\sigma }_{t} {\vartheta }_{t} \big) \bigg) dt
        - \phi ( t,x, \lambda ) d {W}_{t}, \\
          \mathcal{M} ( T, x, \lambda ) 
      & = \zeta + \theta \Big( \lambda - x - \frac{\zeta }{\theta } \Big)_{+}.
      \end{aligned} \right.
      \end{equation}
\item $\{ f ( s, {X}^{ t,x, \bar{\Pi } }_{s}, \lambda ) \}_{ s \in [ t,T ] } \in \mathbb{L}^{2}_{\mathbb{F}} ( t,T; \mathbb{L}^{2} ( \Omega; \mathbb{R} ) )$ for $f = \mathcal{V}_{x}, \psi, \mathcal{M}_{x}, \phi $ and any $( t,x, \lambda )$.
\item With $\Lambda ( t,x )$ satisfying $1 = \mathcal{M} ( t,x, \Lambda ( t,x ) )$,
      \begin{equation}
      \label{eq: equilibrium condition for CNEC :eq}
      \bar{\Pi } ( t,x ) \in \mathop{argmax}_{ \pi \in \mathbb{R} } 
                             \bigg\{ \frac{1}{2} \mathcal{V}_{xx} \big( t,x, \Lambda ( t,x ) \big) | \pi {\sigma }_{t} |^{2} 
                                   + \Big( \mathcal{V}_{x} \big( t,x, \Lambda ( t,x ) \big) {\vartheta }_{t} + {\psi }_{x} \big( t,x, \Lambda ( t,x ) \big) \Big) \pi {\sigma }_{t} \bigg\}.
      \end{equation}  
\end{itemize}
\end{theorem}

\begin{proof}
See \Cref{pf-thm: verification theorem for CNEC}.
\end{proof}

Given that the terminal value $\mathcal{M} ( T, x, \lambda ) = \zeta + \theta ( \lambda - x - \frac{\zeta }{\theta } )_{+}$ is not continuously differentiable,
we have indeed assumed that $\mathcal{M} ( t, \cdot, \lambda )$ is twice differentiable for $t \in [ 0,T )$,
analogous to deriving the Black-Scholes-Merton equation for option pricing (e.g., \cite[Section 4.5]{Shreve-2004}).
In fact, this additional differentiability assumption holds, once ${X}^{ t,x, \bar{\Pi } }_{T}$ has a probability density.
Moreover, it exactly produces meaningful insights to adopt this differentiability assumption.
For example, under additional smoothness and integrability assumptions, the following \Cref{thm: identical trivial NECs} shows that a state-independent CNEC does necessarily produce an ONEC.
Notably, provided that $\bar{\Pi } ( t,x ) \equiv \bar{\pi }_{t}$ on $[ 0,T ] \times \Omega \times \mathbb{R}$ provides a CNEC, 
\begin{equation*}
0 \le \liminf_{\varepsilon \downarrow 0}
      \frac{1}{\varepsilon } \big( J ( t, \bar{\pi } ) - J ( t, \bar{\pi } {1}_{ [ t, t + \varepsilon ) } + \xi {1}_{ [ t, t + \varepsilon ) } ) \big)
\end{equation*}
given by \Cref{def: closed-loop Nash equilibrium control} does not necessarily lead to
\begin{equation*}
0 \le \liminf_{\varepsilon \downarrow 0}
      \frac{1}{\varepsilon } \big( J ( t, \bar{\pi } ) - J ( t, \bar{\pi } + \xi {1}_{ [ t, t + \varepsilon ) } ) \big)
\end{equation*}
as in \Cref{def: open-loop Nash equilibrium control} for ONEC, due to the difference in spike variation.

\begin{theorem}\label{thm: identical trivial NECs}
Suppose that for $\bar{\Pi } ( t,x ) \equiv \bar{\pi }_{t}$ on $[ 0,T ] \times \Omega \times \mathbb{R}$ with $\bar{\pi } \in {C}_{\mathbb{F}} ( 0,T; \mathbb{L}^{2} ( \Omega; \mathbb{R} ) )$,
there are random field pairs 
\begin{equation*}
\left\{ \begin{aligned}
( \mathcal{V}, \psi ) & \in {C}_{\mathbb{F}} ( 0,T; \mathbb{L}^{2} ( \Omega; {C}^{3,1} ( \mathbb{R} \times \mathbb{R}; \mathbb{R} ) ) ) 
                     \times {C}_{\mathbb{F}} ( 0,T; \mathbb{L}^{2} ( \Omega; {C}^{3,1} ( \mathbb{R} \times \mathbb{R}; \mathbb{R} ) ) ), \\
( \mathcal{M}, \phi ) & \in {C}_{\mathbb{F}} ( 0,T; \mathbb{L}^{2} ( \Omega; {C}^{2,0} ( \mathbb{R} \times \mathbb{R}; \mathbb{R} ) ) ) 
                     \times \mathbb{L}^{2}_{\mathbb{F}} ( 0,T; \mathbb{L}^{2} ( \Omega; {C}^{2,0} ( \mathbb{R} \times \mathbb{R}; \mathbb{R} ) ) ),
\end{aligned} \right.
\end{equation*}
fulfilling \cref{eq: BSPDE :eq}, \cref{eq: BSPDE of calM :eq} and 
\begin{equation}
\label{eq: equilibrium condition for trivial NEC :eq}
0 = \mathcal{V}_{xx} \big( t,x, \lambda ( \mathbb{P}^{t}_{ {X}^{ t,x, \bar{\Pi } }_{T}, \zeta } ) \big) \bar{\Pi } ( t,x ) {\sigma }_{t} 
  + \mathcal{V}_{x} \big( t,x, \lambda ( \mathbb{P}^{t}_{ {X}^{ t,x, \bar{\Pi } }_{T}, \zeta } ) \big) {\vartheta }_{t} 
  + {\psi }_{x} \big( t,x, \lambda ( \mathbb{P}^{t}_{ {X}^{ t,x, \bar{\Pi } }_{T}, \zeta } ) \big).
\end{equation}
In addition, assume that $\{ f ( s, {X}^{ t,x, \bar{\Pi } }_{s}, \lambda ) \}_{ s \in [ t,T ] } \in \mathbb{L}^{2}_{\mathbb{F}} ( t,T; \mathbb{L}^{2} ( \Omega; \mathbb{R} ) )$ 
for $f = \mathcal{V}, \psi, \mathcal{V}_{x}, {\psi }_{x}, \mathcal{M}, \phi $ and any $( t,x, \lambda )$.
Then, $\bar{\Pi }$ is a CNEC, and $\bar{\pi }$ is an ONEC.
\end{theorem}

\begin{proof}
See \Cref{pf-thm: identical trivial NECs}.
\end{proof}

\section{Semi-closed-form solution for deterministic parameters}
\label{sec: closed-form solution}

In the case that the model parameters are all deterministic, we obtain a semi-closed-form expression of the ONEC, which is a deterministic function of $t$.
See the upcoming \Cref{thm: path-independent ONEC}.
Moreover, by exploiting \Cref{lem: property of ONEC}, we can show the uniqueness of the ONEC in a certain range.

\begin{theorem}\label{thm: path-independent ONEC}
Assume that $( r, \sigma, \vartheta, \zeta )$ are deterministic.
Then, $\bar{\pi }$ is an ONEC if 
\begin{equation}
\label{eq: ONEC :eq}
\bar{\pi }_{t} = \frac{ {\beta }_{t} {\vartheta }_{t} }{ \theta {\sigma }_{t} } {e}^{ - \int_{t}^{T} {r}_{v} dv }, \quad \mathbb{P}-a.s., ~ a.e. ~ t \in [ 0,T ],
\end{equation} 
where the continuous and decreasing function $\beta : [ 0,T ] \to [ 1, + \infty )$ is the unique solution of the following integral equation with ${\beta }_{T} = 1$:
\begin{equation}
\label{eq: beta equation :eq}
1 = {\beta }_{t} \Phi \circ \Psi \bigg( \frac{ 1 - \zeta }{ ( \int_{t}^{T} | {\beta }_{s} {\vartheta }_{s} |^{2} ds )^{ \frac{1}{2} } } \bigg), \quad t \in [ 0,T ).
\end{equation}
Moreover, for all $t \in [ 0,T )$, ${\beta }_{t}$ is strictly increasing in $\zeta $.
  
Conversely, if $\bar{\pi }$ is an ONEC satisfying that 
\begin{itemize}
\item $( P, \eta )$ given by \cref{eq: martingale representation for indicator :eq} is of
      ${C}_{\mathbb{F}} ( 0,T; \mathbb{L}^{2} ( \Omega; {C}^{2} ( \mathbb{R}; \mathbb{R} ) ) ) \times \mathbb{L}^{2}_{\mathbb{F}} ( 0,T; \mathbb{L}^{2} ( \Omega; {C}^{1} ( \mathbb{R}; \mathbb{R} ) ) )$
\item and that $\rho ( \cdot, 0 )$ and ${\rho }_{y} ( \cdot, 0 )$ given by \cref{eq: rho-varrho :eq} are deterministic functions,
\end{itemize}
then $\bar{\pi }$ satisfies \cref{eq: ONEC :eq} with $\beta $ given by \cref{eq: beta equation :eq}.
\end{theorem}

\begin{proof}
See \Cref{pf-thm: path-independent ONEC}.
\end{proof}

\begin{remark}
Given that $\lambda ( \mathbb{P}^{T}_{ {X}^{ \bar{\pi } }_{T}, \zeta } ) = {X}^{ \bar{\pi } }_{T} + \frac{1}{\theta }$,
a deterministic $\rho ( \cdot, 0 )$ implies that the conditional probability of
$\lambda ( \mathbb{P}^{T}_{ {X}^{ \bar{\pi } }_{T}, \zeta } ) - \lambda ( \mathbb{P}^{t}_{ {X}^{ \bar{\pi } }_{T}, \zeta } ) \le \frac{ 1 - \zeta }{\theta }$ 
(conditioned on $\mathcal{F}_{t}$) is a deterministic function of $t$.
Similarly, a deterministic ${\rho }_{y} ( \cdot, 0 )$ implies that the conditional probability density of
$\lambda ( \mathbb{P}^{T}_{ {X}^{ \bar{\pi } }_{T}, \zeta } ) - \lambda ( \mathbb{P}^{t}_{ {X}^{ \bar{\pi } }_{T}, \zeta } ) = \frac{ 1 - \zeta }{\theta }$ 
(conditioned on $\mathcal{F}_{t}$) is also a deterministic function of $t$.
\end{remark}

For $t \in [ 0,T )$, ${\beta }_{t} > 1$ implies that $\bar{\pi }_{t} > \bar{\pi }^{MV}_{t} := \frac{ {\vartheta }_{t} }{ \theta {\sigma }_{t} } \exp ( - \int_{t}^{T} {r}_{v} dv )$, 
where $\bar{\pi }^{MV}$ is the time-consistent strategy for the MV preference $\mathbb{E}_{t} [ {X}^{\pi }_{T} ] - \frac{\theta }{2} \mathbb{E}_{t} [ ( {X}^{\pi }_{T} - \mathbb{E}_{t} [ {X}^{\pi }_{T} ] )^{2} ]$; 
see also \cite{Basak-Chabakauri-2010,Hu-Jin-Zhou-2012,Bjork-Khapko-Murgoci-2017}. 
And we observe that the equilibrium solution \cref{eq: ONEC :eq} is irrelevant with the wealth level, 
which aligns with the time-consistent solution under MV preferences, as shown in the aforementioned literature.
Since our SMMV preferences are monotone modifications of the MV preference, ${\beta }_{t}$ can be treated as the modification factor for the vanishing penalty on high positive deviations of investment return.
Intuitively speaking, a larger $\zeta$ encourages greater high positive deviations, which in turn increases $\beta_{t}$ and consequently leads to a larger equilibrium investment amount $\bar{\pi}_{t}$.
In addition, in the case with a constant risk premium $\vartheta $, \cref{eq: beta equation :eq} implies that ${\beta }_{t}$ is also strictly increasing in $\vartheta $, 
and hence a larger $\vartheta $ leads to a larger $\bar{\pi }_{t}$; that is, higher returns of risk attract more risky investments.

On the other hand, the following \Cref{thm: state-independent CNEC} gives a state-independent CNEC identical to the path-independent ONEC that we have derived in \Cref{thm: path-independent ONEC}.
Notably, corresponding to this CNEC, the equilibrium value function $\mathcal{V}$ and the random field $\mathcal{M}$ satisfy the differentiability condition stated in \Cref{thm: verification theorem for CNEC}.

\begin{theorem}\label{thm: state-independent CNEC}
Assume that $( r, \sigma, \vartheta, \zeta )$ are deterministic.
Then, $\bar{\Pi }$ is a CNEC if 
\begin{equation}
\label{eq: CNEC :eq}
\bar{\Pi } ( t,x ) = \frac{ {\beta }_{t} {\vartheta }_{t} }{ \theta {\sigma }_{t} } {e}^{ - \int_{t}^{T} {r}_{v} dv }, \quad \forall ( t,x ) \in [ 0,T ] \times \mathbb{R},
\end{equation}
where $\beta $ is given by \cref{eq: beta equation :eq} with ${\beta }_{T} = 1$.
Furthermore, corresponding to this CNEC, 
\begin{equation*}
\bar{J} ( t,x, \bar{\Pi } ) = x {e}^{ \int_{t}^{T} {r}_{v} dv } - \frac{ ( 1 - \zeta )^{2} }{ 2 \theta }
                              + \frac{1}{\theta } \int_{t}^{T} {\beta }_{s} | {\vartheta }_{s} |^{2} ds
                              + \frac{1}{\theta } \int_{t}^{T} | {\beta }_{s} {\vartheta }_{s} |^{2} ds \int_{0}^{ \frac{ 1 - \zeta }{ ( \int_{t}^{T} | {\beta }_{s} {\vartheta }_{s} |^{2} ds )^{\frac{1}{2}} } } \Psi(y) dy.
\end{equation*}
\end{theorem}

\begin{proof}
See \Cref{pf-thm: state-independent CNEC}.
\end{proof}

\begin{remark}
It is found that the CNEC (\Cref{thm: state-independent CNEC}) and the ONEC (\Cref{thm: path-independent ONEC}) are the same and only rely on $t$, which is similar to the findings in \cite[Section 5]{Wang-Liu-Bensoussan-Yiu-Wei-2025},
where the performance functional is a fairly general function of finitely many higher-order central moments of ${X}_{T}$ and is affine in $\mathbb{E}_{t} [ {X}_{T} ]$.
Heuristically, here we show that \cref{eq: SMMV preference with closed-loop control :eq} with a deterministic $\zeta $ can be roughly re-expressed as 
an $\mathbb{E}_{t} [ {X}^{ t,x, \Pi }_{T} ]$-affine objective function with higher-order central moments of ${X}^{ t,x, \Pi }_{T}$.
In fact, $\lambda ( \mathbb{P}^{t}_{ {X}^{ t,x, \Pi }_{T}, \zeta } ) - \mathbb{E}_{t} [ {X}^{ t,x, \Pi }_{T} ]$ is a function of higher-order central moments of ${X}^{ t,x, \Pi }_{T}$, since
\begin{equation*}
1 = \zeta + \theta \mathbb{E}_{t} \Big[ \Big( \big( \lambda ( \mathbb{P}^{t}_{ {X}^{ t,x, \Pi }_{T}, \zeta } ) - \mathbb{E}_{t} [ {X}^{ t,x, \Pi }_{T} ] \big) 
                                            - ( {X}^{ t,x, \Pi }_{T} - \mathbb{E}_{t} [ {X}^{ t,x, \Pi }_{T} ] ) - \frac{\zeta }{\theta } \Big)_{+} \Big].
\end{equation*}
A key underlying assumption is that the distribution of ${X}^{ t,x, \Pi }_{T} - \mathbb{E}_{t} [ {X}^{ t,x, \Pi }_{T} ]$ is uniquely determined by its moments.
Moreover, it follows from \cref{eq: extended preference with closed-loop control :eq} with $\bar{J} ( t,x, \Pi ) = \tilde{J} ( t,x, \Pi, \lambda ( \mathbb{P}^{t}_{ {X}^{ t,x, \Pi }_{T}, \zeta } ) )$ that 
\begin{align*}
    \bar{J} ( t, x, \Pi )
& = \frac{1}{ 2 \theta } ( {\zeta }^{2} - 1 )
  + \mathbb{E}_{t} [ {X}^{ t,x, \Pi }_{T} ] 
  + ( \lambda ( \mathbb{P}^{t}_{ {X}^{ t,x, \Pi }_{T}, \zeta } ) - \mathbb{E}_{t} [ {X}^{ t,x, \Pi }_{T} ] ) ( 1 - \zeta ) \\
& \quad - \frac{\theta }{2} \mathbb{E}_{t} \bigg[ \Big| \Big( \big( \lambda ( \mathbb{P}^{t}_{ {X}^{ t,x, \Pi }_{T}, \zeta } ) - \mathbb{E}_{t} [ {X}^{ t,x, \Pi }_{T} ] \big) 
                                                            - ( {X}^{ t,x, \Pi }_{T} - \mathbb{E}_{t} [ {X}^{ t,x, \Pi }_{T} ] ) - \frac{\zeta }{\theta } \Big)_{+} \Big|^{2} \bigg].
\end{align*}
That is, $\bar{J} ( t, x, \Pi ) - \mathbb{E}_{t} [ {X}^{ t,x, \Pi }_{T} ]$ is also a function of higher-order central moments of ${X}^{ t,x, \Pi }_{T}$.
However, it is still difficult to apply the result of \cite{Wang-Liu-Bensoussan-Yiu-Wei-2025} to derive ONECs and CNECs for our control problem, 
as the aforementioned functions of higher-order central moments of ${X}^{ t,x, \Pi }_{T}$ is implicit in general.
\end{remark}

Taking advantage of the path/state-independence, we can show whether a Nash equilibrium control is a strong equilibrium strategy by straightforward calculation in the case with constant parameters.
This approach is much more clear than using the characterization given by \cite{He-Jiang-2021}.

\begin{proposition}\label{thm: equilibrium value function and strong equilibria}
Assume that $( r, \sigma, \vartheta, \zeta )$ are constant, and $( \bar{\pi }, \bar{\Pi } )$ are the ONEC and CNEC given by \Cref{thm: path-independent ONEC,thm: state-independent CNEC}, respectively.  
Then, for any $t \in [ 0,T )$, $\xi \in \mathbb{L}^{2}_{ \mathcal{F}_{t} } ( \Omega; \mathbb{R} )$ and sufficiently small $( t, \xi )$-dependent $\varepsilon $,
\begin{itemize}
\item $\bar{J} ( t,x, \bar{\Pi } ) > \bar{J} ( t,x, \bar{\Pi }^{ t, \varepsilon, \xi } )$ if $r > \frac{ | {\beta }_{t}' | }{ {\beta }_{t} }$;
\item $\bar{J} ( t,x, \bar{\pi } ) > \bar{J} ( t,x, \bar{\pi }^{ t, \varepsilon, \xi } )$.
\end{itemize}
\end{proposition}

\begin{proof}
See \Cref{pf-thm: equilibrium value function and strong equilibria}.
\end{proof}

\Cref{thm: equilibrium value function and strong equilibria} shows that the CNEC given by \Cref{thm: state-independent CNEC} is a strong equilibrium strategy under the condition $r > \max_{ t \in [ 0,T ] } \frac{ | {\beta }_{t}' | }{ {\beta }_{t} }$.
That is, any spike deviation as in \Cref{def: closed-loop Nash equilibrium control} from the closed-loop equilibrium given by \Cref{thm: state-independent CNEC} will be strictly worse off, only when the interest rate $r$ is sufficiently large.
In comparison, any spike deviation as in \Cref{def: open-loop Nash equilibrium control} from the open-loop equilibrium given by \Cref{thm: path-independent ONEC} must be strictly worse off.

\section{Numerical example}
\label{sec: Numerical example}

This section illustrates the investment strategies derived from Theorems \ref{thm: path-independent ONEC} and \ref{thm: state-independent CNEC}, comparing them against strategies under MV and MMV preferences. 
Assuming the parameters $(r, \sigma, \vartheta, \zeta)$ are constant, we perform a calibration exercise using monthly S\&P 500 stock index from January 1990 to December 2025. 
Data are sourced from the WIND database and supplemented by the 30-year U.S. Treasury bond annual yields for the corresponding period, sourced from the U.S. Department of the Treasury \cite{US-treasury-2026}.

\begin{table}[H]
	\centering
	\caption{Calibrated estimates of the parameters.}
	\begin{tabular}{p{12em}p{12em}p{12em}}
		\toprule
		\textbf{Parameter} & \textbf{Estimate} & \textbf{Standard error}  \\
		\midrule
		$r$ (interest rate) & 0.04601$^{***}$ & 0.00087  \\
		$\sigma$ (volatility) & 0.14775$^{***}$ & 0.00503 \\
		$\vartheta$ (risk premium) & 0.33439$^{*}$ & 0.17146 \\
		\bottomrule
        \multicolumn{3}{l}{\footnotesize $^{***}$, $^{**}$, and $^{*}$ indicate statistical significance at the 1\%, 5\%, and 10\% levels, respectively.}
	\end{tabular}
	\label{tab:calibration}
\end{table}%

As shown in Table \ref{tab:calibration}, the estimated risk-free interest rate $r$ is around 4.60\%, and the volatility $\sigma$ is estimated to be 0.1478. 
The estimated risk premium $\vartheta$ is 0.3344, which implies an expected excess return ($\vartheta \sigma$) of approximately 4.94\%. 
These figures are comparable in magnitude to the estimates reported by \cite{Munk-Sorensen-NygaardVinther-2004}, who utilized historical S\&P 500 data spanning approximately 50 years. 
In our numerical setup, we further consider a long-term investment horizon of $T = 25$ years. 
Regarding the investor's preferences, we set the risk aversion parameter to $\theta = 3$, following the numerical implementation in \cite{Bjork-Murgoci-Zhou-2014}. 
A higher value of $\theta$ indicates greater aversion to negative deviations of terminal wealth from its expectation. 

\begin{figure}[htbp!]
	\centering
	\includegraphics[width=0.8\linewidth]{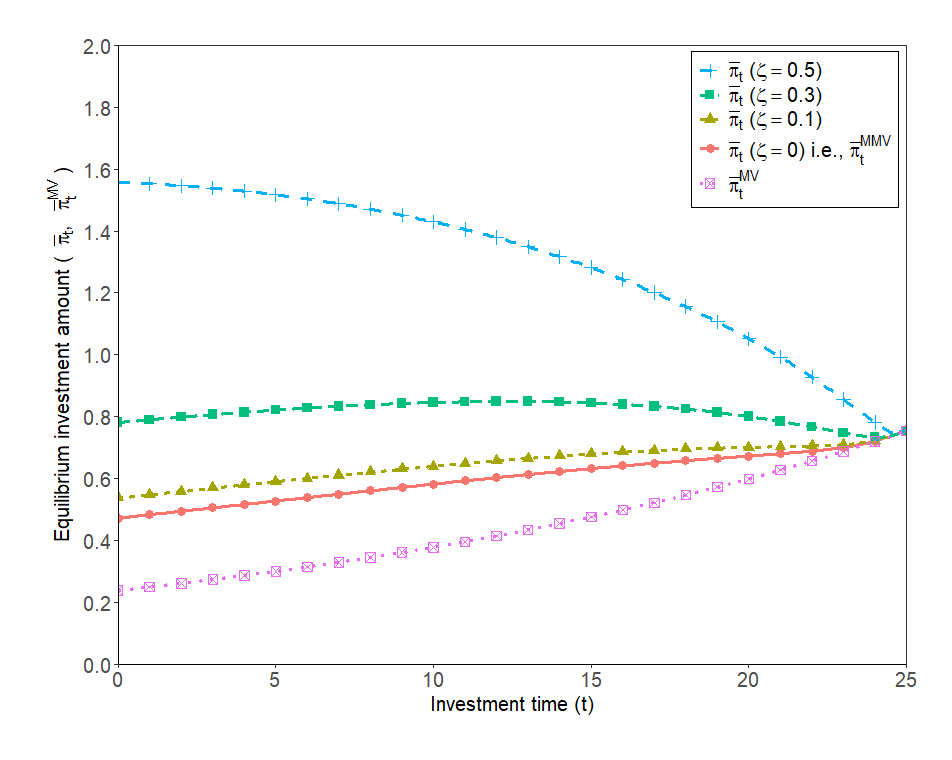}
	\caption{Equilibrium investment amount in the risky asset with different preferences. 
             The parameters are assumed to be deterministic, with $r=0.04601$, $\sigma=0.14775$, $\vartheta=0.33439$, $T = 25$ and $\theta = 3$. $\zeta$ varies from 0 to 0.5.}
	\label{fig: pi plot}
\end{figure}

\begin{figure}[htbp!]
	\centering
	\includegraphics[width=0.8\linewidth]{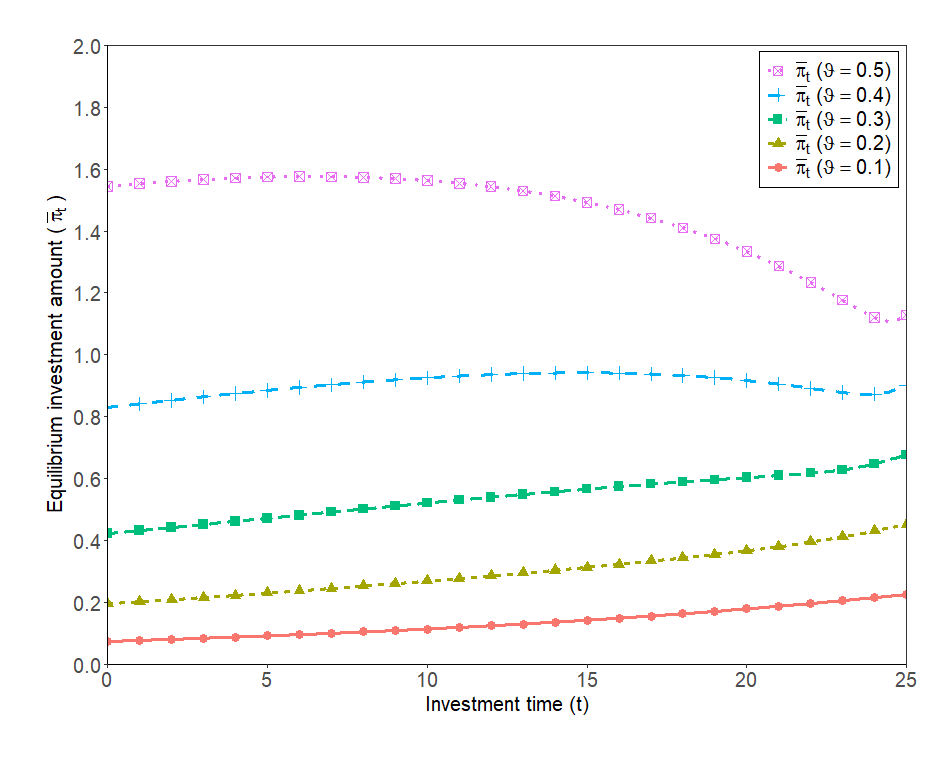}
	\caption{Equilibrium investment amount in the risky asset under the SMMV preference w.r.t. different values of risk premium. 
             The parameters are assumed to be deterministic, with $r=0.04601$, $\sigma=0.14775$, $\zeta = 0.1$, $T = 25$ and $\theta = 3$. The value of the risk premium $\vartheta$ varies from 0.1 to 0.5.}
	\label{fig: pi vartheta plot}
\end{figure}

\Cref{fig: pi plot} illustrates the amount invested in the risky asset with the SMMV and the MV preference, respectively. 
As noted in \Cref{sec: closed-form solution}, since ${\beta }_{t} > 1$, the equilibrium investment in the risky asset with the SMMV preference $\bar{\pi }_{t}$ is higher than that under the MV preference $\bar{\pi }^{MV}_{t}$. 
By sending $\zeta $ to zero, we obtain the equilibrium investment amount under the MMV preference within our time-consistent framework (see the red curve with solid circles in \Cref{fig: pi plot}). 
Based on the estimated parameters, we observe that MV, MMV, and SMMV investors with moderate $\zeta$ (e.g., $\zeta = 0.1$) gradually increase their equilibrium investment amount over time. 
In contrast, a larger value of $\zeta$ induces the investor to utilize significant leverage in the early stages of the horizon, followed by a rapid reduction in risky asset exposure as maturity approaches. 
This shift reflects the SMMV investor's sharply diminishing appetite for risky investments over time. 
By the end of the investment period, at time $T$, the equilibrium investment amounts under both preferences converge and become identical. 
Furthermore, \Cref{fig: pi vartheta plot} exhibits a significant positive effect of the risk premium $\vartheta $ on the equilibrium investment amount $\bar{\pi }$ under SMMV preferences, 
indicating that higher risk premiums encourage greater investment in the risky asset.

\section{Concluding remark}
\label{sec: Concluding remark}

In this paper, we have studied the time-inconsistent portfolio selection problem with MMV preferences, and characterized the time-consistent strategy in two types depending on how game theoretic formulation.
For the problem with deterministic parameters, we obtain analytical expressions of an open-loop equilibrium and a closed-loop equilibrium.
The two equilibrium strategies are found to be identical, which are independent of the wealth level and random path.
Moreover, this time-consistent strategy with the SMMV preference indicates a greater amount of equity investment than the suggested strategy with a conventional MV preference.
We have also studied whether the aforementioned solution is a stronger equilibrium strategy in the case with constant parameters. 
For spike variation in the open-loop type, any deviation for the derived equilibrium will lead to a strict loss in the investor's satisfaction.
But there does not necessarily exist such a strict loss for spike variation in the closed-loop type, unless the interest rate is large enough.

\appendix

\section{Proof of lemmas and theorems}

\subsection{Proof of \texorpdfstring{\Cref{lem: moment estimate for lambda}}{Lemma 2.1}}
\label{pf-lem: moment estimate for lambda}

Here we newly provide a ``mean-variance'' lower bound for $\lambda ( \mathbb{P}^{t}_{ {X}^{\pi }_{T}, \zeta } )$, as the following:
\begin{equation*}
\lambda ( \mathbb{P}^{t}_{ {X}^{\pi }_{T}, \zeta } ) 
\ge \mathbb{E}_{t} [ {X}^{\pi }_{T} ] - \frac{\theta }{\delta } \mathbb{E}_{t} \big[ | {X}^{\pi }_{T} - \mathbb{E}_{t} [ {X}^{\pi }_{T} ] |^{2} \big] + \frac{1}{\theta } - \frac{1}{\theta \delta } \mathbb{E} [ | 1 - \zeta |^{2} ],
\end{equation*}
which arises from 
\begin{align*}
      \frac{\delta }{\theta }
  \le \frac{ 1 - \mathbb{E}_{t} [ \zeta ] }{\theta } 
&   = \int_{ - \infty }^{ \lambda ( \mathbb{P}^{t}_{ {X}^{\pi }_{T}, \zeta } ) - \mathbb{E}_{t} [ {X}^{\pi }_{T} ] - \frac{1}{\theta } } 
      \mathbb{E}_{t} [ {1}_{\{ s \ge {X}^{\pi }_{T} - \mathbb{E}_{t} [ {X}^{\pi }_{T} ] - \frac{ 1 - \zeta }{\theta } \}} ] ds \\
& \le \int_{ - \infty }^{ \lambda ( \mathbb{P}^{t}_{ {X}^{\pi }_{T}, \zeta } ) - \mathbb{E}_{t} [ {X}^{\pi }_{T} ] - \frac{1}{\theta } } 
      \frac{1}{ |s|^{2} } \mathbb{E}_{t} \bigg[ \Big| {X}^{\pi }_{T} - \mathbb{E}_{t} [ {X}^{\pi }_{T} ] - \frac{ 1 - \zeta }{\theta } \Big|^{2} \bigg] ds \\
&   = \frac{1}{ \mathbb{E}_{t} [ {X}^{\pi }_{T} ] + \frac{1}{\theta } - \lambda ( \mathbb{P}^{t}_{ {X}^{\pi }_{T}, \zeta } ) }
      \mathbb{E}_{t} \bigg[ | {X}^{\pi }_{T} - \mathbb{E}_{t} [ {X}^{\pi }_{T} ] |^{2} + \Big| \frac{ 1 - \zeta }{ \theta } \Big|^{2} \bigg].
\end{align*}
Combining the lower bound with the upper bound $\lambda ( \mathbb{P}^{t}_{ {X}^{\pi }_{T}, \zeta } ) \le \mathbb{E}_{t} [ {X}^{\pi }_{T} ] + \frac{1}{\theta }$ immediately delivers the desired result.

\subsection{Proof of \texorpdfstring{\Cref{lem: properties for ONEC}}{Lemma 3.1}}
\label{pf-lem: properties for ONEC}

Write $\Delta = q - p$ for the sake of brevity. It follows from \cref{eq: lambda functional :eq} that
\begin{align*}
0 & = \int_{ \mathbb{R} \times \mathbb{R} } \bigg( \Big( \lambda ( p + \varepsilon \Delta ) - x - \frac{z}{\theta } \Big)_{+} - \Big( \lambda (p) - x -\frac{z}{\theta } \Big)_{+} \bigg) p ( dx, dz ) \\
  & \quad + \varepsilon \int_{ \mathbb{R} \times \mathbb{R} } \bigg( \frac{z}{\theta } + \Big( \lambda ( p + \varepsilon \Delta ) - x -\frac{z}{\theta } \Big)_{+} \bigg) \Delta ( dx, dz ) \\
  & = \int_{ \lambda (p) }^{ \lambda ( p + \varepsilon \Delta ) } ds \int_{\{ x + \frac{z}{\theta } \le s \}} p ( dx, dz ) 
    + \varepsilon \int_{ \mathbb{R} \times \mathbb{R} } \bigg( \frac{z}{\theta } + \Big( \lambda ( p + \varepsilon \Delta ) - x -\frac{z}{\theta } \Big)_{+} \bigg) \Delta ( dx, dz ),
\end{align*}
which implies that $\lambda ( p + \varepsilon \Delta ) \to \lambda (p)$ as $\varepsilon \to 0$.
Moreover, $\int_{\{ x + \frac{z}{\theta } \le s \}} p ( dx, dz )$ as a cumulative distribution function is c\'adl\'ag, 
and hence, for sufficiently small $\varepsilon > 0$, we have
\begin{align*}
D ( p, \varepsilon, \Delta ) := \frac{1}{ \lambda ( p + \varepsilon \Delta ) - \lambda (p) } \int_{ \lambda (p) }^{ \lambda ( p + \varepsilon \Delta ) } ds \int_{\{ x + \frac{z}{\theta } \le s \}} p ( dx, dz ) \in ( 0,1 ],
\end{align*}
and 
\begin{equation*}
0 < \int_{\{ x + \frac{z}{\theta } < \lambda (p) \}} p ( dx, dz ) 
  \le \liminf_{ \varepsilon \downarrow 0 } D ( p, \varepsilon, \Delta )
  \le \limsup_{ \varepsilon \downarrow 0 } D ( p, \varepsilon, \Delta )
  \le \int_{\{ x + \frac{z}{\theta } \le \lambda (p) \}} p ( dx, dz )
  \le 1.
\end{equation*}
Consequently, for sufficiently small $\varepsilon > 0$
\begin{equation*}
  | \lambda ( p + \varepsilon \Delta ) - \lambda (p) | 
= \frac{\varepsilon }{ D ( p, \varepsilon, \Delta ) }
  \bigg| \int_{ \mathbb{R} \times \mathbb{R} } \bigg( \frac{z}{\theta } + \Big( \lambda (p) - x -\frac{z}{\theta } \Big)_{+} \bigg) \Delta ( dx, dz ) \bigg|
+ o ( \varepsilon )
\end{equation*}
On the other hand, for \cref{eq: non-linear objective functional :eq} we have
\begin{align*}
& g ( p + \varepsilon \Delta ) - g (p) \\
& = \big( \lambda ( p + \varepsilon \Delta ) - \lambda (p) \big) \int_{ \mathbb{R} \times \mathbb{R} } ( 1- z ) d p ( dx, dz ) \\
& \quad - \frac{\theta }{2} \int_{ \mathbb{R} \times \mathbb{R} } \Big| \Big( \lambda ( p + \varepsilon \Delta ) - x - \frac{z}{\theta } \Big)_{+} \Big|^{2} p ( dx, dz )
        + \frac{\theta }{2} \int_{ \mathbb{R} \times \mathbb{R} } \Big| \Big( \lambda (p) - x - \frac{z}{\theta } \Big)_{+} \Big|^{2} p ( dx, dz ) \\
& \quad + \varepsilon \int_{ \mathbb{R} \times \mathbb{R} } \bigg( \lambda ( p + \varepsilon \Delta ) ( 1 - z ) 
                                                                 - \frac{\theta }{2} \Big| \Big( \lambda ( p + \varepsilon \Delta ) - x - \frac{z}{\theta } \Big)_{+} \Big|^{2} 
                                                                 + x z + \frac{1}{ 2 \theta } ( {z}^{2} - 1 ) \bigg) \Delta ( dx, dz ) \\
& = \big( \lambda ( p + \varepsilon \Delta ) - \lambda (p) \big) \int_{ \mathbb{R} \times \mathbb{R} } ( 1- z ) d p ( dx, dz ) \\
& \quad - \theta \int_{ - \infty }^{ \lambda (p) } \big( \lambda ( p + \varepsilon \Delta ) - \lambda (p) \big) ds \int_{\{ x + \frac{z}{\theta } \le s \}} p ( dx, dz ) \\
& \quad - \theta \int_{ \lambda (p) }^{ \lambda ( p + \varepsilon \Delta ) } \big( \lambda ( p + \varepsilon \Delta ) - s \big) ds \int_{\{ x + \frac{z}{\theta } \le s \}} p ( dx, dz ) \\
& \quad + \varepsilon \int_{ \mathbb{R} \times \mathbb{R} } \nabla g ( p,x,z ) \Delta ( dx, dz )
        + o ( \varepsilon ). 
\end{align*}
Plugging \cref{eq: lambda functional equivalent expression :eq} with the asymptotic estimation
\begin{align*}
    \bigg| \int_{ \lambda (p) }^{ \lambda ( p + \varepsilon \Delta ) } \big( \lambda ( p + \varepsilon \Delta ) - s \big) ds \int_{\{ x + \frac{z}{\theta } \le s \}} p ( dx, dz ) \bigg|
& \le \bigg| \int_{ \lambda (p) }^{ \lambda ( p + \varepsilon \Delta ) } | \lambda ( p + \varepsilon \Delta ) - s | ds \bigg| \\
&   = \frac{1}{2} | \lambda ( p + \varepsilon \Delta ) - \lambda (p) |^{2},
\end{align*}
into the above statement yields $g ( p + \varepsilon \Delta ) - g (p) = \varepsilon \int_{ \mathbb{R} \times \mathbb{R} } \nabla g ( p,x,z ) \Delta ( dx, dz ) + o ( \varepsilon )$.
Then, the first desired result immediately arises.

Similarly, we have
\begin{align*}
&   g(q) - g(p) \\
& = d g ( p, q-p )
  + \big( \lambda (q) - \lambda (p) \big) \int_{ \mathbb{R} \times \mathbb{R} } ( 1- z ) d q ( dx, dz ) \\
& \quad - \frac{\theta }{2} \int_{ \mathbb{R} \times \mathbb{R} } \Big| \Big( \lambda (q) - x - \frac{z}{\theta } \Big)_{+} \Big|^{2} q ( dx, dz )
        + \frac{\theta }{2} \int_{ \mathbb{R} \times \mathbb{R} } \Big| \Big( \lambda (p) - x - \frac{z}{\theta } \Big)_{+} \Big|^{2} q ( dx, dz ) \\
& = d g ( p, q-p ) + \int_{ \lambda (q) }^{ \lambda (p) } \big( \lambda (p) - s \big) ds \int_{\{ x + \frac{z}{\theta } \le s \}} q ( dx, dz ) \\
& \quad + \big( \lambda (q) - \lambda (p) \big) \int_{ \mathbb{R} \times \mathbb{R} } ( 1- z ) d q ( dx, dz ) 
        - \theta \int_{ - \infty }^{ \lambda (q) } \big( \lambda (q) - \lambda (p) \big) ds \int_{\{ x + \frac{z}{\theta } \le s \}} q ( dx, dz ),
\end{align*}
where the last line vanishes, with
\begin{equation*}
    \bigg| \int_{ \lambda (q) }^{ \lambda (p) } \big( \lambda (p) - s \big) ds \int_{\{ x + \frac{z}{\theta } \le s \}} q ( dx, dz ) \bigg| 
\le \bigg| \int_{ \lambda (q) }^{ \lambda (p) } | \lambda (p) - s | ds \bigg|
  = \frac{1}{2} | \lambda (q) - \lambda (p) |^{2}.
\end{equation*}
So the second desired result follows.

In terms of the third assertion, we proceed with the following estimate arising from \cref{eq: lambda functional :eq}:
\begin{align*}
&   \lambda ( \mathbb{P}^{t}_{ {X}^{\pi }_{T} + \chi, \zeta } ) - \lambda ( \mathbb{P}^{t}_{ {X}^{\pi }_{T}, \zeta } ) \\
& = \mathbb{E}_{t} \Big[ \lambda ( \mathbb{P}^{t}_{ {X}^{\pi }_{T} + \chi, \zeta } ) \wedge \Big( {X}^{\pi }_{T} + \chi + \frac{\zeta }{\theta } \Big)  
                       - \lambda ( \mathbb{P}^{t}_{ {X}^{\pi }_{T}, \zeta } ) \wedge \Big( {X}^{\pi }_{T} + \frac{\zeta }{\theta } \Big) \Big] \\
& = \mathbb{E}_{t} \Big[ \Big( \lambda ( \mathbb{P}^{t}_{ {X}^{\pi }_{T} + \chi, \zeta } ) - {X}^{\pi }_{T} - \frac{\zeta }{\theta } \Big) {1}_{ {R}_{0} \cap ( \Omega \setminus {R}_{1} ) }
                       + \Big( {X}^{\pi }_{T} + \chi + \frac{\zeta }{\theta } - \lambda ( \mathbb{P}^{t}_{ {X}^{\pi }_{T}, \zeta } ) \Big) {1}_{ ( \Omega \setminus {R}_{0} ) \cap {R}_{1} } \Big] \\
& \quad + \mathbb{E}_{t} [ \chi {1}_{ {R}_{0} \cap {R}_{1} }  ]
        + \big( \lambda ( \mathbb{P}^{t}_{ {X}^{\pi }_{T} + \chi, \zeta } ) - \lambda ( \mathbb{P}^{t}_{ {X}^{\pi }_{T}, \zeta } ) \big) \mathbb{E}_{t} [ {1}_{ ( \Omega \setminus {R}_{0} ) \cap ( \Omega \setminus {R}_{1} ) } ] \\
& \le \mathbb{E}_{t} [ \chi {1}_{ {R}_{0} } ] + \big( \lambda ( \mathbb{P}^{t}_{ {X}^{\pi }_{T} + \chi, \zeta } ) - \lambda ( \mathbb{P}^{t}_{ {X}^{\pi }_{T}, \zeta } ) \big) \mathbb{E}_{t} [ {1}_{ \Omega \setminus {R}_{0} } ].
\end{align*}
Consequently, $\lambda ( \mathbb{P}^{t}_{ {X}^{\pi }_{T} + \chi, \zeta } ) - \lambda ( \mathbb{P}^{t}_{ {X}^{\pi }_{T}, \zeta } ) \le \frac{ \mathbb{E}_{t} [ \chi {1}_{ {R}_{0} } ] }{ \mathbb{E}_{t} [ {1}_{ {R}_{0} } ] }$,
as ${R}_{0}$ is not a $\mathbb{P}$-null set.
In the same manner, we can show that $\lambda ( \mathbb{P}^{t}_{ {X}^{\pi }_{T} + \chi, \zeta } ) - \lambda ( \mathbb{P}^{t}_{ {X}^{\pi }_{T}, \zeta } ) \ge \frac{ \mathbb{E}_{t} [ \chi {1}_{ {R}_{1} } ] }{ \mathbb{E}_{t} [ {1}_{ {R}_{1} } ] }$.
Thus, the proof is completed.

\subsection{Proof of \texorpdfstring{\Cref{thm: asymptotic estimate for ONEC}}{Theorem 3.2}}
\label{pf-thm: asymptotic estimate for ONEC}

Fix $t \in [ 0,T )$ and $\xi \in \mathbb{L}^{2}_{ \mathcal{F}_{t} } ( \Omega; \mathbb{R} )$.
Applying the third property in \Cref{lem: properties for ONEC} to ${X}_{T} = {X}^{ \bar{\pi } }_{T}$ and 
$\chi = \xi \int_{t}^{ t + \varepsilon } \exp ( \int_{s}^{T} {r}_{v} dv ) {\sigma }_{s} ( d {W}_{s} + {\vartheta }_{s} ds ) = \bar{X}^{\varepsilon }_{T} - {X}^{ \bar{\pi } }_{T}$,
where $\bar{X}^{\varepsilon } := {X}^{ \bar{\pi }^{ t, \varepsilon, \xi } }$ for short, yields
\begin{equation*}
    | \lambda ( \mathbb{P}^{t}_{ \bar{X}^{\varepsilon }_{T}, \zeta } ) - \lambda ( \mathbb{P}^{t}_{ {X}^{ \bar{\pi } }_{T}, \zeta } ) | 
\le \frac{ | \mathbb{E}_{t} [ \chi {1}_{ {R}_{0} } ] | }{ \mathbb{E}_{t} [ {1}_{ {R}_{0} } ] }
  + \bigg| \frac{ \mathbb{E}_{t} [ \chi {1}_{ {R}_{0} } ] }{ \mathbb{E}_{t} [ {1}_{ {R}_{0} } ] } - \frac{ \mathbb{E}_{t} [ \chi {1}_{ {R}_{1} } ] }{ \mathbb{E}_{t} [ {1}_{ {R}_{1} } ] } \bigg|.
\end{equation*}
As $\varepsilon \downarrow 0$, we have $\chi \to 0$, $\mathbb{P}$-a.s., and hence, ${1}_{ {R}_{1} } \to {1}_{ {R}_{0} }$, $\mathbb{P}$-a.s.
By the dominated convergence theorem, $\mathbb{E}_{t} [ {1}_{ {R}_{1} } ] \to \mathbb{E}_{t} [ {1}_{ {R}_{0} } ]$, $\mathbb{P}$-a.s.
Consequently, by Cauchy-Schwarz inequality, 
\begin{equation*}
\bigg| \frac{ \mathbb{E}_{t} [ \chi {1}_{ {R}_{0} } ] }{ \mathbb{E}_{t} [ {1}_{ {R}_{0} } ] } - \frac{ \mathbb{E}_{t} [ \chi {1}_{ {R}_{1} } ] }{ \mathbb{E}_{t} [ {1}_{ {R}_{1} } ] } \bigg|^{2}
\le \frac{ \mathbb{E}_{t} [ {\chi }^{2} ] }{ | \mathbb{E}_{t} [ {1}_{ {R}_{0} } ] \mathbb{E}_{t} [ {1}_{ {R}_{1} } ] |^{2} } 
    \mathbb{E}_{t} \big[ | {1}_{ {R}_{0} } \mathbb{E}_{t} [ {1}_{ {R}_{1} } ] - {1}_{ {R}_{1} } \mathbb{E}_{t} [ {1}_{ {R}_{0} } ] |^{2} \big]
  = o ( \varepsilon ).
\end{equation*}
Since ${R}_{0}$ is independent of $( \varepsilon, \xi )$, mirroring the proof of \cite[Lemma 4.1]{Wang-Liu-Bensoussan-Yiu-Wei-2025}, one can obtain $| \mathbb{E}_{t} [ \chi {1}_{ {R}_{0} } ] |^{2} = o ( \varepsilon )$.
Therefore, $| \lambda ( \mathbb{P}^{t}_{ \bar{X}^{\varepsilon }_{T}, \zeta } ) - \lambda ( \mathbb{P}^{t}_{ {X}^{ \bar{\pi } }_{T}, \zeta } ) |^{2} = o ( \varepsilon )$.
Furthermore, from the second property in \Cref{lem: properties for ONEC}, together with applying Taylor expansion to $\nabla g ( p, \cdot, z )$ 
that is continuously differentiable on $\mathbb{R}$ and twice continuously differentiable on $( - \infty, \lambda (p) - \frac{z}{\theta } )$ and on $( \lambda (p) - \frac{z}{\theta }, + \infty )$
and applying It\^o's rule to ${Y}^{t}_{T} ( \bar{X}^{\varepsilon }_{T} - {X}^{ \bar{\pi } }_{T} )$, 
we obtain
\begin{align*}
    g ( \mathbb{P}^{t}_{ \bar{X}^{\varepsilon }_{T}, \zeta } ) - g ( \mathbb{P}^{t}_{ {X}^{ \bar{\pi } }_{T}, \zeta } )
& = \int_{ \mathbb{R} \times \mathbb{R} } \nabla g ( \mathbb{P}^{t}_{ {X}^{ \bar{\pi } }_{T}, \zeta }, x,z ) 
                                          \big( \mathbb{P}^{t}_{ \bar{X}^{\varepsilon }_{T}, \zeta } ( dx, dz ) - \mathbb{P}^{t}_{ {X}^{ \bar{\pi } }_{T}, \zeta } ( dx, dz ) \big) + o ( {\varepsilon } ) \\
& = \mathbb{E}_{t} [ \nabla g ( \mathbb{P}^{t}_{ {X}^{ \bar{\pi } }_{T}, \zeta }, \bar{X}^{\varepsilon }_{T}, \zeta ) - \nabla g ( \mathbb{P}^{t}_{ {X}^{ \bar{\pi } }_{T}, \zeta }, {X}^{ \bar{\pi } }_{T}, \zeta ) ] 
  + o ( \varepsilon ) \\
& = \mathbb{E}_{t} \bigg[ \bigg( \zeta + \theta \Big( \lambda ( \mathbb{P}^{t}_{ {X}^{ \bar{\pi } }_{T}, \zeta } ) - {X}^{ \bar{\pi } }_{T} - \frac{\zeta }{\theta } \Big)_{+} \bigg) 
                          ( \bar{X}^{\varepsilon }_{T} - {X}^{ \bar{\pi } }_{T} ) \bigg] \\
& \quad 
  - \mathbb{E}_{t} \bigg[ \theta ( \bar{X}^{\varepsilon }_{T} - {X}^{ \bar{\pi } }_{T} )^{2} 
                          \int_{0}^{1} \rho {1}_{\{ \rho {X}^{ \bar{\pi } }_{T} + ( 1 - \rho ) \bar{X}^{\varepsilon }_{T} + \frac{\zeta }{\theta } \le \lambda ( \mathbb{P}^{t}_{ {X}^{ \bar{\pi } }_{T}, \zeta } ) \}} d \rho \bigg]
  + o ( \varepsilon ) \\
& = \xi \int_{t}^{ t + \varepsilon } \mathbb{E}_{t} [ ( {Y}^{t}_{s} {\vartheta }_{s} + \mathcal{Y}^{t}_{s} ) {\sigma }_{s} ] ds \\
& \quad
  - \frac{\theta }{2} {\xi }^{2} \mathbb{E}_{t} \bigg[ \bigg( \int_{t}^{ t + \varepsilon } {e}^{ \int_{s}^{T} {r}_{v} dv } {\sigma }_{s} d {W}_{s} \bigg)^{2} 
                                                       {1}_{\{ {X}^{ \bar{\pi } }_{T} + \frac{\zeta }{\theta } \le \lambda ( \mathbb{P}^{t}_{ {X}^{ \bar{\pi } }_{T}, \zeta } ) \}} \bigg]
  + o ( \varepsilon ),
\end{align*}
where the last equality arise from the following estimate 
\begin{align*}
& \bigg| \mathbb{E}_{t} \bigg[ ( \bar{X}^{\varepsilon }_{T} - {X}^{ \bar{\pi } }_{T} )^{2} 
                               \bigg( \int_{0}^{1} 2 \rho {1}_{\{ \rho {X}^{ \bar{\pi } }_{T} + ( 1 - \rho ) \bar{X}^{\varepsilon }_{T} + \frac{\zeta }{\theta } \le \lambda ( \mathbb{P}^{t}_{ {X}^{ \bar{\pi } }_{T}, \zeta } ) \}} d \rho
                                    - {1}_{\{ {X}^{ \bar{\pi } }_{T} + \frac{\zeta }{\theta } \le \lambda ( \mathbb{P}^{t}_{ {X}^{ \bar{\pi } }_{T}, \zeta } ) \}} \bigg) \bigg] \bigg|^{2} \\
& \le \mathbb{E}_{t} \big[ ( \bar{X}^{\varepsilon }_{T} - {X}^{ \bar{\pi } }_{T} )^{4} \big] \cdot
      \mathbb{E}_{t} \bigg[ \int_{0}^{1} 4 {\rho }^{2} \big| {1}_{\{ \rho {X}^{ \bar{\pi } }_{T} + ( 1 - \rho ) \bar{X}^{\varepsilon }_{T} + \frac{\zeta }{\theta } \le \lambda ( \mathbb{P}^{t}_{ {X}^{ \bar{\pi } }_{T}, \zeta } ) \}} 
                                                                  - {1}_{\{ {X}^{ \bar{\pi } }_{T} + \frac{\zeta }{\theta } \le \lambda ( \mathbb{P}^{t}_{ {X}^{ \bar{\pi } }_{T}, \zeta } ) \}} \big|^{2} d \rho \bigg] 
\end{align*}
with the Burkholder inequality
\begin{align*}
    \mathbb{E}_{t} [ | \bar{X}^{\varepsilon }_{T} - {X}^{ \bar{\pi } }_{T} |^{4} ]
& \le 8 {\xi }^{4} {e}^{ 4 T \esssup\limits_{ [ 0,T ] \times \Omega } |r| }
      \mathbb{E}_{t} \bigg[ \Big| \int_{t}^{ t + \varepsilon } {e}^{ - \int_{t}^{s} {r}_{v} dv } {\sigma }_{s} d {W}_{s} \Big|^{4} \bigg] 
    + 8 {\xi }^{4} \mathbb{E}_{t} \bigg[ \Big| \int_{t}^{ t + \varepsilon } {e}^{ \int_{s}^{T} {r}_{v} dv } {\sigma }_{s} {\vartheta }_{s} ds \Big|^{4} \bigg] \\
& \le \frac{ {2}^{21} }{ {3}^{6} } {\xi }^{4} {e}^{ 4 T \esssup\limits_{ [ 0,T ] \times \Omega } |r| } 
      \mathbb{E}_{t} \bigg[ \Big| \int_{t}^{ t + \varepsilon } {e}^{ - 2 \int_{t}^{s} {r}_{v} dv } | {\sigma }_{s} |^{2} ds \Big|^{2} \bigg] + O ( {\varepsilon }^{4} ) \\
& = O ( {\varepsilon }^{2} ).
\end{align*}
In addition, 
\begin{align*}
& \mathbb{E}_{t} \bigg[ \bigg| \bigg( \int_{t}^{ t + \varepsilon } {e}^{ - \int_{t}^{s} {r}_{v} dv } {\sigma }_{s} d {W}_{s} \bigg)^{2} - \bigg( \int_{t}^{ t + \varepsilon } {\sigma }_{s} d {W}_{s} \bigg)^{2} \bigg|
                               {e}^{ 2 \int_{t}^{T} {r}_{v} dv } {1}_{\{ {X}^{ \bar{\pi } }_{T} + \frac{\zeta }{\theta } \le \lambda ( \mathbb{P}^{t}_{ {X}^{ \bar{\pi } }_{T}, \zeta } ) \}} \bigg] \\
& \le {e}^{ 2 T \esssup\limits_{ [ 0,T ] \times \Omega } |r| }   
      \mathbb{E}_{t} \bigg[ \bigg| \int_{t}^{ t + \varepsilon } \Big( {e}^{ - \int_{t}^{s} {r}_{v} dv } + 1 \Big) {\sigma }_{s} d {W}_{s} \bigg| 
                            \bigg| \int_{t}^{ t + \varepsilon } \Big( {e}^{ - \int_{t}^{s} {r}_{v} dv } - 1 \Big) {\sigma }_{s} d {W}_{s} \bigg| \bigg] \\
& \le {e}^{ 2 T \esssup\limits_{ [ 0,T ] \times \Omega } |r| }   
      \bigg( \mathbb{E}_{t} \bigg[ \int_{t}^{ t + \varepsilon } \Big( {e}^{ - \int_{t}^{s} {r}_{v} dv } + 1 \Big)^{2} | {\sigma }_{s} |^{2} ds \bigg] \bigg)^{ \frac{1}{2} }
      \bigg( \mathbb{E}_{t} \bigg[ \int_{t}^{ t + \varepsilon } \Big( {e}^{ - \int_{t}^{s} {r}_{v} dv } - 1 \Big)^{2} | {\sigma }_{s} |^{2} ds \bigg] \bigg)^{ \frac{1}{2} } \\
& \le 2 \varepsilon {e}^{ 4 T \esssup\limits_{ [ 0,T ] \times \Omega } |r| } \bigg( \esssup\limits_{ [ 0,T ] \times \Omega } | \sigma |^{2} \bigg) 
      \bigg( \mathbb{E}_{t} \bigg[ \frac{1}{\varepsilon } \int_{t}^{ t + \varepsilon } \Big( {e}^{ - \int_{t}^{s} {r}_{v} dv } - 1 \Big)^{2} ds \bigg] \bigg)^{ \frac{1}{2} } \\
&   = o ( \varepsilon ).                  
\end{align*}
Given the relationship between $J ( t, \bar{\pi } )$ and $g ( \mathbb{P}^{t}_{ {X}^{ \bar{\pi } }_{T}, \zeta } )$, the asymptotic estimate \cref{eq: asymptotic estimate for ONEC :eq} follows.

\subsection{Proof of \texorpdfstring{\Cref{thm: verification theorem for ONEC}}{Theorem 3.3}}
\label{pf-thm: verification theorem for ONEC}

For the sake of brevity, we let $K$ represent a sufficiently large generic constant that only relies on $( T, \theta, \delta, p )$ and the essential supremum of $( r, \sigma, \vartheta )$, 
which usually serves for H\"older's inequality, Minkowski's inequality, Burkholder-Davis-Gundy inequality and a priori estimate for $\mathbb{L}^{p}$-solution of BSDE (see \cite{Briand-Delyon-Hu-Pardoux-Stoica-2003}) and may differ at different places.
To formally ``separate'' the variables $( {X}^{ \bar{\pi } }_{T}, \lambda )$ in the terminal value of \cref{eq: BSDE :eq}, we write 
$( \lambda - x - \frac{\zeta }{\theta } )_{+} = \int_{ - \infty }^{\lambda } {1}_{\{ y \ge x + \frac{\zeta }{\theta } \}} dy$.
Let us introduce the martingale representation
\begin{equation*}
\mathbb{E}_{t} [ {1}_{\{ y \ge {X}^{ \bar{\pi } }_{T} + \frac{\zeta }{\theta } \}} ] 
= \mathbb{E} [ {1}_{\{ y \ge {X}^{ \bar{\pi } }_{T} + \frac{\zeta }{\theta } \}} ] + \int_{0}^{t} \eta ( s,y ) d {W}_{s}, \\
\end{equation*}
where $\eta ( \cdot, y ) \in \mathbb{L}^{2}_{\mathbb{F}} ( 0,T ; \mathbb{L}^{p} ( \Omega; \mathbb{R} ) )$ for any $p > 1$, leading to
\begin{equation*}
  \mathbb{E}_{t} \Big[ \Big( \lambda - {X}^{ \bar{\pi } }_{T} - \frac{\zeta }{\theta } \Big)_{+} \Big] 
= \mathbb{E} \Big[ \Big( \lambda - {X}^{ \bar{\pi } }_{T} - \frac{\zeta }{\theta } \Big)_{+} \Big] 
+ \int_{ - \infty }^{\lambda } dy \int_{0}^{t} \eta ( s,y ) d {W}_{s}, \quad \forall \lambda \in \mathbb{R}. 
\end{equation*}
Notably, for any fixed $\lambda \in \mathbb{R}$,
\begin{align*}
  \bigg| \mathbb{E} \bigg[ \int_{0}^{\lambda } \bigg( \int_{0}^{T} | \eta ( s,y ) |^{2} ds \bigg)^{ \frac{1}{2} } dy \bigg] \bigg|
& \le \bigg| \int_{0}^{\lambda } \bigg( \mathbb{E} \bigg[ \int_{0}^{T} | \eta ( s,y ) |^{2} ds \bigg] \bigg)^{ \frac{1}{2} } dy \bigg| \\
& \le \bigg| \int_{0}^{\lambda } \bigg( \mathbb{E} \bigg[ \Big| \int_{0}^{T} \eta ( s,y ) d {W}_{s} \Big|^{2} \bigg] \bigg)^{ \frac{1}{2} } dy \bigg| \\
& \le K | \lambda |^{ \frac{1}{2} }
      \bigg( \int_{ - \infty }^{ | \lambda | } \mathbb{E} [ {1}_{\{ y \ge {X}^{ \bar{\pi } }_{T} + \frac{\zeta }{\theta } \}} ] dy \bigg)^{ \frac{1}{2} } \\
&   = K | \lambda |^{ \frac{1}{2} }
      \bigg( \mathbb{E} \bigg[ \Big( | \lambda | - {X}^{ \bar{\pi } }_{T} - \frac{\zeta }{\theta } \Big)_{+} \bigg] \bigg)^{ \frac{1}{2} } \\
& \le K | \lambda |^{ \frac{1}{2} } \bigg( | \lambda | + \mathbb{E} \Big[ \Big| {X}^{ \bar{\pi } }_{T} + \frac{\zeta }{\theta } \Big| \Big] \bigg)^{ \frac{1}{2} } < \infty,
\end{align*}
and for any fixed $p \in ( 1,2 )$,
\begin{align*}
  \mathbb{E} \bigg[ \int_{ - \infty }^{-1} \bigg( \int_{0}^{T} | \eta ( s,y ) |^{2} ds \bigg)^{ \frac{1}{2} } dy \bigg]
& \le \int_{ - \infty }^{-1} \bigg( \mathbb{E} \bigg[ \Big( \int_{0}^{T} | \eta ( s,y ) |^{2} ds \Big)^{ \frac{p}{2} } \bigg] \bigg)^{ \frac{1}{p} } dy \\
& \le K \int_{ - \infty }^{-1} ( \mathbb{E} [ {1}_{\{ y \ge {X}^{ \bar{\pi } }_{T} + \frac{\zeta }{\theta } \}} ] )^{ \frac{1}{p} } dy \\
& \le K \int_{1}^{\infty } {y}^{ - \frac{2}{p} } ( {y}^{2} \mathbb{E} [ {1}_{\{ {X}^{ \bar{\pi } }_{T} + \frac{\zeta }{\theta } \le -y \}} ] )^{ \frac{1}{p} } dy < \infty,
\end{align*}
where the last inequality is because $\int_{1}^{\infty } {y}^{ - \frac{2}{p} } dy = \frac{p}{2-p} < - \infty $ and 
\begin{equation*}
    {y}^{2} \mathbb{E} [ {1}_{\{ {X}^{ \bar{\pi } }_{T} + \frac{\zeta }{\theta } \le - y \}} ]
\le \int_{0}^{ {y}^{2} } \mathbb{E} [ {1}_{\{ | {X}^{ \bar{\pi } }_{T} + \frac{\zeta }{\theta } |^{2} \ge {y}^{2} \ge z \}} ] dz
\le \int_{0}^{\infty } \mathbb{E} [ {1}_{\{ | {X}^{ \bar{\pi } }_{T} + \frac{\zeta }{\theta } |^{2} \ge z \}} ] dz
  = \mathbb{E} \Big[ \Big| {X}^{ \bar{\pi } }_{T} + \frac{\zeta }{\theta } \Big|^{2} \Big] < \infty.
\end{equation*}
As a result, $\mathbb{E} [ \int_{ - \infty }^{\lambda } ( \int_{0}^{T} | \eta ( s,y ) |^{2} ds )^{ \frac{1}{2} } dy ] < \infty$, 
implying that $\int_{ - \infty }^{\lambda } ( \int_{0}^{T} | \eta ( s,y ) |^{2} ds )^{ \frac{1}{2} } dy < \infty$, $\mathbb{P}$-a.s., 
and hence one can refer to the stochastic Fubini theorem (see \cite{Veraar-2012}) for the interchange of Lebesgue and It\^o integrals
\begin{equation}
\label{eq: interchange of integrals :eq}
\int_{ - \infty }^{\lambda } dy \int_{0}^{t} \eta ( s,y ) d {W}_{s} = \int_{0}^{t} \bigg( \int_{ - \infty }^{\lambda } \eta ( s,y ) dy \bigg) d {W}_{s}, \quad \forall t \in [ 0,T ], ~ \mathbb{P}-a.s.
\end{equation}
Therefore, we obtain the martingale representation
\begin{equation*}
  \mathbb{E}_{t} \Big[ \Big( \lambda - {X}^{ \bar{\pi } }_{T} - \frac{\zeta }{\theta } \Big)_{+} \Big] 
= \mathbb{E} \Big[ \Big( \lambda - {X}^{ \bar{\pi } }_{T} - \frac{\zeta }{\theta } \Big)_{+} \Big] 
+ \int_{0}^{t} \bigg( \int_{ - \infty }^{\lambda } \eta ( s,y ) dy \bigg) d {W}_{s}, \quad \forall \lambda \in \mathbb{R}. 
\end{equation*}

For $\lambda \in \mathbb{L}^{1}_{ \mathcal{F}_{T} } ( \Omega; \mathbb{R} )$, by using the law of iterated conditioning $\mathbb{E} [ \cdot ] = \mathbb{E} [ \mathbb{E} [ \cdot | \lambda ] ]$, one obtains
\begin{align*}
&   \mathbb{E} \bigg[ \int_{0}^{T} ds \int_{ - \infty }^{\lambda } | \eta ( s,y ) | dy \bigg] 
\le K \mathbb{E} \bigg[ \int_{ - \infty }^{\lambda } \Big( \int_{0}^{T} | \eta ( s,y ) |^{2} ds \Big)^{ \frac{1}{2} } dy \bigg] < \infty, \\
&   \mathbb{E} \bigg[ \int_{0}^{T} ds \int_{ - \infty }^{\lambda } | \eta ( s,y ) |^{2} dy \bigg] 
  = \mathbb{E} \bigg[ \int_{ - \infty }^{\lambda } dy \int_{0}^{T} | \eta ( s,y ) |^{2} ds \bigg] 
\le K \mathbb{E} \Big[ \Big( \lambda - {X}^{ \bar{\pi } }_{T} - \frac{\zeta }{\theta } \Big)_{+} \Big] < \infty.
\end{align*}
Furthermore, with writing ${\lambda }_{t} = \lambda ( \mathbb{P}^{t}_{ {X}^{ \bar{\pi } }_{T}, \zeta } )$ for notational simplicity, we have
\begin{align*}
0 & = \mathbb{E}_{\tau } \Big[ \zeta + \theta \Big( {\lambda }_{\tau } - {X}^{ \bar{\pi } }_{T} - \frac{\zeta }{\theta } \Big)_{+} \Big]
    - \mathbb{E}_{t} \Big[ \zeta + \theta \Big( {\lambda }_{t} - {X}^{ \bar{\pi } }_{T} - \frac{\zeta }{\theta } \Big)_{+} \Big] \\
  & = \int_{t}^{\tau } \bigg( {\eta }_{s} + \theta \int_{ - \infty }^{ {\lambda }_{t} } \eta ( s,y ) dy \bigg) d {W}_{s}
    + \theta \int_{ {\lambda }_{t} }^{ {\lambda }_{\tau } } \mathbb{E}_{\tau } [ {1}_{\{ y \ge {X}^{ \bar{\pi } }_{T} + \frac{\zeta }{\theta } \}} ] dy,
\end{align*}
implying the right-continuity of $\lambda $.
Notably, it is not easy to derive a general moment estimate for ${\lambda }_{\tau } - {\lambda }_{t}$ in general, because
\begin{align*}
  \mathbb{E}_{t} \bigg[ \sup_{ \tau \in [ t, t + \varepsilon ] } | {\lambda }_{\tau } - {\lambda }_{t} |^{p} \bigg]
& \ge \mathbb{E}_{t} \bigg[ \sup_{ \tau \in [ t, t + \varepsilon ] } 
                            \bigg| \int_{ {\lambda }_{t} }^{ {\lambda }_{\tau } } \mathbb{E}_{\tau } [ {1}_{\{ y \ge {X}^{ \bar{\pi } }_{T} + \frac{\zeta }{\theta } \}} ] dy \bigg|^{p} \bigg] \\
&   = \mathbb{E}_{t} \bigg[ \sup_{ \tau \in [ t, t + \varepsilon ] } 
                            \bigg| \int_{t}^{\tau } \bigg( {\eta }_{s} + \theta \int_{ - \infty }^{ {\lambda }_{t} } \eta ( s,y ) dy \bigg) d {W}_{s} \bigg|^{p} \bigg] \\
& \ge K \mathbb{E}_{t} \bigg[ \bigg( \int_{t}^{ t + \varepsilon } \bigg| {\eta }_{s} + \theta \int_{ - \infty }^{ {\lambda }_{t} } \eta ( s,y ) dy \bigg|^{2} ds \bigg)^{ \frac{p}{2} } \bigg],
\quad \forall p > 0.
\end{align*}

Now we show the sufficiency of \cref{eq: equilibrium condition for ONEC :eq} for \cref{eq: primal equilibrium condition for ONEC :eq}.
Let us introduce the following martingale representations:
\begin{align*}
  \mathbb{E}_{t} \Big[ {e}^{ \int_{0}^{T} {r}_{v} dv } \zeta \Big] 
& = \mathbb{E} \Big[ {e}^{ \int_{0}^{T} {r}_{v} dv } \zeta \Big] + \int_{0}^{t} \bar{\eta }_{s} d {W}_{s} \quad \text{and} \\
  \mathbb{E}_{t} \Big[ {e}^{ \int_{0}^{T} {r}_{v} dv } {1}_{\{ y \ge {X}^{ \bar{\pi } }_{T} + \frac{\zeta }{\theta } \}} \Big] 
& = \mathbb{E} \Big[ {e}^{ \int_{0}^{T} {r}_{v} dv } {1}_{\{ y \ge {X}^{ \bar{\pi } }_{T} + \frac{\zeta }{\theta } \}} \Big] + \int_{0}^{t} \hat{\eta } ( s,y ) d {W}_{s},
\end{align*}
where $\bar{\eta }, \hat{\eta } ( \cdot, y ) \in \mathbb{L}^{2}_{\mathbb{F}} ( 0,T ; \mathbb{L}^{p} ( \Omega; \mathbb{R} ) )$ for any $p > 1$.
Given the boundedness of $r$, mirroring the previous discussions on ${\eta }^{(y)}$, 
we can show that 
\begin{equation*}
\int_{ - \infty }^{\lambda } dy \int_{0}^{t} \hat{\eta } ( s,y ) d {W}_{v} = \int_{0}^{t} \bigg( \int_{ - \infty }^{\lambda } \hat{\eta } ( s,y ) dy \bigg) d {W}_{s}, \quad 
\forall t \in [ 0,T ], ~ \mathbb{P}-a.s., \quad \forall \lambda \in \mathbb{R};
\end{equation*}
and
\begin{equation*}
\mathbb{E} \bigg[ \int_{0}^{T} ds \int_{ - \infty }^{\lambda } | \hat{\eta } ( s,y ) | dy \bigg] < \infty, \quad
\mathbb{E} \bigg[ \int_{0}^{T} ds \int_{ - \infty }^{\lambda } | \hat{\eta } ( s,y ) |^{2} dy \bigg] < \infty, \quad 
\forall \lambda \in \mathbb{L}^{1}_{ \mathcal{F}_{T} } ( \Omega; \mathbb{R} ).
\end{equation*}
Consequently, for the pair $( {Y}^{t}, \mathcal{Y}^{t} ) = ( Y ( {\lambda }_{t} ), \mathcal{Y} ( {\lambda }_{t} ) )$ given by the BSDE \cref{eq: BSDE :eq}, one obtains
\begin{equation*}
\left\{ \begin{aligned}
{Y}^{t}_{s} & = \mathbb{E}_{s} \Big[ {e}^{ \int_{s}^{T} {r}_{v} dv } \zeta \Big]
              + \theta \int_{ - \infty }^{ {\lambda }_{t} } \mathbb{E}_{s} \Big[ {e}^{ \int_{s}^{T} {r}_{v} dv } {1}_{\{ y \ge {X}^{ \bar{\pi } }_{T} + \frac{\zeta }{\theta } \}} \Big] dy, \\
\mathcal{Y}^{t}_{s} & = {e}^{ - \int_{0}^{s} {r}_{v} dv } \bigg( \bar{\eta }_{s} + \theta \int_{ - \infty }^{ {\lambda }_{t} } \hat{\eta } ( s,y ) dy \bigg),
\end{aligned} \right.
\end{equation*}
$\mathbb{P}$-a.s. for a.e. $s \in [ 0,T ]$.
Therefore, for any $\varepsilon \in ( 0, T-t ]$,
\begin{align*}
& \bigg| \int_{t}^{ t + \varepsilon } \mathbb{E}_{t} [ ( {Y}^{t}_{s} {\vartheta }_{s} + \mathcal{Y}^{t}_{s} ) {\sigma }_{s} ] ds 
       - \int_{t}^{ t + \varepsilon } \mathbb{E}_{t} [ ( {Y}^{s}_{s} {\vartheta }_{s} + \mathcal{Y}^{s}_{s} ) {\sigma }_{s} ] ds \bigg| \\
& \le K \bigg( \int_{t}^{ t + \varepsilon } \mathbb{E}_{t} [ | {Y}^{s}_{s} - {Y}^{t}_{s} | ] ds
             + \int_{t}^{ t + \varepsilon } \mathbb{E}_{t} [ | \mathcal{Y}^{s}_{s} - \mathcal{Y}^{t}_{s} | ] ds \bigg) \\
&   = K \bigg( \int_{t}^{ t + \varepsilon } \mathbb{E}_{t} \bigg[ \bigg| \int_{ {\lambda }_{t} }^{ {\lambda }_{s} } 
                                                                         {e}^{ \int_{s}^{T} {r}_{v} dv } {1}_{\{ y \ge {X}^{ \bar{\pi } }_{T} + \frac{\zeta }{\theta } \}} dy \bigg| \bigg] ds 
             + \int_{t}^{ t + \varepsilon } \mathbb{E}_{t} \bigg[ \bigg| \int_{ {\lambda }_{t} }^{ {\lambda }_{s} } \hat{\eta } ( s,y ) dy \bigg| \bigg] ds \bigg) \\
& \le K \bigg( \int_{t}^{ t + \varepsilon } \mathbb{E}_{t} [ | {\lambda }_{s} - {\lambda }_{t} | ] ds
             + \int_{t}^{ t + \varepsilon } \mathbb{E}_{t} \bigg[ | {\lambda }_{s} - {\lambda }_{t} |^{ \frac{1}{2} }
                                                                  \Big( \int_{ - \infty }^{ \sup_{ \tau \in [ 0,s ] } | {\lambda }_{\tau } | } | \hat{\eta } ( s,y ) |^{2} dy \Big)^{ \frac{1}{2} } \bigg] ds \bigg) \\
& \le \varepsilon K \bigg( \frac{1}{\varepsilon } \int_{t}^{ t + \varepsilon } \mathbb{E}_{t} [ | {\lambda }_{s} - {\lambda }_{t} | ] ds \\
& \qquad \qquad          + \bigg| \frac{1}{\varepsilon } \int_{t}^{ t + \varepsilon } \mathbb{E}_{t} [ | {\lambda }_{s} - {\lambda }_{t} | ] ds \bigg|^{ \frac{1}{2} } \cdot
                           \bigg| \mathbb{E}_{t} \bigg[ \frac{1}{\varepsilon } \int_{t}^{ t + \varepsilon } 
                                                        \Big( \int_{ - \infty }^{ \sup_{ \tau \in [ 0,s ] } | {\lambda }_{\tau } | } | \hat{\eta } ( s,y ) |^{2} dy \Big) ds \bigg] \bigg|^{ \frac{1}{2} } \bigg).
\end{align*}
Sending $\varepsilon $ to $0$, we obtain
\begin{equation*}
\frac{1}{\varepsilon } \int_{t}^{ t + \varepsilon } \mathbb{E}_{t} [ | {\lambda }_{s} - {\lambda }_{t} | ] ds \to 0, \quad \mathbb{P}-a.s., ~ \forall t \in [ 0,T ],
\end{equation*}
due to $\lim_{ \varepsilon \downarrow 0 }\mathbb{E}_{t} [ | {\lambda }_{s} - {\lambda }_{t} | ] = 0$;
and
\begin{equation*}
\frac{1}{\varepsilon } \int_{t}^{ t + \varepsilon } \Big( \int_{ - \infty }^{ \sup_{ \tau \in [ 0,s ] } | {\lambda }_{\tau } | } | \hat{\eta } ( s,y ) |^{2} dy \Big) ds
\to \int_{ - \infty }^{ \sup_{ \tau \in [ 0,t ] } | {\lambda }_{\tau } | } | \hat{\eta } ( t,y ) |^{2} dy
  < \infty, \quad a.e. ~ t \in [ 0,T ], ~ \mathbb{P}-a.s.,
\end{equation*}
due to Lebesgue differentiation theorem with the integrality of $\int_{ - \infty }^{ \sup_{ \tau \in [ 0,t ] } | {\lambda }_{\tau } | } | \hat{\eta } ( t,y ) |^{2} dy$ 
(see \Cref{lem: moment estimate for lambda} for the integrability of $\sup_{ \tau \in [ 0,T ] } | {\lambda }_{\tau } |$).
Summing up, we arrive at
\begin{equation*}
  \lim_{ \varepsilon \downarrow 0 } \frac{1}{\varepsilon } \int_{t}^{ t + \varepsilon } \mathbb{E}_{t} [ ( {Y}^{t}_{s} {\vartheta }_{s} + \mathcal{Y}^{t}_{s} ) {\sigma }_{s} ] ds
= \lim_{ \varepsilon \downarrow 0 } \frac{1}{\varepsilon } \int_{t}^{ t + \varepsilon } \mathbb{E}_{t} [ ( {Y}^{s}_{s} {\vartheta }_{s} + \mathcal{Y}^{s}_{s} ) {\sigma }_{s} ] ds,
\quad a.e. ~ t \in [ 0,T ], ~ \mathbb{P}-a.s.,
\end{equation*}
with the moment estimate
\begin{align*}
& \mathbb{E} \bigg[ \int_{0}^{T} | ( {Y}^{s}_{s} {\vartheta }_{s} + \mathcal{Y}^{s}_{s} ) {\sigma }_{s} | ds \bigg] \\
& \le K \mathbb{E} \bigg[ \int_{0}^{T} ( | {Y}^{s}_{s} | + | \mathcal {Y}^{s}_{s} | ) ds \bigg] \\
& \le K \mathbb{E} \bigg[ 1 + \int_{0}^{T} \mathbb{E}_{s} \Big[ \Big( {\lambda }_{s} - {X}^{ \bar{\pi } }_{T} - \frac{\zeta }{\theta } \Big)_{+} \Big] ds
                            + \int_{0}^{T} | \bar{\eta }_{s} | ds + \int_{0}^{T} ds \int_{ - \infty }^{ {\lambda }_{s} } | \hat{\eta } ( s,y ) | dy \bigg] < \infty.
\end{align*}
Hence, the sufficiency of \cref{eq: equilibrium condition for ONEC :eq} for \cref{eq: primal equilibrium condition for ONEC :eq} immediately follows;
and thus, $\bar{\pi }$ is an ONEC if \cref{eq: equilibrium condition for ONEC :eq} holds.

Conversely, we assume that $\bar{\pi } \in \mathbb{L}^{2}_{\mathbb{F}} ( 0,T ; \mathbb{L}^{2} ( \Omega; \mathbb{R} ) )$ is an ONEC. 
For notational simplicity, we introduce the martingale representation
\begin{equation*}
\mathbb{E}_{t} \Big[ {e}^{ 2 \int_{0}^{T} {r}_{v} dv } {1}_{\{ y \ge {X}^{ \bar{\pi } }_{T} + \frac{\zeta }{\theta } \}} \Big] 
= \mathbb{E} \Big[ {e}^{ 2 \int_{0}^{T} {r}_{v} dv } {1}_{\{ y \ge {X}^{ \bar{\pi } }_{T} + \frac{\zeta }{\theta } \}} \Big] + \int_{0}^{t} \tilde{\eta } ( s,y ) d {W}_{s}.
\end{equation*}
Then, 
\begin{align*}
& \mathbb{E}_{t} \bigg[ \bigg( \int_{t}^{ t + \varepsilon } {\sigma }_{s} d {W}_{s} \bigg)^{2} {e}^{ 2 \int_{t}^{T} {r}_{v} dv } {1}_{\{ {X}^{ \bar{\pi } }_{T} + \frac{\zeta }{\theta } \le {\lambda }_{t} \}} \bigg] \\
& = \mathbb{E}_{t} \bigg[ \int_{t}^{ t + \varepsilon } | {\sigma }_{s} |^{2} ds \bigg]
    \mathbb{E}_{t} \Big[ {e}^{ 2 \int_{t}^{T} {r}_{v} dv } {1}_{\{ {X}^{ \bar{\pi } }_{T} + \frac{\zeta }{\theta } \le {\lambda }_{t} \}} \Big]
  + \mathbb{E}_{t} \bigg[ \bigg( \int_{t}^{ t + \varepsilon } {\sigma }_{s} d {W}_{s} \bigg)^{2} \int_{t}^{ t + \varepsilon } \tilde{\eta } ( s,y ) d {W}_{s} \bigg] \\
& = \mathbb{E}_{t} \bigg[ \int_{t}^{ t + \varepsilon } | {\sigma }_{s} |^{2} {e}^{ 2 \int_{s}^{T} {r}_{v} dv } {1}_{\{ {X}^{ \bar{\pi } }_{T} + \frac{\zeta }{\theta } \le {\lambda }_{s} \}} ds \bigg]
  + o ( \varepsilon ),
\end{align*}
because of
\begin{align*}
  \bigg| \mathbb{E}_{t} \bigg[ \bigg( \int_{t}^{ t + \varepsilon } {\sigma }_{s} d {W}_{s} \bigg)^{2} \int_{t}^{ t + \varepsilon } \tilde{\eta } ( s,y ) d {W}_{s} \bigg] \bigg|^{2}
& \le \mathbb{E}_{t} \bigg[ \bigg( \int_{t}^{ t + \varepsilon } {\sigma }_{s} d {W}_{s} \bigg)^{4} \bigg]
      \mathbb{E}_{t} \bigg[ \bigg( \int_{t}^{ t + \varepsilon } \tilde{\eta } ( s,y ) d {W}_{s} \bigg)^{2} \bigg] \\
& \le K \mathbb{E}_{t} \bigg[ \bigg( \int_{t}^{ t + \varepsilon } | {\sigma }_{s} |^{2} ds \bigg)^{2} \bigg]
        \mathbb{E}_{t} \bigg[ \int_{t}^{ t + \varepsilon } | \tilde{\eta } ( s,y ) |^{2} ds \bigg] \\
& \le K {\varepsilon }^{2} \mathbb{E}_{t} \bigg[ \int_{t}^{ t + \varepsilon } | \tilde{\eta } ( s,y ) |^{2} ds \bigg]
    = o ( {\varepsilon }^{2} ),
\end{align*}
\begin{align*}
& \bigg| \mathbb{E}_{t} \bigg[ \int_{t}^{ t + \varepsilon } | {\sigma }_{s} |^{2} {e}^{ 2 \int_{s}^{T} {r}_{v} dv } {1}_{\{ {X}^{ \bar{\pi } }_{T} + \frac{\zeta }{\theta } \le {\lambda }_{t} \}} ds \bigg]
       - \mathbb{E}_{t} \bigg[ \int_{t}^{ t + \varepsilon } | {\sigma }_{s} |^{2} ds \bigg]
         \mathbb{E}_{t} \Big[ {e}^{ 2 \int_{t}^{T} {r}_{v} dv } {1}_{\{ {X}^{ \bar{\pi } }_{T} + \frac{\zeta }{\theta } \le {\lambda }_{t} \}} \Big] \bigg| \\
& \le \mathbb{E}_{t} \bigg[ \int_{t}^{ t + \varepsilon } | {\sigma }_{s} |^{2} {e}^{ 2 \int_{s}^{T} {r}_{v} dv } \Big( 1 - {e}^{ 2 \int_{t}^{s} {r}_{v} dv } \Big) {1}_{\{ {X}^{ \bar{\pi } }_{T} + \frac{\zeta }{\theta } \le {\lambda }_{t} \}} ds \bigg] \\
& \quad + \mathbb{E}_{t} \bigg[ \int_{t}^{ t + \varepsilon } | {\sigma }_{s} |^{2} 
                                \bigg| \mathbb{E}_{s} \Big[ {e}^{ 2 \int_{t}^{T} {r}_{v} dv } {1}_{\{ {X}^{ \bar{\pi } }_{T} + \frac{\zeta }{\theta } \le {\lambda }_{t} \}} \Big]
                                     - \mathbb{E}_{t} \Big[ {e}^{ 2 \int_{t}^{T} {r}_{v} dv } {1}_{\{ {X}^{ \bar{\pi } }_{T} + \frac{\zeta }{\theta } \le {\lambda }_{t} \}} \Big] \bigg| ds \bigg]
    = o ( \varepsilon )
\end{align*}
and the following estimate arising from the right-continuity of $\lambda $ and the dominated convergence theorem:
\begin{equation*}
\mathbb{E}_{t} \bigg[ \int_{t}^{ t + \varepsilon } | {\sigma }_{s} |^{2} {e}^{ 2 \int_{s}^{T} {r}_{v} dv } 
                                                   | {1}_{\{ {X}^{ \bar{\pi } }_{T} + \frac{\zeta }{\theta } \le {\lambda }_{s} \}} - {1}_{\{ {X}^{ \bar{\pi } }_{T} + \frac{\zeta }{\theta } \le {\lambda }_{t} \}} | ds \bigg]
\le K \int_{t}^{ t + \varepsilon } \mathbb{E}_{t} [ {1}_{\{ {\lambda }_{t} \wedge {\lambda }_{s} < {X}^{ \bar{\pi } }_{T} + \frac{\zeta }{\theta } \le {\lambda }_{t} \vee {\lambda }_{s} \}} ] ds 
  = o ( \varepsilon ).
\end{equation*}
Consequently, \cref{eq: asymptotic estimate for ONEC :eq} is equivalent to
\begin{align*}
  J ( t, \bar{\pi }^{ t, \varepsilon, \xi } ) - J ( t, \bar{\pi } )
& = \mathbb{E}_{t} \bigg[ \int_{t}^{ t + \varepsilon } \bigg( \xi ( {Y}^{s}_{s} {\vartheta }_{s} + \mathcal{Y}^{s}_{s} ) {\sigma }_{s}
                        - \frac{\theta }{2} {\xi }^{2} | {\sigma }_{s} |^{2} \mathbb{E}_{s} \Big[ {e}^{ 2 \int_{s}^{T} {r}_{v} dv } {1}_{\{ {X}^{ \bar{\pi } }_{T} + \frac{\zeta }{\theta } \le {\lambda }_{s} \}} \Big] \bigg) ds \bigg] 
  + o ( \varepsilon ).
\end{align*}
Notably, the integrand does not depend on $t$, and this is the key property for deriving the necessity of \cref{eq: equilibrium condition for ONEC :eq} for \cref{eq: primal equilibrium condition for ONEC :eq}.
In view of the arbitrariness of $\xi \in \mathbb{L}^{2}_{ \mathcal{F}_{t} } ( \Omega; \mathbb{R} )$, 
let us consider the simple function $\xi = z {1}_{A}$ with $z \in \mathbb{R}$ and $A \in \mathcal{F}_{t} \subseteq \mathcal{F}_{T}$ both arbitrarily fixed.
Since $\bar{\pi }$ is an ONEC, it follows from \Cref{def: open-loop Nash equilibrium control} that
\begin{align*}
0 & \ge \mathbb{E} \Bigg[ \limsup_{ \varepsilon \downarrow 0 } 
                          \mathbb{E}_{t} \bigg[ \frac{ {1}_{A} }{\varepsilon } 
                                                \int_{t}^{ t + \varepsilon } \bigg( z ( {Y}^{s}_{s} {\vartheta }_{s} + \mathcal{Y}^{s}_{s} ) {\sigma }_{s}
                                                                                  - \frac{\theta }{2} {z}^{2} | {\sigma }_{s} |^{2} 
                                                                                    \mathbb{E}_{s} \Big[ {e}^{ 2 \int_{s}^{T} {r}_{v} dv } {1}_{\{ {X}^{ \bar{\pi } }_{T} + \frac{\zeta }{\theta } \le {\lambda }_{s} \}} \Big] \bigg) ds \bigg] \Bigg] \\
  & \ge \limsup_{ \varepsilon \downarrow 0 } \frac{1}{\varepsilon } 
        \int_{t}^{ t + \varepsilon } \mathbb{E} \bigg[ {1}_{A} \bigg( z ( {Y}^{s}_{s} {\vartheta }_{s} + \mathcal{Y}^{s}_{s} ) {\sigma }_{s}
                                                                    - \frac{\theta }{2} {z}^{2} | {\sigma }_{s} |^{2} 
                                                                      \mathbb{E}_{s} \Big[ {e}^{ 2 \int_{s}^{T} {r}_{v} dv } {1}_{\{ {X}^{ \bar{\pi } }_{T} + \frac{\zeta }{\theta } \le {\lambda }_{s} \}} \Big] \bigg) \bigg] ds \\
  &   =  \mathbb{E} \bigg[ {1}_{A} \bigg( z ( {Y}^{t}_{t} {\vartheta }_{t} + \mathcal{Y}^{t}_{t} ) {\sigma }_{t}
                                        - \frac{\theta }{2} {z}^{2} | {\sigma }_{t} |^{2} \mathbb{E}_{t} \Big[ {e}^{ 2 \int_{t}^{T} {r}_{v} dv } {1}_{\{ {X}^{ \bar{\pi } }_{T} + \frac{\zeta }{\theta } \le {\lambda }_{t} \}} \Big] \bigg) \bigg],
\quad a.e. ~ t \in [ 0,T ].
\end{align*}
The arbitrariness of $A \in \mathcal{F}_{t}$ implies that
\begin{equation*}
0 \ge z ( {Y}^{t}_{t} {\vartheta }_{t} + \mathcal{Y}^{t}_{t} ) {\sigma }_{t}
    - \frac{\theta }{2} {z}^{2} | {\sigma }_{t} |^{2} \mathbb{E}_{t} \Big[ {e}^{ 2 \int_{t}^{T} {r}_{v} dv } {1}_{\{ {X}^{ \bar{\pi } }_{T} + \frac{\zeta }{\theta } \le {\lambda }_{t} \}} \Big], \quad \mathbb{P}-a.s., ~ a.e. ~ t \in [ 0,T ],
\end{equation*}
and then \cref{eq: equilibrium condition for ONEC :eq} arises due to the arbitrariness of $z \in \mathbb{R}$.
Thus, $\bar{\pi }$ is an ONEC only if \cref{eq: equilibrium condition for ONEC :eq} holds.
In conclusion, $\bar{\pi }$ is an ONEC if and only if \cref{eq: equilibrium condition for ONEC :eq} holds, and hence \cref{eq: equilibrium condition for ONEC :eq} is equivalent to \cref{eq: primal equilibrium condition for ONEC :eq}.

\subsection{Proof of \texorpdfstring{\Cref{lem: property of ONEC}}{Theorem 3.5}}
\label{pf-lem: property of ONEC}

Since $r$ is deterministic, from \Cref{pf-thm: verification theorem for ONEC} we know that
\begin{equation*}
\left\{ \begin{aligned}
{Y}^{t}_{s} & = {e}^{ \int_{s}^{T} {r}_{v} dv } \bigg( \mathbb{E}_{s} [ \zeta ]
                                                     + \theta \int_{ - \infty }^{ {\lambda }_{t} } P ( s,y ) dy \bigg), \\
\mathcal{Y}^{t}_{s} & = {e}^{ \int_{s}^{T} {r}_{v} dv } \bigg( {e}^{ - \int_{0}^{T} {r}_{v} dv } \bar{\eta }_{s} + \theta \int_{ - \infty }^{ {\lambda }_{t} } \eta ( s,y ) dy \bigg).
\end{aligned} \right.
\end{equation*}
Moreover, ${Y}^{t}_{t} = {e}^{ \int_{t}^{T} {r}_{v} dv }$ follows from the definition of ${\lambda }_{t}$, 
and hence the equilibrium condition \cref{eq: equilibrium condition for ONEC :eq} can be re-expressed as
\begin{equation*}
0 = {\vartheta }_{t} + {e}^{ - \int_{0}^{T} {r}_{v} dv } \bar{\eta }_{t} + \theta \int_{ - \infty }^{ {\lambda }_{t} } \eta ( t,y ) dy.
\end{equation*}
Now we consider the ansatz that $\{ {\lambda }_{t} \}_{ t \in [ 0,T ] }$ is a semi-martingale, since ${\lambda }_{t}$ is the unique solution of $1 = \mathbb{E}_{t} [ \zeta ] + \theta \int_{ - \infty }^{ {\lambda }_{t} } P ( t,y ) dy$.
Applying the It\^o-Kunita-Ventzel formula (see \cite[Theorem 1.5.3.2]{Jeanblanc-Yor-Chesney-2009}) to differentiate the both sides of $1 = \mathbb{E}_{t} [ \zeta ] + \theta \int_{ - \infty }^{ {\lambda }_{t} } P ( t,y ) dy$, 
with exploiting the semi-martingale decomposition \cref{eq: martingale representation for indicator :eq} and the interchange of integrals \cref{eq: interchange of integrals :eq}, yields
\begin{equation*}
0 = \bigg( {e}^{ - \int_{0}^{T} {r}_{v} dv } \bar{\eta }_{t} + \theta \int_{ - \infty }^{ {\lambda }_{t} } \eta ( t,y ) dy \bigg) d {W}_{t}
  + \theta P ( t, {\lambda }_{t} ) d {\lambda }_{t}
  + \theta \eta ( t, {\lambda }_{t} ) d \langle \lambda, W \rangle_{t} 
  + \frac{\theta }{2} {P}_{y} ( t, {\lambda }_{t} ) d \langle \lambda, \lambda \rangle_{t}.
\end{equation*}
Consequently, we obtain the following semi-martingale decomposition of $\{ {\lambda }_{t} \}_{ t \in [ 0,T ] }$:
\begin{align*}
  d {\lambda }_{t} 
& = \frac{ {\vartheta }_{t} }{ \theta P ( t, {\lambda }_{t} ) } d {W}_{t} 
  - \frac{ \eta ( t, {\lambda }_{t} ) }{ P ( t, {\lambda }_{t} ) } \frac{ {\vartheta }_{t} }{ \theta P ( t, {\lambda }_{t} ) } dt 
  - \frac{ {P}_{y} ( t, {\lambda }_{t} ) }{ 2 P ( t, {\lambda }_{t} ) } \bigg| \frac{ {\vartheta }_{t} }{ \theta P ( t, {\lambda }_{t} ) } \bigg|^{2} dt \\
& = \frac{ {\vartheta }_{t} }{ \theta \rho ( t,0 ) } d {W}_{t} 
  - \frac{ \varrho ( t,0 ) }{ \rho ( t,0 ) } d \langle \lambda, W \rangle_{t}
  + \frac{ {\rho }_{y} ( t,0 ) }{ 2 \rho ( t,0 ) } d \langle \lambda, \lambda \rangle_{t}.
\end{align*}
Likewise, differentiating $\rho ( t,y ) = P ( t, y + {\lambda }_{t} )$ w.r.t. $t$ yields
\begin{equation*}
d \rho ( t,y ) = \eta ( t, y + {\lambda }_{t} ) d {W}_{t} 
               + {\rho }_{y} ( t,y ) d {\lambda }_{t}
               + {\eta }_{y} ( t, y + {\lambda }_{t} ) d \langle \lambda, W \rangle_{t}
               + \frac{1}{2} {\rho }_{yy} ( t,y ) d \langle \lambda, \lambda \rangle_{t}.
\end{equation*} 
By rearrangement with $\rho ( T,y ) = {1}_{\{ y + {\lambda }_{T} \ge {X}^{ \bar{\pi } }_{T} + \frac{\zeta }{\theta } \}}$ and ${\lambda }_{T} = {X}^{ \bar{\pi } }_{T} + \frac{1}{\theta }$, 
we immediately obtain \cref{eq: BSPDE of rho :eq}, i.e.,
\begin{align*}               
\rho ( t,y ) & = {1}_{\{ y \ge - \frac{ 1 - \zeta }{\theta } \}}
               + \frac{1}{2} \int_{t}^{T} \bigg( {\rho }_{yy} ( v,y ) - {\rho }_{y} ( v,y ) \frac{ {\rho }_{y} ( v,0 ) }{ \rho ( v,0 ) } \bigg) d \langle \lambda, \lambda \rangle_{v} \\
             & \quad - \int_{t}^{T} \bigg( {\varrho }_{y} ( v,y ) - {\rho }_{y} ( v,y ) \frac{ \varrho ( v,0 ) }{ \rho ( v,0 ) } \bigg) d \langle \lambda, W \rangle_{v}
                     - \int_{t}^{T} \varrho ( v,y ) d {W}_{v}
\end{align*}

Conversely, if \cref{eq: BSPDE of rho :eq} admits a solution $( \rho, \varrho )$ satisfying the given conditions, 
then integrating the both sides of the above BSPDE on $y \in ( - \infty, 0 ]$ yields
\begin{equation*}
\int_{ - \infty }^{0} \rho ( t,y ) dy = \frac{ 1 - \zeta }{\theta } - \int_{t}^{T} \bigg( \int_{ - \infty }^{0} \varrho ( v,y ) dy \bigg) d {W}_{v},
\end{equation*}
implying that
\begin{equation*}
1 = \mathbb{E}_{t} [ \zeta ] + \theta \int_{ - \infty }^{0} \rho ( t,y ) dy, \quad
0 = {e}^{ - \int_{0}^{T} {r}_{v} dv } \bar{\eta }_{t} + \theta \int_{ - \infty }^{0} \varrho ( t,y ) dy.
\end{equation*}
With a slight abuse of notation, we consider the given hedging portfolio $\bar{\pi } \in \mathbb{L}^{2}_{\mathbb{F}} ( 0,T ; \mathbb{L}^{2} ( \Omega; \mathbb{R} ) )$ and
\begin{equation*}
\left\{ \begin{aligned}
\bar{Y}_{s} ( \lambda ) & = {e}^{ \int_{s}^{T} {r}_{v} dv } \bigg( \mathbb{E}_{s} [ \zeta ]
                                                                 + \theta \int_{ - \infty }^{ \lambda - \bar{\lambda }_{s} - c - \frac{1}{\theta } } \rho ( s,y ) dy \bigg), \\
\bar{\mathcal{Y}}_{s} ( \lambda ) & = {e}^{ \int_{s}^{T} {r}_{v} dv } \bigg( {e}^{ - \int_{0}^{T} {r}_{v} dv } \bar{\eta }_{s} 
                                                                           + \theta \int_{ - \infty }^{ \lambda - \bar{\lambda }_{s} - c - \frac{1}{\theta } } 
                                                                                    \Big( \varrho ( s,y ) - {\rho }_{y} ( s,y ) \frac{ {\vartheta }_{s} }{ \theta \rho ( s,0 ) } \Big) dy \bigg), \\
\bar{\lambda }_{t} & = \int_{0}^{t} \frac{ {\vartheta }_{s} }{ \theta \rho ( s,0 ) } d {W}_{s} 
                     + \int_{0}^{t} \bigg( \frac{ {\rho }_{y} ( s,0 ) }{ 2 \rho ( s,0 ) } \bigg| \frac{ {\vartheta }_{s} }{ \theta \rho ( s,0 ) } \bigg|^{2} 
                                         - \frac{ \varrho ( s,0 ) }{ \rho ( s,0 ) } \frac{ {\vartheta }_{s} }{ \theta \rho ( s,0 ) } \bigg) ds.
\end{aligned} \right.
\end{equation*}
One can find that $( \bar{Y} ( \lambda ), \bar{\mathcal{Y}} ( \lambda ) )$ solves the BSDE \cref{eq: BSDE :eq}, and hence ${\lambda }_{t} = \bar{\lambda }_{t} + c + \frac{1}{\theta }$.
Moreover, 
\begin{align*}
& \big( \bar{Y}_{t} ( {\lambda }_{t} ) {\vartheta }_{t} + \bar{\mathcal{Y}}_{t} ( {\lambda }_{t} ) \big) {e}^{ - \int_{t}^{T} {r}_{v} dv } \\
& = \bigg( \mathbb{E}_{t} [ \zeta ] + \theta \int_{ - \infty }^{0} \rho ( s,y ) dy \bigg) {\vartheta }_{t}
  + \bigg( {e}^{ - \int_{0}^{T} {r}_{v} dv } \bar{\eta }_{s} 
         + \theta \int_{ - \infty }^{0} \varrho ( s,y ) dy
         - \frac{ {\vartheta }_{s} }{ \rho ( s,0 ) } \int_{ - \infty }^{0} {\rho }_{y} ( s,y ) dy \bigg) \\
& = 0.
\end{align*}
Therefore, $\bar{\pi }$ is an ONEC according to \Cref{thm: verification theorem for ONEC}.

\subsection{Proof of \texorpdfstring{\Cref{thm: verification theorem for CNEC}}{Theorem 4.1}}
\label{pf-thm: verification theorem for CNEC}

Applying the It\^o-Kunita-Ventzel formula (see \cite[Theorem 1.5.3.2]{Jeanblanc-Yor-Chesney-2009}) to differentiate $\{ \mathcal{V} ( s, {X}^{ t,x, \Pi }_{s}, \lambda ) \}_{ s \in [ t,T ] }$
with the semi-martingale decomposition \cref{eq: BSPDE :eq} of $\mathcal{V} ( \cdot, x, \lambda )$ yields
\begin{equation*}
d \mathcal{V} ( s, {X}^{ t,x, \Pi }_{s}, \lambda ) = \big( \psi ( s, {X}^{ t,x, \Pi }_{s}, \lambda ) + \mathcal{V}_{x} ( s, {X}^{ t,x, \Pi }_{s}, \lambda ) \bar{\Pi } ( s, {X}^{ t,x, \Pi }_{s} ) {\sigma }_{s} \big) d {W}_{s},
\end{equation*}
which implies that $\mathcal{V} ( t,x, \lambda ) = \mathbb{E}_{t} [ \mathcal{V} ( T, {X}^{ t,x, \Pi }_{T}, \lambda ) ] = \tilde{J} ( t,x, \bar{\Pi }, \lambda )$.
In the same manner, 
\begin{equation*}
\mathcal{M} ( t,x, \lambda ) = \mathbb{E}_{t} [ \mathcal{M} ( T, {X}^{ t,x, \bar{\Pi } }_{T}, \lambda ) ] = \mathbb{E}_{t} \Big[ \zeta + \theta \Big( \lambda - {X}^{ t,x, \bar{\Pi } }_{T} - \frac{\zeta }{\theta } \Big)_{+} \Big],
\end{equation*}
and hence $\Lambda ( t,x ) = \lambda ( \mathbb{P}^{t}_{ {X}^{ t,x, \bar{\Pi } }_{T}, \zeta } )$.
Next, we proceed with the perturbation argument for an arbitrarily fixed $( t,x )$.
By applying iterated conditioning to \cref{eq: extended preference with closed-loop control :eq}, with ${X}^{ t,x, \bar{\Pi }^{ t, \varepsilon, \xi } }_{T} = {X}^{ t + \varepsilon, {X}^{ t,x, \xi }_{ t + \varepsilon }, \bar{\Pi } }_{T}$, we have
\begin{equation*}
  \bar{J} ( t,x, \bar{\Pi }^{ t, \varepsilon, \xi } ) 
= \tilde{J} ( t, x, \bar{\Pi }^{ t, \varepsilon, \xi }, {\lambda }_{\varepsilon } )
= \mathbb{E}_{t} [ \tilde{J} ( t + \varepsilon, {X}^{ t,x, \xi }_{ t + \varepsilon }, \bar{\Pi }, {\lambda }_{\varepsilon } ) ]
= \mathbb{E}_{t} [ \mathcal{V} ( t + \varepsilon, {X}^{ t,x, \xi }_{ t + \varepsilon }, {\lambda }_{\varepsilon } ) ],
\end{equation*}
where ${\lambda }_{\varepsilon } := \lambda ( \mathbb{P}^{t}_{ {X}^{ t,x, \bar{\Pi }^{ t, \varepsilon, \xi } }_{T}, \zeta } )$ with a slight abuse of notation.
So it is supposed to make an asymptotic estimation for ${\lambda }_{\varepsilon } - {\lambda }_{0}$.
By iterated conditioning,
\begin{equation*}
  \mathbb{E}_{t} [ \mathcal{M} ( t + \varepsilon, {X}^{ t,x, \xi }_{ t + \varepsilon }, {\lambda }_{\varepsilon } ) ]
= \mathbb{E}_{t} \Big[ \zeta + \theta \Big( {\lambda }_{\varepsilon } - {X}^{ t,x, \bar{\Pi }^{ t, \varepsilon, \xi } }_{T} - \frac{\zeta }{\theta } \Big)_{+} \Big] 
= 1
= \mathcal{M} ( t,x, {\lambda }_{0} ),
\end{equation*}
which implies that ${\lambda }_{\varepsilon } \to {\lambda }_{0}$ as $\varepsilon \downarrow 0$.
Furthermore, by It\^o-Kunita-Ventzel formula,
\begin{align*}
&   | \mathcal{M} ( t,x, {\lambda }_{0} ) - \mathcal{M} ( t,x, {\lambda }_{\varepsilon } ) | \\
& = | \mathbb{E}_{t} [ \mathcal{M} ( t + \varepsilon, {X}^{ t,x, \xi }_{ t + \varepsilon }, {\lambda }_{\varepsilon } ) - \mathcal{M} ( t,x, {\lambda }_{\varepsilon } ) ] | \\
& = \bigg| \mathbb{E}_{t} \bigg[ \int_{t}^{ t + \varepsilon } \bigg( \frac{1}{2} \mathcal{M}_{xx} ( s, {X}^{ t,x, \xi }_{s}, {\lambda }_{\varepsilon } ) \big( {\xi }^{2} - | \bar{\Pi } ( s, {X}^{ t,x, \xi }_{s} ) |^{2} \big) | {\sigma }_{s} |^{2} \\
& \qquad \qquad \qquad                                             + \big( \mathcal{M}_{x} ( s, {X}^{ t,x, \xi }_{s}, {\lambda }_{s} ) {\vartheta }_{s} + {\phi }_{x} ( s, {X}^{ t,x, \xi }_{s}, {\lambda }_{\varepsilon } ) \big) 
                                                                     \big( \xi - \bar{\Pi } ( s, {X}^{ t,x, \xi }_{s} ) \big) {\sigma }_{s} \bigg) ds \bigg] \bigg| \\
& \le \bigg| \frac{1}{2} \mathcal{M}_{xx} ( t,x, {\lambda }_{0} ) \big( {\xi }^{2} - | \bar{\Pi } ( t,x ) |^{2} \big) \\
& \qquad + \big( \mathcal{M}_{x} ( t,x, {\lambda }_{0} ) {\vartheta }_{t} + {\phi }_{x} ( t,x, {\lambda }_{0} ) \big) 
           \big( \xi - \bar{\Pi } ( t,x ) \big) \bigg| \bigg( \esssup_{ [ 0,T ] \times \Omega } | \sigma |^{2} \vee 1 \bigg) \varepsilon + o ( \varepsilon ).
\end{align*}
Notably, due to the convexity of $\mathcal{M} ( t,x, \cdot )$,
\begin{align*}
    ( {\lambda }_{\varepsilon } - {\lambda }_{0} ) \int_{\{ y + \frac{z}{\theta } < {\lambda }_{0} \}} \mathbb{P}^{t}_{ {X}^{ t,x, \bar{\Pi } }_{T}, \zeta } ( dy, dz )
& \le \mathcal{M} ( t,x, {\lambda }_{\varepsilon } ) - \mathcal{M} ( t,x, {\lambda }_{0} ) \\
& \le ( {\lambda }_{\varepsilon } - {\lambda }_{0} ) \int_{\{ y + \frac{z}{\theta } < {\lambda }_{\varepsilon } \}} \mathbb{P}^{t}_{ {X}^{ t,x, \bar{\Pi } }_{T}, \zeta } ( dy, dz ).
\end{align*}
Summing up, we obtain ${\lambda }_{\varepsilon } - {\lambda }_{0} = O ( \varepsilon )$.
Then, by It\^o-Kunita-Ventzel formula and Taylor expansion,
\begin{align*}
&   \bar{J} ( t,x, \bar{\Pi }^{ t, \varepsilon, \xi } ) - \bar{J} ( t,x, \bar{\Pi } ) \\
& = \mathbb{E}_{t} [ \mathcal{V} ( t + \varepsilon, {X}^{ t,x, \xi }_{ t + \varepsilon }, {\lambda }_{\varepsilon } ) ] - \mathcal{V} ( t,x, {\lambda }_{0} ) \\
& = \mathbb{E}_{t} [ \mathcal{V} ( t + \varepsilon, {X}^{ t,x, \xi }_{ t + \varepsilon }, {\lambda }_{\varepsilon } ) - \mathcal{V} ( t,x, {\lambda }_{\varepsilon } ) ] 
  + \big( \mathcal{V} ( t,x, {\lambda }_{\varepsilon } ) - \mathcal{V} ( t,x, {\lambda }_{0} ) \big) \\
& = \mathbb{E}_{t} \bigg[ \int_{t}^{ t + \varepsilon } \bigg( \frac{1}{2} \mathcal{V}_{xx} ( s, {X}^{ t,x, \xi }_{s}, {\lambda }_{\varepsilon } ) \big( {\xi }^{2} - | \bar{\Pi } ( s, {X}^{ t,x, \xi }_{s} ) |^{2} \big) | {\sigma }_{s} |^{2} \\
& \qquad \qquad \qquad                                      + \big( \mathcal{V}_{x} ( s, {X}^{ t,x, \xi }_{s}, {\lambda }_{\varepsilon } ) {\vartheta }_{s} + {\psi }_{x} ( s, {X}^{ t,x, \xi }_{s}, {\lambda }_{\varepsilon } ) \big) 
                                                              \big( \xi - \bar{\Pi } ( s, {X}^{ t,x, \xi }_{s} ) \big) {\sigma }_{s} \bigg) ds \bigg] \\
& \quad + \mathcal{V}_{\lambda } ( t,x, {\lambda }_{0} ) ( {\lambda }_{\varepsilon } - {\lambda }_{0} ) + o ( \varepsilon ) \\
& = \mathbb{E}_{t} \bigg[ \int_{t}^{ t + \varepsilon } \bigg( \frac{1}{2} \mathcal{V}_{xx} ( s, {X}^{ t,x, \xi }_{s}, {\lambda }_{0} ) \big( {\xi }^{2} - | \bar{\Pi } ( s, {X}^{ t,x, \xi }_{s} ) |^{2} \big) | {\sigma }_{s} |^{2} \\
& \qquad \qquad \qquad                                      + \big( \mathcal{V}_{x} ( s, {X}^{ t,x, \xi }_{s}, {\lambda }_{0} ) {\vartheta }_{s} + {\psi }_{x} ( s, {X}^{ t,x, \xi }_{s}, {\lambda }_{0} ) \big) 
                                                              \big( \xi - \bar{\Pi } ( s, {X}^{ t,x, \xi }_{s} ) \big) {\sigma }_{s} \bigg) ds \bigg]
  + o ( \varepsilon ), 
\end{align*}
where the last equality follows from 
\begin{equation*}
\mathcal{V}_{\lambda } ( t,x, {\lambda }_{0} ) = 1 - \mathbb{E}_{t} [ \zeta ] - \theta \mathbb{E}_{t} \Big[ \Big( {\lambda }_{0} - {X}^{ t,x, \bar{\Pi } }_{T} - \frac{\zeta }{\theta } \Big)_{+} \Big] = 0.
\end{equation*}
Therefore, \cref{eq: equilibrium condition for CNEC :eq} gives $\bar{J} ( t,x, \bar{\Pi }^{ t, \varepsilon, \xi } ) - \bar{J} ( t,x, \bar{\Pi } ) \le o ( \varepsilon )$,
and hence $\bar{\Pi }$ is a CNEC.

\subsection{Proof of \texorpdfstring{\Cref{thm: identical trivial NECs}}{Theorem 4.2}}
\label{pf-thm: identical trivial NECs}

Mirroring the proof of \Cref{thm: verification theorem for CNEC}, one can obtain that $\bar{\Pi }$ is a CNEC.
To show that $\bar{\pi }$ is a CNEC, we differentiate the both sides of each equation in \cref{eq: BSPDE :eq} w.r.t. $x$, together with $\bar{\Pi } ( t,x ) \equiv \bar{\pi }_{t}$, to arrive at the following BSPDE:
\begin{equation*}
\left\{ \begin{aligned}
  - d \mathcal{V}_{x} ( t,x, \lambda )
& = \bigg( \frac{1}{2} \mathcal{V}_{xxx} ( t,x, \lambda ) | \bar{\pi }_{t} {\sigma }_{t} |^{2} 
         + {\psi }_{xx} ( t,x, \lambda ) \bar{\pi }_{t} {\sigma }_{t} \\
& \qquad + \mathcal{V}_{xx} ( t,x, \lambda ) ( {r}_{t} x + \bar{\pi }_{t} {\sigma }_{t} {\vartheta }_{t} )
         + \mathcal{V}_{x} ( t,x, \lambda ) {r}_{t} \bigg) dt
  - {\psi }_{x} ( t,x, \lambda ) d {W}_{t}, \\
  \mathcal{V}_{x} ( T,x, \lambda ) 
& = \zeta + \theta \Big( \lambda - x - \frac{\zeta }{\theta } \Big)_{+},
\end{aligned} \right.
\end{equation*}
for $( \mathcal{V}_{x}, {\psi }_{x} ) \in {C}_{\mathbb{F}} ( 0,T; \mathbb{L}^{2} ( \Omega; {C}^{2,1} ( \mathbb{R} \times \mathbb{R}; \mathbb{R} ) ) ) 
                                   \times {C}_{\mathbb{F}} ( 0,T; \mathbb{L}^{2} ( \Omega; {C}^{2,1} ( \mathbb{R} \times \mathbb{R}; \mathbb{R} ) ) )$.
Applying the It\^o-Kunita-Ventzel formula to $\{ \mathcal{V}_{x} ( s, {X}^{ \bar{\pi } }_{s}, \lambda ) \}_{ s \in [ 0,T ] }$, we obtain
\begin{equation*}
\left\{ \begin{aligned}
  - d \mathcal{V}_{x} ( s, {X}^{ \bar{\pi } }_{s}, \lambda )
& = \mathcal{V}_{x} ( s, {X}^{ \bar{\pi } }_{s}, \lambda ) {r}_{s} ds
  - \big( {\psi }_{x} ( s, {X}^{ \bar{\pi } }_{s}, \lambda ) + \mathcal{V}_{xx} ( s, {X}^{ \bar{\pi } }_{s}, \lambda ) \bar{\pi }_{s} {\sigma }_{s} \big) d {W}_{s}, \\
    \mathcal{V}_{x} ( T, {X}^{ \bar{\pi } }_{T}, \lambda ) 
& = \zeta + \theta \Big( \lambda - {X}^{ \bar{\pi } }_{T} - \frac{\zeta }{\theta } \Big)_{+}.
\end{aligned} \right.
\end{equation*}
Thus, $( Y ( \lambda ), \mathcal{Y} ( \lambda ) ) \in {C}_{\mathbb{F}} ( 0,T; \mathbb{L}^{2} ( \Omega; \mathbb{R} ) ) \times {C}_{\mathbb{F}} ( 0,T; \mathbb{L}^{2} ( \Omega; \mathbb{R} ) )$ given by
\begin{equation*}
{Y}_{s} ( \lambda ) = \mathcal{V}_{x} ( s, {X}^{ \bar{\pi } }_{s}, \lambda ), \quad 
\mathcal{Y}_{s} ( \lambda ) = \mathcal{V}_{xx} ( s, {X}^{ \bar{\pi } }_{s}, \lambda ) \bar{\pi }_{s} {\sigma }_{s} + {\psi }_{x} ( s, {X}^{ \bar{\pi } }_{s}, \lambda )
\end{equation*}
solves the BSDE \cref{eq: BSDE :eq}.
Plugging $x = {X}^{ \bar{\pi } }_{t}$ into \cref{eq: equilibrium condition for trivial NEC :eq} yields ${Y}^{t}_{t} {\vartheta }_{t} + \mathcal{Y}^{t}_{t} = 0$, implying that $\bar{\pi }$ is an ONEC.

\subsection{Proof of \texorpdfstring{\Cref{thm: path-independent ONEC}}{Theorem 5.1}}
\label{pf-thm: path-independent ONEC}

Instead of applying Picard's theorem to prove the existence and uniqueness of solution for \cref{eq: beta equation :eq}, we adopt a more detailed and insightful approach that also serves for discontinuous parameters. 
This method not only provides additional information about $\beta_t$, but also demonstrates an alternative way -- the Picard iteration algorithm -- that can be utilized for numerical simulations.

At first, we show the existence of solution for \cref{eq: beta equation :eq}, or equivalently, for
\begin{equation}
\int_{ - \infty }^{ Q ( \frac{1}{ {\beta }_{t} } ) } \Phi (x) dx = \frac{ 1 - \zeta }{ ( \int_{t}^{T} | {\beta }_{s} {\vartheta }_{s} |^{2} ds )^{ \frac{1}{2} } }.
\label{eq: beta equation for proof :eq}
\end{equation}
where $Q (p) := \inf \{ x: \Phi (x) \ge p \}$ is the quantile function of standard normal distribution.
Let us investigate the following Picard iteration:
\begin{equation*}
\frac{ 1 - \zeta }{ ( \int_{t}^{T} | {\beta }^{(n)}_{s} {\vartheta }_{s} |^{2} ds )^{ \frac{1}{2} } } = \int_{ - \infty }^{ Q ( \frac{1}{ {\beta }^{(n+1)}_{t} } ) } \Phi (x) dx,
\quad {\beta }^{(0)} \equiv 1,
\end{equation*}
which implies that each of $\{ {\beta }^{(n)} \}_{ n \in \mathbb{N}_{+} }$ is continuous and strictly decreasing with ${\beta }^{(n)}_{T-} = 1$.
Notably, it is more convenient for numerical computations to use the following identity: 
\begin{equation}
\label{eq: identity for normal distribution :eq}
  \int_{ - \infty }^{q} \Phi (x) dx
= \int_{ - \infty }^{q} dx \int_{ - \infty }^{x} \Phi' (z) dz
= \int_{ - \infty }^{q} (q-z) \Phi' (z) dz 
= q \Phi (q) + \Phi' (q), \quad \forall q \in \mathbb{R}.
\end{equation}
Denote ${\Theta }_{t} := ( 1 - \zeta ) ( \int_{t}^{T} | {\vartheta }_{s} |^{2} ds )^{ - \frac{1}{2} }$, 
so that $( 1 - \zeta ) ( \int_{t}^{T} | {\beta }^{(n)}_{s} {\vartheta }_{s} |^{2} ds )^{ - \frac{1}{2} } \in [ \frac{ {\Theta }_{t} }{ {\beta }^{(n)}_{t} }, {\Theta }_{t} ]$.
For the auxiliary function ${g}_{t} (z) := \int_{ - \infty }^{ Q(z) } \Phi (x) dx - {\Theta }_{t} z$,
since ${g}_{t}' ( \Phi (x) ) = \frac{ \Phi (x) }{ \Phi'(x) } - {\Theta }_{t}$ is strictly increasing with a unique zero, together with ${g}_{t} (0+) = 0$ and ${g}_{t} (1-) = + \infty $,
we conclude that ${g}_{t} (z)$ has a unique zero (denoted by ${m}_{t}$), which is also the unique zero of $Q(z) - \Psi ( {\Theta }_{t} z )$ and is increasing in $t$.
As a result, ${g}_{t} (x) \ge 0$ if and only if $x \ge {m}_{t}$.
On the one hand, if $\frac{1}{ {\beta }^{(n)}_{t} } \ge {m}_{t}$, then the Picard iteration gives
\begin{equation*}
  \int_{ - \infty }^{ Q ( \frac{1}{ {\beta }^{(n+1)}_{t} } ) } \Phi (x) dx 
\ge \frac{ 1 - \zeta }{ ( \int_{t}^{T} | {\beta }^{(n)}_{t} {\vartheta }_{s} |^{2} ds )^{ \frac{1}{2} } } 
\ge \frac{ {\Theta }_{t} }{ {\beta }^{(n)}_{t} }
\ge {\Theta }_{t} {m}_{t}
= \int_{ - \infty }^{ Q ( {m}_{t} ) } \Phi (x) dx,
\end{equation*}
and hence $\frac{1}{ {\beta }^{(n+1)}_{t} } \ge {m}_{t}$.
As $\frac{1}{ {\beta }^{(0)}_{t} } = 1 > {m}_{t}$, by induction one obtains $\frac{1}{ {\beta }^{(n)}_{t} } \ge {m}_{t}$ for every $n$.
On the other hand, if ${\beta }^{(n)}_{t} \ge {\beta }^{(n-1)}_{t}$ for all $t \in [ 0,T ]$, then the Picard iteration gives
\begin{equation*}
  \int_{ Q ( \frac{1}{ {\beta }^{(n)}_{t} } ) }^{ Q ( \frac{1}{ {\beta }^{(n+1)}_{t} } ) } \Phi (x) dx 
= \frac{ 1 - \zeta }{ ( \int_{t}^{T} | {\beta }^{(n)}_{s} {\vartheta }_{s} |^{2} ds )^{ \frac{1}{2} } } 
- \frac{ 1 - \zeta }{ ( \int_{t}^{T} | {\beta }^{(n-1)}_{s} {\vartheta }_{s} |^{2} ds )^{ \frac{1}{2} } }
\le 0,
\end{equation*}
and hence ${\beta }^{(n+1)}_{t} \ge {\beta }^{(n)}_{t}$ for all $t \in [ 0,T ]$.
As ${\beta }^{(1)}_{t} \ge 1 = {\beta }^{(0)}_{t}$, by induction we obtain $\{ {\beta }^{(n)}_{t} \}_{n}$ is an increasing sequence.
Due to monotone convergence theorem, ${\beta }^{(\infty )}_{t}$ must exist.
Furthermore, from the Picard iteration, one obtains
\begin{align*}
0 & = \frac{ 1 - \zeta }{ ( \int_{ t + \varepsilon }^{T} | {\beta }^{(n)}_{s} {\vartheta }_{s} |^{2} ds )^{ \frac{1}{2} } } - \frac{ 1 - \zeta }{ ( \int_{t}^{T} | {\beta }^{(n)}_{s} {\vartheta }_{s} |^{2} ds )^{ \frac{1}{2} } }
    - \int_{ Q ( \frac{1}{ {\beta }^{(n+1)}_{t} } ) }^{ Q ( \frac{1}{ {\beta }^{(n+1)}_{ t + \varepsilon } } ) } \Phi (x) dx \\
  & = \frac{ ( 1 - \zeta ) \int_{t}^{ t + \varepsilon } | {\beta }^{(n)}_{s} {\vartheta }_{s} |^{2} ds }{ 2 ( \int_{t}^{T} | {\beta }^{(n)}_{s} {\vartheta }_{s} |^{2} ds )^{ \frac{3}{2} } }
    - \frac{ {\beta }^{(n+1)}_{t} - {\beta }^{(n+1)}_{ t + \varepsilon } }{ | {\beta }^{(n+1)}_{t} |^{3} } Q' \Big( \frac{1}{ {\beta }^{(n+1)}_{t} } \Big)
    + o ( \varepsilon ),
\end{align*}
and hence
\begin{align*}
  \limsup_{ \varepsilon \to 0 } \frac{ | {\beta }^{(n+1)}_{ t + \varepsilon } - {\beta }^{(n+1)}_{t} | }{\varepsilon }
& \le \frac{ ( 1 - \zeta ) | {\beta }^{(n+1)}_{t} |^{3} }{ 2 ( \int_{t}^{T} | {\beta }^{(n)}_{s} {\vartheta }_{s} |^{2} ds )^{ \frac{3}{2} } } \Phi' \bigg( Q \Big( \frac{1}{ {\beta }^{(n+1)}_{t} } \Big) \bigg)
      \limsup_{ \varepsilon \to 0 } \bigg| \frac{1}{\varepsilon } \int_{t}^{ t + \varepsilon } | {\beta }^{(n)}_{s} {\vartheta }_{s} |^{2} ds \bigg| \\
& \le \frac{ | {\Theta }_{0} |^{3} }{ 2 ( 1 - \zeta )^{2} | {m}_{0} |^{5} } \bigg( \sup_{ [ 0,T ] } | \vartheta |^{2} \bigg) =: M.
\end{align*}
Thus, ${\beta }^{(n)}$ is Lipschitz continuous, and the Lipschitz continuity parameter $M$ is uniform for all $n$.
For an arbitrarily fixed $\varepsilon > 0$ and $\tau \in [ 0,T ]$, 
there exists a sufficiently large ${N}_{ \varepsilon, \tau } \in \mathbb{N}_{+}$ such that $| {\beta }^{(m)}_{\tau } - {\beta }^{(n)}_{\tau } | < \frac{\varepsilon }{2}$ for all $m > n > {N}_{ \varepsilon, \tau }$.
Consequently, if $| t - \tau | < \frac{\varepsilon }{4M}$, then
\begin{equation*}
    | {\beta }^{(m)}_{t} - {\beta }^{(n)}_{t} | 
\le | {\beta }^{(m)}_{t} - {\beta }^{(m)}_{\tau } | + | {\beta }^{(n)}_{t} - {\beta }^{(m)}_{\tau } | + | {\beta }^{(m)}_{\tau } - {\beta }^{(n)}_{\tau } |
\le 2M | t - \tau | + \frac{\varepsilon }{2} < \varepsilon.
\end{equation*}
By Heine-Borel theorem, there exists a finite set $\{ {\tau }_{k} \}_{k=1,2,\ldots,K}$ such that $\cup_{k=1}^{K} ( {\tau }_{k} - \frac{\varepsilon }{4M}, {\tau }_{k} + \frac{\varepsilon }{4M} ) \supseteq [ 0,T ]$.
Therefore, $| {\beta }^{(m)}_{t} - {\beta }^{(n)}_{t} | < \varepsilon $ for $m > n \ge \max_{k} {N}_{ \varepsilon, {\tau }_{k} }$ and $t \in [ 0,T ]$.
This implies that ${\beta }^{(n)} \rightrightarrows {\beta }^{(\infty )}$ as $n \to \infty $, so that ${\beta }^{(\infty )}$ is continuous and fulfills \cref{eq: beta equation :eq}, and hence the existence of $\beta $ follows.

Now we show the uniqueness of solution for \cref{eq: beta equation :eq}. 
Suppose that $\beta $ and ${\beta }^{*}$ are two solutions of \cref{eq: beta equation :eq}.
By Weierstrass Theorem, $\max_{ t \in [ 0,T ] } {\beta }_{t} \vee {\beta }^{*}_{t} \le {M}_{*}$ for some ${M}_{*} > 1$. 
Consequently, for any $t \in [ 0,T )$ with some sufficiently large constant $K$, 
\begin{align*}
      | {\beta }_{t} - {\beta }^{*}_{t} |
&   = {\beta }_{t} {\beta }^{*}_{t}
      \bigg| \int_{ ( 1 - \zeta ) ( \int_{t}^{T} | {\beta }_{s} {\vartheta }_{s} |^{2} ds )^{ - \frac{1}{2} } }^{ ( 1 - \zeta ) ( \int_{t}^{T} | {\beta }^{*}_{s} {\vartheta }_{s} |^{2} ds )^{ - \frac{1}{2} } }
             \Phi' \big( \Psi (z) \big) \Psi' (z) dz \bigg| \\
& \le {M}_{*}^{2} \frac{ \Phi' \big( \Psi ( {M}_{*}^{-1} {\Theta}_{t} ) \big) }{ \Phi \big( \Psi ( {M}_{*}^{-1} {\Theta}_{t} ) \big) }
      \frac{ ( 1 - \zeta ) \int_{t}^{T} | {\vartheta }_{s} |^{2} ( {\beta }_{s} + {\beta }^{*}_{s} ) | {\beta }_{s} - {\beta }^{*}_{s} | ds }
           { ( \int_{t}^{T} | {\beta }_{s} {\vartheta }_{s} |^{2} ds )^{ \frac{1}{2} } + ( \int_{t}^{T} | {\beta }^{*}_{s} {\vartheta }_{s} |^{2} ds )^{ \frac{1}{2} } } \\
& \le {M}_{*}^{3} \frac{ \Phi' \big( \Psi ( {M}_{*}^{-1} {\Theta}_{t} ) \big) }{ \Phi \big( \Psi ( {M}_{*}^{-1} {\Theta}_{t} ) \big) } {\Theta }_{t} 
      \bigg( \max_{ t \in [ 0,T ] } | {\vartheta }_{t} |^{2} \bigg) \int_{t}^{T} | {\beta }_{s} - {\beta }^{*}_{s} | ds \\
& \le K \int_{t}^{T} | {\beta }_{s} - {\beta }^{*}_{s} | ds,
\end{align*}
where the last inequality is due to the uniform boundedness of $\Phi'(x) \int_{-\infty }^{x} \Phi (z) dz$ on $x \in \mathbb{R}$.
Therefore, by Gr\"onwall's inequality, $| {\beta }_{t} - {\beta }^{*}_{t} | \le 0$, and hence $| {\beta }_{t} - {\beta }^{*}_{t} | = 0$.

To show that ${\beta }_{t}$ is strictly increasing in $\zeta $ for every $t$, 
we differentiate the both sides of \cref{eq: beta equation for proof :eq} w.r.t. $\zeta $, with a slight abuse of notation and rearrangement, to arrive at
\begin{equation*}
  \frac{\partial {\beta }_{t} }{\partial \zeta }
= \frac{ {y}_{t} | {\beta }_{t} |^{3} }{ 1 - \zeta } \Phi' \bigg( Q \Big( \frac{1}{ {\beta }_{t} } \Big) \bigg)
  \bigg( 1 + \frac{ | {y}_{t} |^{2} }{ 1 - \zeta } \int_{t}^{T} | {\vartheta }_{s} |^{2} {\beta }_{s} \frac{\partial {\beta }_{s} }{\partial \zeta } ds \bigg),
\end{equation*}
where ${y}_{t} := ( 1 - \zeta ) ( \int_{t}^{T} | {\beta }_{s} {\vartheta }_{s} |^{2} ds )^{-\frac{1}{2}}$ for ease of expression.
By inverse Gr\"onwall inequality, 
\begin{equation*}
\frac{\partial {\beta }_{t} }{\partial \zeta } \ge \frac{ {y}_{t} | {\beta }_{t} |^{3} }{ 1 - \zeta } \Phi' \bigg( Q \Big( \frac{1}{ {\beta }_{t} } \Big) \bigg) > 0, \quad \forall t \in [ 0,T ).
\end{equation*}

Now we show that \cref{eq: ONEC :eq} gives an ONEC.
As $( r, \sigma, \vartheta, \zeta )$ are deterministic, the BSDE \cref{eq: BSDE :eq} with the identity \cref{eq: lambda functional :eq} gives
\begin{equation*}
  {e}^{ - \int_{s}^{T} {r}_{v} dv } {Y}^{t}_{s} 
= \zeta + \theta \mathbb{E}_{s} \Big[ \Big( \lambda ( \mathbb{P}^{t}_{ {X}^{ \bar{\pi } }_{T}, \zeta } ) - {X}^{ \bar{\pi } }_{T} - \frac{\zeta }{\theta } \Big)_{+} \Big]
= 1 + \int_{t}^{s} {e}^{ - \int_{u}^{T} {r}_{v} dv } \mathcal{Y}^{t}_{u} d {W}_{u}.
\end{equation*}
Corresponding to $\bar{\pi }$ given by \cref{eq: ONEC :eq}, 
\begin{equation*}
{X}^{ \bar{\pi } }_{T} = {X}^{ \bar{\pi } }_{t} {e}^{ \int_{t}^{T} {r}_{v} dv } + \frac{1}{\theta } \int_{t}^{T} {\beta }_{s} | {\vartheta }_{s} |^{2} ds + \frac{1}{\theta } \int_{t}^{T} {\beta }_{s} {\vartheta }_{s} d {W}_{s}.
\end{equation*}
Thus, conditioned on given ${X}^{ \bar{\pi } }_{t}$, 
${X}^{ \bar{\pi } }_{T}$ is normally distributed with mean ${\mu }_{t} := {X}^{ \bar{\pi } }_{t} \exp ( \int_{t}^{T} {r}_{v} dv ) + \frac{1}{\theta } \int_{t}^{T} {\beta }_{s} | {\vartheta }_{s} |^{2} ds$
and standard deviation ${S}_{t} := \frac{1}{\theta } ( \int_{t}^{T} | {\beta }_{s} {\vartheta }_{s} |^{2} ds )^{ \frac{1}{2} } > 0$.
Consequently, by straightforward calculation,
\begin{equation*}
  {e}^{ - \int_{s}^{T} {r}_{v} dv } {Y}_{s} ( \lambda )
= \zeta + \theta \mathbb{E}_{s} \Big[ \Big( \lambda - {X}^{ \bar{\pi } }_{T} - \frac{\zeta }{\theta } \Big)_{+} \Big]
= \zeta + \theta {S}_{s} \int_{ - \infty }^{ \frac{1}{ {S}_{s} } ( \lambda - \frac{\zeta }{\theta } - {\mu }_{s} ) } \Phi (z) dz,
\end{equation*}
and then applying the chain rule and It\^o's rule yields
\begin{equation*}
{e}^{ - \int_{s}^{T} {r}_{v} dv } \mathcal{Y}_{s} ( \lambda ) 
= - \Phi \Big( \frac{ \lambda - \frac{\zeta }{\theta } - {\mu }_{s} }{ {S}_{s} } \Big) \theta \bar{\pi }_{s} {\sigma }_{s} {e}^{ \int_{s}^{T} {r}_{v} dv }
= - \Phi \Big( \frac{ \lambda - \frac{\zeta }{\theta } - {\mu }_{s} }{ {S}_{s} } \Big) {\beta }_{s} {\vartheta }_{s}.
\end{equation*}
In particular, $\lambda ( \mathbb{P}^{t}_{ {X}^{ \bar{\pi } }_{T}, \zeta } ) = \frac{\zeta }{\theta } + {\mu }_{t} + \Psi ( \frac{ 1 - \zeta }{ \theta {S}_{t} } ) {S}_{t}$ follows from
$1 = \zeta + \theta \mathbb{E}_{t} [ ( \lambda ( \mathbb{P}^{t}_{ {X}^{ \bar{\pi } }_{T}, \zeta } ) - {X}^{ \bar{\pi } }_{T} - \frac{\zeta }{\theta } )_{+} ]$.
Hence, the deterministic continuous function \cref{eq: ONEC :eq} gives an ONEC according to \Cref{thm: verification theorem for ONEC} with
\begin{equation*}
  {e}^{ - \int_{t}^{T} {r}_{v} dv } ( {Y}^{t}_{t} {\vartheta }_{t} + \mathcal{Y}^{t}_{t} ) 
= {\vartheta }_{t} - \Phi \circ \Psi \Big( \frac{ 1 - \zeta }{ \theta {S}_{t} } \Big) {\beta }_{t} {\vartheta }_{t}
= 0.
\end{equation*}

Finally, we show the desired uniqueness of ONEC with a slight abuse of notation.
Suppose that $\bar{\pi }$ is another ONEC satisfying the given conditions.
Since $\rho ( \cdot, 0 )$ is deterministic, we have $\varrho ( \cdot, 0 ) \equiv 0$.
Consequently, \cref{eq: BSPDE of rho :eq} reduced to the following linear BSPDE on $[ 0,T ] \times \Omega \times \mathbb{R}$:
\begin{equation*}
\left\{ \begin{aligned}               
- d \rho ( t,y ) & = \frac{1}{2} {\rho }_{yy} ( t,y ) \Big| \frac{ {\alpha }_{t} {\vartheta }_{t} }{\theta } \Big|^{2} dt
                   + {\rho }_{y} ( t,y ) \frac{ {\Gamma }_{t} }{ 2 {H}_{t} } \Big| \frac{ {\alpha }_{t} {\vartheta }_{t} }{\theta } \Big|^{2} dt
                   - {\varrho }_{y} ( t,y ) \frac{ {\alpha }_{t} {\vartheta }_{t} }{\theta } dt
                   - \varrho ( t,y ) d {W}_{t}, \\
    \rho ( T,y ) & = {1}_{\{ y + \frac{ 1 - \zeta }{\theta } \ge 0 \}},
\end{aligned} \right.
\end{equation*}
where ${\alpha }_{t} = \frac{1}{ \rho ( t,0 ) }$, ${H}_{t} = \frac{1}{\theta } ( \int_{t}^{T} | {\alpha }_{s} {\vartheta }_{s} |^{2} ds )^{ \frac{1}{2} }$ and ${\Gamma }_{t} = - {H}_{t} \frac{ {\rho }_{y} ( t,0 ) }{ \rho ( t,0 ) }$.
Using the stochastic Feynman-Kac formula, one can conclude that the above linear BSPDE admits the unique solution
\begin{equation*}
\rho ( t,y ) = \mathbb{E}_{t} \bigg[ {1}_{\{ y + \frac{ 1 - \zeta }{\theta } - \int_{t}^{T} {\Gamma }_{s} d {H}_{s} - \frac{1}{\theta } \int_{t}^{T} {\alpha }_{s} {\vartheta }_{s} d {W}_{s} \ge 0 \}} \bigg] 
             = \Phi \bigg( \frac{ y + \frac{ 1 - \zeta }{\theta } - \int_{t}^{T} {\Gamma }_{s} d {H}_{s} }{ {H}_{t} } \bigg).
\end{equation*}
Let
\begin{equation*}
{G}_{t} = \frac{ 1 - \zeta }{\theta } \bigg( \int_{ - \infty }^{ \frac{1}{ {H}_{t} } ( \frac{ 1 - \zeta }{\theta } - \int_{t}^{T} {\Gamma }_{s} d {H}_{s} ) } \Phi (z) dz \bigg)^{-1}, \quad i.e. \quad
{H}_{t} \Psi \Big( \frac{ 1 - \zeta }{ \theta {G}_{t} } \Big) = \frac{ 1 - \zeta }{\theta } - \int_{t}^{T} {\Gamma }_{s} d {H}_{s}.
\end{equation*}
Then,
\begin{equation*}
{\Gamma }_{t} = - {H}_{t} \frac{ {\rho }_{y} ( t,0 ) }{ \rho ( t,0 ) } 
              = - \frac{ \Phi' \circ \Psi ( \frac{ 1 - \zeta }{ \theta {G}_{t} } ) }{ \Phi \circ \Psi ( \frac{ 1 - \zeta }{ \theta {G}_{t} } ) }
              = \Psi \Big( \frac{ 1 - \zeta }{ \theta {G}_{t} } \Big) - \frac{ 1 - \zeta }{ \theta {G}_{t} } \Psi' \Big( \frac{ 1 - \zeta }{ \theta {G}_{t} } \Big),
\end{equation*}
where the last equality is due to the identity \cref{eq: identity for normal distribution :eq} with $q = \Psi ( \frac{ 1 - \zeta }{ \theta {G}_{t} } )$,
and hence differentiating the both sides of ${H}_{t} \Psi ( \frac{ 1 - \zeta }{ \theta {G}_{t} } ) = \frac{ 1 - \zeta }{\theta } - \int_{t}^{T} {\Gamma }_{s} d {H}_{s}$ w.r.t. $t$ yields
\begin{equation*}
  \Psi \Big( \frac{ 1 - \zeta }{ \theta {G}_{t} } \Big) d {H}_{t} - \frac{ 1 - \zeta }{\theta } \frac{ {H}_{t} }{ | {G}_{t} |^{2} } \Psi' \Big( \frac{ 1 - \zeta }{ \theta {G}_{t} } \Big) d {G}_{t}
= \bigg( \Psi \Big( \frac{ 1 - \zeta }{ \theta {G}_{t} } \Big) - \frac{ 1 - \zeta }{ \theta {G}_{t} } \Psi' \Big( \frac{ 1 - \zeta }{ \theta {G}_{t} } \Big) \bigg) d {H}_{t},
\end{equation*}
implying that $d \frac{ {H}_{t} }{ {G}_{t} } = 0$, i.e., $\frac{G}{H}$ is a constant.
Sending $t$ on the both sides of $1 - \frac{ \theta }{ 1 - \zeta } \int_{t}^{T} {\Gamma }_{s} d {H}_{s} = \frac{ {H}_{t} }{ {G}_{t} } \cdot \frac{ \theta {G}_{t} }{ 1 - \zeta } \Psi \Big( \frac{ 1 - \zeta }{ \theta {G}_{t} } \Big)$ to $T$
delivers $\lim_{t \uparrow T} \frac{ {H}_{t} }{ {G}_{t} } = 1$, and hence $H \equiv G$.
Therefore,
\begin{equation*}
1 = {\alpha }_{t} \rho ( t,0 )
              = {\alpha }_{t} \Phi \bigg( \frac{ \frac{ 1 - \zeta }{\theta } - \int_{t}^{T} d ( {H}_{s} \Psi ( \frac{ 1 - \zeta }{ \theta {H}_{s} } ) ) }{ {H}_{t} } \bigg)
              = {\alpha }_{t} \Phi \circ \Psi \Big( \frac{ 1 - \zeta }{ \theta {H}_{s} } \Big).
\end{equation*}
Given ${\alpha }_{T} = \frac{1}{ \rho ( T,0 ) } = 1$ and the uniqueness of solution for \cref{eq: beta equation :eq}, one concludes that $( \alpha, H ) \equiv ( \beta, S )$.
Then, from \Cref{pf-lem: property of ONEC}, one obtains
\begin{equation*}
d {\lambda }_{t} = \frac{ {\vartheta }_{t} }{ \theta \rho ( t,0 ) } d {W}_{t} 
                 + \frac{ {\rho }_{y} ( t,0 ) }{ 2 \rho ( t,0 ) } \Big| \frac{ {\vartheta }_{t} }{ \theta \rho ( t,0 ) } \Big|^{2} dt
                 = \frac{ {\beta }_{t} {\vartheta }_{t} }{\theta } d {W}_{t} 
                 - \frac{ {\Gamma }_{t} }{ 2 {S}_{t} } \Big| \frac{ {\beta }_{t} {\vartheta }_{t} }{ \theta } \Big|^{2} dt
                 = \frac{ {\beta }_{t} {\vartheta }_{t} }{\theta } d {W}_{t} + d \bigg( {S}_{t} \Psi \Big( \frac{ 1 - \zeta }{ \theta {S}_{t} } \Big) \bigg).
\end{equation*}
Hence, according to \Cref{lem: property of ONEC},
\begin{equation*}
{X}^{ \bar{\pi } }_{T} = {\lambda }_{T} - \frac{1}{\theta } 
                       = {\lambda }_{0} - \frac{\zeta }{\theta } - {S}_{0} \Psi \Big( \frac{ 1 - \zeta }{ \theta {S}_{0} } \Big) + \frac{1}{\theta } \int_{0}^{T} {\beta }_{t} {\vartheta }_{t} d {W}_{t},
\end{equation*}
which corresponds to the unique hedging portfolio $\bar{\pi }$ given by \cref{eq: ONEC :eq} and \cref{eq: beta equation :eq}. 
So we complete this proof.

\subsection{Proof of \texorpdfstring{\Cref{thm: state-independent CNEC}}{Theorem 5.3}}
\label{pf-thm: state-independent CNEC}

We have shown in \Cref{pf-thm: path-independent ONEC} that corresponding to \cref{eq: CNEC :eq},
\begin{equation*}
{X}^{ t,x, \bar{\Pi } }_{T} = x {e}^{ \int_{t}^{T} {r}_{v} dv } + \frac{1}{\theta } \int_{t}^{T} {\beta }_{s} | {\vartheta }_{s} |^{2} ds + \frac{1}{\theta } \int_{t}^{T} {\beta }_{s} {\vartheta }_{s} d {W}_{s}
\end{equation*}
is normally distributed with mean ${\mu }_{t} = x \exp ( \int_{t}^{T} {r}_{v} dv ) + \frac{1}{\theta } \int_{t}^{T} {\beta }_{s} | {\vartheta }_{s} |^{2} ds$
and standard deviation ${S}_{t} = \frac{1}{\theta } ( \int_{t}^{T} | {\beta }_{s} {\vartheta }_{s} |^{2} ds )^{\frac{1}{2}}$.
According to \Cref{pf-thm: verification theorem for CNEC}, 
\begin{align*}
    \Lambda ( t,x ) 
& = \lambda ( \mathbb{P}^{t}_{ {X}^{ t,x, \bar{\Pi } }_{T}, \zeta } ) 
  = \frac{\zeta }{\theta } + {\mu }_{t} + \Psi \Big( \frac{ 1 - \zeta }{ \theta {S}_{t} } \Big) {S}_{t}, \\
    \mathcal{V} ( t,x, \lambda )
& = \lambda ( 1 - \zeta )
  - \frac{\theta }{2} \mathbb{E}_{t} \bigg[ \Big| \Big( \lambda - {X}^{ t,x, \bar{\Pi } }_{T} - \frac{\zeta }{\theta } \Big)_{+} \Big|^{2} \bigg]
  + \zeta \mathbb{E}_{t} [ {X}^{ t,x, \bar{\Pi } }_{T} ] + \frac{1}{ 2 \theta } {\zeta }^{2} - \frac{1}{ 2 \theta }, \\
    \mathcal{M} ( t,x, \lambda ) 
& = \zeta + \theta \mathbb{E}_{t} \Big[ \Big( \lambda - {X}^{ t,x, \bar{\Pi } }_{T} - \frac{\zeta }{\theta } \Big)_{+} \Big]
  = \zeta + \theta {S}_{t} \int_{ - \infty }^{ \frac{1}{ {S}_{t} } ( \lambda - \frac{\zeta }{\theta } - {\mu }_{t} ) } \Phi (y) dy,
\end{align*}
and hence $\psi \equiv 0$ and $\phi \equiv 0$.
Obviously, $\mathcal{M} ( t, \cdot, \lambda )$ is twice continuously differentiable.
By the equilibrium condition \cref{eq: equilibrium condition for CNEC :eq}, to verify that \cref{eq: CNEC :eq} provides a CNEC, it suffices to show 
\begin{equation}
  \mathcal{V}_{xx} \big( t,x, \lambda ( \mathbb{P}^{t}_{ {X}^{ t,x, \bar{\Pi } }_{T}, \zeta } ) \big) \frac{ {\beta }_{t} }{\theta } {e}^{ \int_{t}^{T} {r}_{v} dv } 
+ \mathcal{V}_{x} \big( t,x, \lambda ( \mathbb{P}^{t}_{ {X}^{ t,x, \bar{\Pi } }_{T}, \zeta } ) \big) = 0
\label{eq: equilibrium condition for CNEC with deterministic parameters :eq}
\end{equation}
and $\mathcal{V}_{xx} < 0$.
In fact, by straightforward calculation, 
\begin{align*}
  \mathcal{V}_{x} ( t,x, \lambda ) 
& = {e}^{ \int_{t}^{T} {r}_{v} dv } \bigg( \theta \mathbb{E}_{t} \Big[ \Big( \lambda - {X}^{ t,x, \bar{\Pi } }_{T} - \frac{\zeta }{\theta } \Big)_{+} \Big] + \zeta \bigg), \\
  \mathcal{V}_{xx} ( t,x, \lambda ) 
& = - {e}^{ 2 \int_{t}^{T} {r}_{v} dv } \theta \mathbb{P} \Big( \lambda - {X}^{ t,x, \bar{\Pi } }_{T} - \frac{\zeta }{\theta } \ge 0 \Big)
  = - {e}^{ 2 \int_{t}^{T} {r}_{v} dv } \theta \Phi \bigg( \frac{ \lambda - \frac{\zeta }{\theta } - {\mu }_{t} }{ {S}_{t} } \bigg) < 0.
\end{align*}
Plugging $\lambda = \lambda ( \mathbb{P}^{t}_{ {X}^{ t,x, \bar{\Pi } }_{T}, \zeta } ) = \frac{\zeta }{\theta } + {\mu }_{t} + \Psi ( \frac{ 1 - \zeta }{ \theta {S}_{t} } ) {S}_{t}$ into the above results yields
\begin{equation*}
  \mathcal{V}_{x} \big( t,x, \lambda ( \mathbb{P}^{t}_{ {X}^{ t,x, \bar{\Pi } }_{T}, \zeta } ) \big) 
= {e}^{ \int_{t}^{T} {r}_{v} dv }, \quad 
  \mathcal{V}_{xx} \big( t,x, \lambda ( \mathbb{P}^{t}_{ {X}^{ t,x, \bar{\Pi } }_{T}, \zeta } ) \big) 
= - {e}^{ 2 \int_{t}^{T} {r}_{v} dv } \theta \Phi \circ \Psi \Big( \frac{ 1 - \zeta }{ \theta {S}_{t} } \Big).
\end{equation*}
As $\beta $ satisfies \cref{eq: beta equation :eq}, \cref{eq: equilibrium condition for CNEC with deterministic parameters :eq} holds,
and hence $\bar{\Pi }$ given by \cref{eq: CNEC :eq} and \cref{eq: beta equation :eq} is a CENC.

Furthermore, by straightforward calculation with a slight abuse of notation,
for any $\Pi $ given by \cref{eq: CNEC :eq} with any fixed piecewise continuous function 
$\beta $ and $\lambda ( \mathbb{P}^{t}_{ {X}^{ t,x, \Pi }_{T}, \zeta } ) = \frac{\zeta }{\theta } + {\mu }_{t} + \Psi ( \frac{ 1 - \zeta }{ \theta {S}_{t} } ) {S}_{t}$, we have
\begin{align*}
    \bar{J} ( t,x, \Pi )
& = \mathcal{V} \big( t,x, \lambda ( \mathbb{P}^{t}_{ {X}^{ t,x, \Pi }_{T}, \zeta } ) \big) \\
& = {\mu }_{t} + ( 1 - \zeta ) \Psi \Big( \frac{ 1 - \zeta }{ \theta {S}_{t} } \Big) {S}_{t} \\
& \quad - \frac{\theta }{2} | {S}_{t} |^{2} \mathbb{E} \bigg[ \Big| \Big( \Psi \Big( \frac{ 1 - \zeta }{ \theta {S}_{t} } \Big) - \frac{1}{\theta {S}_{t} } \int_{t}^{T} {\beta }_{s} {\vartheta }_{s} d {W}_{s} \Big)_{+} \Big|^{2} \bigg]
  - \frac{ ( 1 - \zeta )^{2} }{ 2 \theta } \\
& = {\mu }_{t} + ( 1 - \zeta ) \Psi \Big( \frac{ 1 - \zeta }{ \theta {S}_{t} } \Big) {S}_{t} \\
& \quad - \theta | {S}_{t} |^{2} \int_{ - \infty }^{ \Psi ( \frac{ 1 - \zeta }{ \theta {S}_{t} } ) } d \Phi (x) \int_{x}^{ \Psi ( \frac{ 1 - \zeta }{ \theta {S}_{t} } ) } \bigg( \Psi \Big( \frac{ 1 - \zeta }{ \theta {S}_{t} } \Big) - y \bigg) dy
        - \frac{ ( 1 - \zeta )^{2} }{ 2 \theta } \\
& = {\mu }_{t} + ( 1 - \zeta ) \Psi \Big( \frac{ 1 - \zeta }{ \theta {S}_{t} } \Big) {S}_{t}
  - \theta | {S}_{t} |^{2} \int_{ - \infty }^{ \Psi ( \frac{ 1 - \zeta }{ \theta {S}_{t} } ) } \bigg( \Psi \Big( \frac{ 1 - \zeta }{ \theta {S}_{t} } \Big) - y \bigg) \Phi (y) dy
  - \frac{ ( 1 - \zeta )^{2} }{ 2 \theta } \\
& = {\mu }_{t} + \theta | {S}_{t} |^{2} \int_{ - \infty }^{ \Psi ( \frac{ 1 - \zeta }{ \theta {S}_{t} } ) } y \Phi (y) dy - \frac{ ( 1 - \zeta )^{2} }{ 2 \theta } \\
& = {\mu }_{t} + \theta | {S}_{t} |^{2} \int_{0}^{ \frac{ 1 - \zeta }{ \theta {S}_{t} } } \Psi (z) dz - \frac{ ( 1 - \zeta )^{2} }{ 2 \theta },
\end{align*}
where the last line arises from $dz = \Phi \circ \Psi (z) d \Psi (z)$.
As a result, the second claim of \Cref{thm: state-independent CNEC} follows.

\subsection{Proof of \texorpdfstring{\Cref{thm: equilibrium value function and strong equilibria}}{Proposition 5.5}}
\label{pf-thm: equilibrium value function and strong equilibria}

For the first claim, it suffices to show $\bar{J} ( t,x, \bar{\Pi } ) > \bar{J} ( t,x, \bar{\Pi }^{ t, \varepsilon, \bar{\Pi } ( t,x ) } )$ for any sufficiently small $\varepsilon $,
since \Cref{pf-thm: verification theorem for CNEC,pf-thm: state-independent CNEC} have provided $\mathcal{V}_{xx} < 0$ and
\begin{align*}
& \bar{J} ( t,x, \bar{\Pi }^{ t, \varepsilon, \xi } ) - \bar{J} ( t,x, \bar{\Pi } ) \\
& = \bigg( \frac{1}{2} \mathcal{V}_{xx} ( t,x, {\lambda }_{0} ) \big( {\xi }^{2} - | \bar{\Pi } ( t,x ) |^{2} \big) | \sigma |^{2} 
         + \mathcal{V}_{x} ( t,x, {\lambda }_{0} ) \vartheta \big( \xi - \bar{\Pi } ( t,x ) \big) \sigma \bigg) \varepsilon 
  + o ( \varepsilon ),
\end{align*}
which by \cref{eq: equilibrium condition for CNEC :eq}, it implies that
\begin{equation*}
\lim_{ \varepsilon \downarrow 0 } \frac{1}{\varepsilon } \big( \bar{J} ( t,x, \bar{\Pi }^{ t, \varepsilon, \xi } ) - \bar{J} ( t,x, \bar{\Pi } ) \big) < 0, \quad \forall \xi \ne \bar{\Pi } ( t,x ).
\end{equation*}
Due to the proof of the second claim of \Cref{thm: state-independent CNEC} and
\begin{equation*}
\bar{\Pi }^{ t, \varepsilon, \bar{\Pi } ( t,x ) } ( s,y ) = \Big( {1}_{\{ s \in [ t, t + \varepsilon ) \}} {\beta }_{t} {e}^{ - r (s-t) } + {1}_{\{ s \notin [ t, t + \varepsilon ) \}} {\beta }_{s} \Big) \frac{\vartheta }{ \theta \sigma } {e}^{ - r (T-s) },
\end{equation*} 
we obtain
\begin{align*}
& \theta | \vartheta |^{-2} \big( \bar{J} ( t,x, \bar{\Pi } ) - \bar{J} ( t,x, \bar{\Pi }^{ t, \varepsilon, \bar{\Pi } ( t,x ) } ) \big) \\
& = \int_{t}^{ t + \varepsilon } {\beta }_{s} ds - {\beta }_{t} \int_{0}^{\varepsilon } {e}^{ - r s } ds
  + \int_{t}^{T} | {\beta }_{s} |^{2} ds \int_{0}^{ \frac{ 1 - \zeta }{\vartheta } ( \int_{t}^{T} | {\beta }_{s} |^{2} ds )^{ - \frac{1}{2} } } \Psi(y) dy \\
& \quad - \bigg( \int_{ t + \varepsilon }^{T} | {\beta }_{s} |^{2} ds + | {\beta }_{t} |^{2} \int_{0}^{\varepsilon } {e}^{ - 2 r s } ds \bigg)
          \int_{0}^{ \frac{ 1 - \zeta }{\vartheta } ( \int_{ t + \varepsilon }^{T} | {\beta }_{s} |^{2} ds + | {\beta }_{t} |^{2} \int_{0}^{\varepsilon } {e}^{ - 2 r s } ds )^{ - \frac{1}{2} } } \Psi(y) dy. 
\end{align*}
By Taylor expansion, 
\begin{align*}
  \int_{t}^{ t + \varepsilon } {\beta }_{s} ds - {\beta }_{t} \int_{0}^{\varepsilon } {e}^{ - r s } ds
& = \frac{1}{2} ( {\beta }_{t}' + r {\beta }_{t} ) {\varepsilon }^{2} 
  + \frac{1}{6} ( {\beta }_{t}'' - {r}^{2} {\beta }_{t} ) {\varepsilon }^{3} 
  + O ( {\varepsilon }^{4} ), \\
  \int_{t}^{ t + \varepsilon } | {\beta }_{s} |^{2} ds - | {\beta }_{t} |^{2} \int_{0}^{\varepsilon } {e}^{ - 2 r s } ds 
& = ( {\beta }_{t}' {\beta }_{t} + r | {\beta }_{t} |^{2} ) {\varepsilon }^{2} 
  + \frac{1}{3} ( {\beta }_{t}'' {\beta }_{t} + | {\beta }_{t}' |^{2} - 2 {r}^{2} | {\beta }_{t} |^{2} ) {\varepsilon }^{3} 
  + O ( {\varepsilon }^{4} ),
\end{align*}
the latter of which implies that
\begin{align*}
& \bigg( \int_{ t + \varepsilon }^{T} | {\beta }_{s} |^{2} ds + | {\beta }_{t} |^{2} \int_{0}^{\varepsilon } {e}^{ - 2 r s } ds \bigg)
         \int_{0}^{ \frac{ 1 - \zeta }{\vartheta } ( \int_{ t + \varepsilon }^{T} | {\beta }_{s} |^{2} ds + | {\beta }_{t} |^{2} \int_{0}^{\varepsilon } {e}^{ - 2 r s } ds )^{ - \frac{1}{2} } } \Psi(y) dy \\
& = \int_{t}^{T} | {\beta }_{s} |^{2} ds \int_{0}^{ \frac{ 1 - \zeta }{\vartheta ( \int_{t}^{T} | {\beta }_{s} |^{2} ds )^{ \frac{1}{2} } } } \Psi(y) dy \\
& \quad - \Big( ( {\beta }_{t}' {\beta }_{t} + r | {\beta }_{t} |^{2} ) {\varepsilon }^{2} 
              + \frac{1}{3} ( {\beta }_{t}'' {\beta }_{t} + | {\beta }_{t}' |^{2} - 2 {r}^{2} | {\beta }_{t} |^{2} ) {\varepsilon }^{3} \Big) \\
& \qquad  \times \Bigg( \int_{0}^{ \frac{ 1 - \zeta }{ \vartheta ( \int_{t}^{T} | {\beta }_{s} |^{2} ds )^{ \frac{1}{2} } } } \Psi(y) dy
                      - \frac{ 1 - \zeta }{ 2 \vartheta ( \int_{t}^{T} | {\beta }_{s} |^{2} ds )^{ \frac{1}{2} } } 
                        \Psi \bigg( \frac{ 1 - \zeta }{ \vartheta ( \int_{t}^{T} | {\beta }_{s} |^{2} ds )^{ \frac{1}{2} } } \bigg) \Bigg)
        + O ( {\varepsilon }^{4} ).
\end{align*}
Hence, with ${B}_{t} := \frac{ 1 - \zeta }{\vartheta } ( \int_{t}^{T} | {\beta }_{s} |^{2} ds )^{ - \frac{1}{2} }$ for the sake of brevity, we obtain
\begin{align}
& \theta | \vartheta |^{-2} \big( \bar{J} ( t,x, \bar{\Pi } ) - \bar{J} ( t,x, \bar{\Pi }^{ t, \varepsilon, \bar{\Pi } ( t,x ) } ) \big) \notag \\
& = \Big( \frac{1}{2} ( {\beta }_{t}' + r {\beta }_{t} ) {\varepsilon }^{2} + \frac{1}{6} ( {\beta }_{t}'' - {r}^{2} {\beta }_{t} ) {\varepsilon }^{3} \Big)
    \bigg( 1 + {\beta }_{t} \Big( 2 \int_{0}^{ {B}_{t} } \Psi(y) dy - {B}_{t} \Psi ( {B}_{t} ) \Big) \bigg) \notag \\
& \quad
  - \frac{1}{6} ( | r {\beta }_{t} |^{2} - | {\beta }_{t}' |^{2} ) 
                \bigg( 2 \int_{0}^{ {B}_{t} } \Psi(y) dy - {B}_{t} \Psi ( {B}_{t} ) \bigg) {\varepsilon }^{3}
  + O ( {\varepsilon }^{4} ).
\label{eq: asymptotic estimation for strong CNEC}
\end{align}
Now we show that $\chi (x) := \Phi (x) + 2 \int_{ - \infty }^{x} y \Phi (y) dy - x \int_{ - \infty }^{x} \Phi (y) dy \equiv 0$,
so that plugging $x = \Psi ( {B}_{t} )$ with \cref{eq: beta equation :eq} and $dy = \Phi \circ \Psi (y) d \Psi (y)$ into $\chi (x)$ 
immediately yields $1 + {\beta }_{t} ( 2 \int_{0}^{ {B}_{t} } \Psi(y) dy - {B}_{t} \Psi ( {B}_{t} ) ) = 0$.
In fact, by $\chi'(x) = \Phi' (x) + x \Phi (x) - \int_{ - \infty }^{x} \Phi (y) dy = \Phi' (x) + \int_{ - \infty }^{x} y \Phi' (y) dy = 0$ and $\chi ( - \infty ) = 0$, we obtain $\chi \equiv 0$.
Consequently, 
\begin{equation*}
  \bar{J} ( t,x, \bar{\Pi } ) - \bar{J} ( t,x, \bar{\Pi }^{ t, \varepsilon, \bar{\Pi } ( t,x ) } ) 
= \frac{ | \vartheta |^{2} }{ 6 \theta {\beta }_{t} } ( r {\beta }_{t} - {\beta }_{t}' ) ( r {\beta }_{t} + {\beta }_{t}' ) {\varepsilon }^{3} + O ( {\varepsilon }^{4} ).
\end{equation*}
As $\beta $ is strictly decreasing with ${\beta }_{T} = 1$ (see \Cref{thm: path-independent ONEC}), under the condition $r > \frac{ | {\beta }_{t}' | }{ {\beta }_{t} }$ we have
\begin{equation*}
\lim_{ \varepsilon \downarrow 0 } \frac{1}{ {\varepsilon }^{3} } \big( \bar{J} ( t,x, \bar{\Pi }^{ t, \varepsilon, \bar{\Pi } ( t,x ) } ) - \bar{J} ( t,x, \bar{\Pi } ) \big) < 0.
\end{equation*}
Hence, the first claim follows.

In terms of the second claim, we proceed with the following general result:
\begin{align*}
& \mathbb{E}_{t} \bigg[ {1}_{A} \bigg| \Big( \int_{t}^{ t + \varepsilon } {\gamma }_{s} {\rho }_{s} d {W}_{s} \Big)^{2} - \Big( \int_{t}^{ t + \varepsilon } {\gamma }_{s} d {W}_{s} \Big)^{2} \bigg| \bigg] \\
& \le \bigg( \int_{t}^{ t + \varepsilon } \mathbb{E}_{t} [ | {\gamma }_{s} ( {\rho }_{s} + 1 ) |^{2} ] ds \bigg)^{ \frac{1}{2} }
      \bigg( \int_{t}^{ t + \varepsilon } \mathbb{E} [ | {\gamma }_{s} ( {\rho }_{s} - 1 ) |^{2} ] ds \bigg)^{ \frac{1}{2} } \\
& \le \varepsilon \bigg( \esssup_{ [ t,T ] \times \Omega } | \gamma | \bigg)^{2} \bigg( 1 + \esssup_{ [ t,T ] \times \Omega } | \rho | \bigg)^{2}
      \bigg( \frac{1}{\varepsilon } \int_{t}^{ t + \varepsilon } \mathbb{E} [ ( {\rho }_{s} - 1 )^{2} ] ds \bigg)^{ \frac{1}{2} } \\
&   = o ( \varepsilon )
\end{align*}
for any $A \in \mathcal{F}_{T}$ and $\gamma, \rho \in {C}_{\mathbb{F}} ( t,T ; \mathbb{L}^{\infty } ( \Omega; \mathcal{R} ) )$ with ${\rho }_{t} = 1$.
Therefore, plugging \cref{eq: equilibrium condition for ONEC :eq} with
\begin{align*}
& \lim_{\varepsilon \downarrow 0} \frac{1}{\varepsilon }
  \mathbb{E}_{t} \bigg[ \bigg( \int_{t}^{ t + \varepsilon } {e}^{ (T-s) r } \sigma d {W}_{s} \bigg)^{2} 
                                                            {1}_{\{ {X}^{ \bar{\pi } }_{T} + \frac{\zeta }{\theta } \le \lambda ( \mathbb{P}^{t}_{ {X}^{ \bar{\pi } }_{T}, \zeta } ) \}} \bigg] \\
& = {e}^{ 2 (T-t) r } {\sigma }^{2}
    \lim_{\varepsilon \downarrow 0} \frac{1}{\varepsilon }
    \mathbb{E}_{t} \bigg[ \Big( \int_{t}^{ t + \varepsilon } \frac{ {\beta }_{s} }{ {\beta }_{t} } d {W}_{s} \Big)^{2} 
                          \mathbb{E}_{ t + \varepsilon } \Big[ {1}_{\big\{ \frac{ \int_{t}^{T} {\beta }_{s} d {W}_{s} }{ ( \int_{t}^{T} | {\beta }_{s} |^{2} ds )^{ \frac{1}{2} } }
                                                                           \le \Psi \big( \frac{ 1 - \zeta }{ \vartheta ( \int_{t}^{T} | {\beta }_{s} |^{2} ds )^{ \frac{1}{2} } } \big) \big\}} \Big] \bigg] \\
& = {e}^{ 2 (T-t) r } \frac{ {\sigma }^{2} }{ | {\beta }_{t} |^{2} }
    \lim_{\varepsilon \downarrow 0} \frac{ \int_{t}^{ t + \varepsilon } | {\beta }_{s} |^{2} ds }{\varepsilon } \\
& \quad \times \lim_{\varepsilon \downarrow 0}
    \int_{ - \infty }^{ + \infty } {z}^{2} \Phi \bigg( \Psi \Big( \frac{ 1 - \zeta }{ \vartheta ( \int_{t}^{T} | {\beta }_{s} |^{2} ds )^{ \frac{1}{2} } } \Big) 
                                                       \frac{ ( \int_{t}^{T} | {\beta }_{s} |^{2} ds )^{ \frac{1}{2} } }{ ( \int_{ t + \varepsilon }^{T} | {\beta }_{s} |^{2} ds )^{ \frac{1}{2} } }
                                                     - z \frac{ ( \int_{t}^{ t + \varepsilon } | {\beta }_{s} |^{2} ds )^{ \frac{1}{2} } }
                                                              { ( \int_{ t + \varepsilon }^{T} | {\beta }_{s} |^{2} ds )^{ \frac{1}{2} } } \bigg) d \Phi (z) \\
& = {e}^{ 2 (T-t) r } {\sigma }^{2} \Phi \circ \Psi \bigg( \frac{ 1 - \zeta }{ \vartheta ( \int_{t}^{T} | {\beta }_{s} |^{2} ds )^{ \frac{1}{2} } } \bigg)
\end{align*} 
into \cref{eq: asymptotic estimate for ONEC :eq} yields 
\begin{equation*}
  \lim_{\varepsilon \downarrow 0} \frac{1}{\varepsilon } \big( J ( t, \bar{\pi }^{ t, \varepsilon, \xi } ) - J ( t, \bar{\pi } ) \big)
= \lim_{ \varepsilon \downarrow 0 } \frac{1}{\varepsilon } \big( g ( \mathbb{P}^{t}_{ \bar{X}^{\varepsilon }_{T}, \zeta } ) - g ( \mathbb{P}^{t}_{ {X}^{ \bar{\pi } }_{T}, \zeta } ) \big) < 0,
\quad \forall \xi \ne 0.
\end{equation*} 
So the proof is completed.

\bibliographystyle{apacite}
\bibliography{SMMV-references}
\end{document}